\pdfoutput=1
\RequirePackage{ifpdf}
\ifpdf 
\documentclass[pdftex]{sigma}
\else
\documentclass{sigma}
\fi

\usepackage{stmaryrd, mathtools}
\usepackage{tikz-cd}
\tikzset{%
    symbol/.style={%
        ,draw=none
        ,every to/.append style={%
            edge node={node [sloped, allow upside down, auto=false]{$#1$}}}
    }
}
\usepackage[makeroom]{cancel}
\usepackage{bbm}

\usepackage{xfrac}
\usepackage[shortlabels]{enumitem}

\usepackage{mathtools}
\usepackage{mathrsfs}
\usepackage{comment}

\usepackage{tikz}
\usetikzlibrary{arrows,decorations.markings,matrix}
\usetikzlibrary{decorations.pathmorphing}

\numberwithin{equation}{section}

\newtheorem{Theorem}{Theorem}[section]
\newtheorem*{Theorem*}{Theorem}
\newtheorem{Corollary}[Theorem]{Corollary}
\newtheorem{Lemma}[Theorem]{Lemma}
\newtheorem{Proposition}[Theorem]{Proposition}
 { \theoremstyle{definition}
\newtheorem{Definition}[Theorem]{Definition}

\newtheorem{Example}[Theorem]{Example}
\newtheorem{Remark}[Theorem]{Remark} }

\DeclareMathOperator{\action}{\curvearrowright}

\newcommand{\R}{\mathbb{R}}

\newcommand{\g}{\mathfrak{g}}

\renewcommand{\epsilon}{\varepsilon}

\newcommand{\Ham}{\mathsf{Ham}}
\newcommand{\ham}{\text{ham}}

\newcommand{\fgmodule}{\mathfrak{X}_{\g}}
\newcommand{\tgtvanform}{I_{\Omega}(N)}
\newcommand{\tgtvf}{\mathfrak X_N(M)}
\newcommand{\vanvf}{I_{\mathfrak{X}}(N)}

\newcommand{\A}{\mathcal{A}}

\newcommand{\G}{\mathcal{G}}
\renewcommand{\L}{\mathscr{L}}
\newcommand{\M}{\mathcal{M}}

\newcommand{\X}{\mathfrak{X}}
\newcommand{\ad}{\mathrm{ad}}
\newcommand{\Ad}{\mathrm{Ad}}

\renewcommand{\subset}{\subseteq}
\renewcommand{\supset}{\supseteq}
\renewcommand{\d}{\mathrm{d}}

\usepackage{tikz-cd}
\tikzset{%
    symbol/.style={%
        ,draw=none
        ,every to/.append style={%
            edge node={node [sloped, allow upside down, auto=false]{$#1$}}}
    }
}
\tikzset{/tikz/commutative diagrams/diagrams={column sep =large}}
\usetikzlibrary{decorations.pathmorphing}
\usetikzlibrary{cd, fit}

\newcommand{\morphism}[5]{%
  \begin{tikzcd}[
    column sep=2em,
    row sep=-.5ex,
    ampersand replacement=\&
  ]
  #1\colon~ \&[-3em]
  #2\vphantom{#3} \arrow[r] \&
  #3\vphantom{#2} \\
  \&
  #4\vphantom{#5} \arrow[r,mapsto] \&
  #5\vphantom{#4}
  \end{tikzcd}%
}

\begin{document}

\allowdisplaybreaks

\newcommand{\arXivNumber}{2206.03137}

\renewcommand{\PaperNumber}{061}

\FirstPageHeading

\ShortArticleName{Reduction of $L_\infty$-Algebras of Observables on Multisymplectic Manifolds}

\ArticleName{Reduction of $\boldsymbol{L_\infty}$-Algebras of Observables\\ on Multisymplectic Manifolds}

\Author{Casey BLACKER~$^{\rm a}$, Antonio Michele MITI~$^{\rm b}$ and Leonid RYVKIN~$^{\rm c}$}

\AuthorNameForHeading{C.~Blacker, A.M.~Miti and L.~Ryvkin}

\Address{$^{\rm a)}$~Department of Mathematical Sciences, George Mason University,\\
\hphantom{$^{\rm a)}$}~4400 University Dr, Fairfax, VA 22030, USA}
\EmailD{\href{mailto:cblacke@gmu.edu}{cblacke@gmu.edu}}
\URLaddressD{\url{https://math.gmu.edu/~cblacke/}}

\Address{$^{\rm b)}$~Dipartimento di Matematica, Sapienza Universit\`a di Roma,\\
\hphantom{$^{\rm b)}$}~Piazzale Aldo Moro 5, 00185 Roma, Italy}
\EmailD{\href{mailto:antoniomichele.miti@uniroma1.it}{antoniomichele.miti@uniroma1.it}}
\URLaddressD{\url{https://www.antoniomiti.it/}}

\Address{$^{\rm c)}$~Institut Camille Jordan, Universit\'e Claude Bernard Lyon 1,\\
\hphantom{$^{\rm c)}$}~43 boulevard du 11 novembre 1918, 69622 Villeurbann, France}
\EmailD{\href{mailto:leonid.ryvkin@math.univ-lyon1.fr}{leonid.ryvkin@math.univ-lyon1.fr}}
\URLaddressD{\url{https://ryvkin.eu/}}

\ArticleDates{Received October 24, 2023, in final form June 24, 2024; Published online July 03, 2024}

\Abstract{We develop a reduction scheme for the $L_\infty$-algebra of observables on a premultisymplectic manifold $(M,\omega)$ in the presence of a compatible Lie algebra action $\mathfrak{g}\curvearrowright M$ and subset $N\subset M$. This reproduces in the symplectic setting the Poisson algebra of observables on the Marsden--Weinstein--Meyer symplectic reduced space, whenever the reduced space exists, but is otherwise distinct from the Dirac, \'Sniatycki--Weinstein, and Arms--Cushman--Gotay observable reduction schemes. We examine various examples, including multicotangent bundles and multiphase spaces, and we conclude with a discussion of applications to classical field theories and quantization.}


\Classification{53D05; 53D20; 70S05; 70S10}


\section{Introduction}

Reduction has its origins in mechanics, specifically in the problem of reducing the degrees of freedom of a mechanical system with symmetries. In this regard, it is rooted in the work of Euler and Lagrange on motions of the rigid body. We refer to \cite{CendraMarsdenRatiu01, MarsdenWeinstein01} for a historical review. The manifestation of these ideas in symplectic geometry is the celebrated Marsden--Weinstein--Meyer theorem \cite{MarsdenWeinstein74,Meyer73}. The underlying ideas of symplectic reduction have been extended throughout differential geometry. Indeed, the range of adaptations includes Poisson manifolds \cite{MarsdenRatiu86,Willett02}, contact structures \cite{deLeonLainzValcazar19,Willett02}, cosymplectic manifolds \cite{Albert89}, polysymplectic manifolds \cite{MarreroRoman-RoySalgadoVilarino15}, higher Poisson structures \cite{BursztynMartinezAlbaRubio19}, Courant algebroids and generalized complex structures \cite{BursztynCavalcantiGualtieri07,StienonXu08,Vaisman07}, and quasi-Hamiltonian $G$-spaces \cite{AlekseevMalkinMeinrenken98}.

It is possible to formally extend the symplectic formalism to continuous physical systems by defining a presymplectic structure on the infinite-dimensional space of solutions of the field equations.
However, this procedure encounters several technicalities, for instance, those arising from introducing suitably weaker notions of smooth structures apt to deal with infinite-dimensional manifolds.
The covariant formalism avoids these issues by working not with a space of solutions, but rather on a finite-dimensional \emph{multiphase space}, a multisymplectic manifold associated to the configuration bundle.
We refer to \cite{GotayIsenbergMarsden04, GotayIsenbergMarsdenMontgomery98,Helein12,Kijowski73} for background on the mathematical physics that precipitated much of the development of the multisymplectic formalism.

The question of multisymplectic reduction is first addressed in \cite{Sniatycki04}, where an extension of the Marsden--Weinstein quotient is defined for multiphase spaces associated to classical field theories. The prospect of a general multisymplectic reduction procedure is given in \cite{MarsdenWeinstein01,OrtegaRatiu04} and a more thorough examination in \cite{Echeverria-Enriquez18}. Building on this work, a reduction for general multisymplectic Hamiltonian $G$-spaces is proposed in \cite{Blacker21}. We note that an interesting alternative perspective, reflecting a mechanical notion of reduction and distinct from the Marsden--Weinstein quotient in the symplectic setting, appears in \cite[Definition 4.2]{deLucasGraciaRivasRoman-RoyVilarino22}.

The Marsden--Weinstein--Meyer theorem \cite{MarsdenWeinstein74,Meyer73} states that under suitable regularity conditions a symplectic Hamiltonian action $G\curvearrowright(M,\omega)$ with moment map $\mu\colon M\to\g^*$ determines a~canonical symplectic structure on the reduced space $M_0=\mu^{-1}(0)/G$. Instead, we could consider the Poisson algebra $C^\infty(M,\omega)$ and try to reduce it directly to an algebra $C^\infty(M,\omega)_0$ isomorphic to $C^\infty(M_0,\omega_0)$, without passing through the reduced space $(M_0,\omega_0)$. It turns out that the conditions under which the existence of $(M_0,\omega_0)$ is ensured are more restrictive than those admitting a natural algebraic reduction  $C^\infty(M,\omega)_0$. A symplectic \emph{observable reduction scheme} is a procedure for defining a reduced algebra of observables $C^\infty(M,\omega)_0$ without assuming the existence of $(M_0,\omega_0)$.

In this paper, we explore observable reduction in the multisymplectic framework, that is, to the setting of a smooth manifold $M$ equipped with a closed and nondegenerate $(n+1)$-form $\omega\in\Omega^{n+1}(M)$. In particular, given a multisymplectic Hamiltonian action $G\curvearrowright(M,\omega)$ and an associated covariant moment map $\mu\in\Omega^{n-1}(M,\g^*)$, we reduce the $L_\infty$-algebra of observables $\Ham_\infty(M,\omega)$ to obtain a reduced space $\Ham_\infty(M,\omega)_0$ that canonically includes in $\Ham_\infty(M_0,\omega_0)$ whenever the geometric reduced space $(M_0,\omega_0)$ exists. In fact, our construction is rather more general than this. We define a reduction of $\Ham_\infty(M,\omega)$ with respect to any Lie algebra action $\g\curvearrowright M$ and subset $N\subset M$ satisfying mild compatibility conditions.

The main difficulty in the construction stems from the fact that the space of observables $\Ham_\infty(M,\omega)$ on an $n$-plectic manifold $(M,\omega)$ exhibits the rather technical structure of an \emph{$L_\infty$-algebra}, or \emph{homotopy Lie algebra} and not the structure of a Poisson algebra as in the classical symplectic setting. The $L_\infty$-algebra structure of $\Ham_\infty(M,\omega)$ was first identified in \cite{Rogers10}. Currently, the commutative multiplication of the observable algebra in symplectic geometry has no analogue in the multisymplectic world, and the search for it is an open problem in multisymplectic geometry. However, the $L_\infty$-structure alone has many interesting consequences and has since become an object of interest in its own right, for example, in its role in the construction of \emph{homotopy moment maps} \cite{CalliesFregierRogersZambon16,FregierLaurent-GengouxZambon15,Herman18,ShabaziZambon16}, conserved quantities \cite{RyvkinWurzbacherZambon20}, and prequantizations \cite{FiorenzaRogersSchreiber14, SevestreWurzbacher21}.

In the symplectic (i.e., $1$-plectic) setting, our construction agrees with the Dirac, \'Sniatycki--Weinstein, and Arms--Cushman--Gotay reduction schemes whenever the Marsden--Weinstein reduced space $(M_0,\omega_0)$ exists and, in this situation, returns the Poisson algebra $C^\infty(M_0,\omega_0)$. We~review these reduction schemes in Appendix \ref{Sec:SymSingRed}, and emphasize that even in the symplectic case our construction is distinct from each of them.

\subsection*{Types of reduction}

We now review a few terminological conventions of symplectic reduction that carry over in a~natural way to the multisymplectic setting:
\begin{itemize}\itemsep=0pt
	\item A reduction scheme is said to be \emph{geometric} when it produces a reduced symplectic mani\-fold~$(M_0,\omega_0)$.
	\item It is said to be \emph{algebraic}, or \emph{observable}, when it returns a reduced space of observables $C^\infty(M,\omega)_0$, without necessarily exhibiting an underlying reduced symplectic manifold.
\end{itemize}
In mechanical terms, the former is a reduction of \emph{states}, while the latter is a reduction of \emph{observables}.

Every reduction scheme we consider involves in an essential way a subset $N\subset M$. Traditionally, $N$ is the vanishing locus of a family of constraint functions $\varphi_a \in C^\infty(M)$. Such a~scheme may additionally utilize an action $G\curvearrowright M$, historically encoding a group of symmetries of a~generalized momentum phase space. This brings us to a second distinction:
\begin{itemize}\itemsep=0pt
	\item We call a reduction scheme \emph{constraints-based} when it takes as input a subset $N\subset M$.
	\item We call a constraints-based reduction scheme \emph{symmetry-based} if it additionally takes as input a compatible Lie group action $G\curvearrowright (M,\omega)$.
\end{itemize}

Finally, we make a distinction based on the smoothness conditions imposed on the constraint set $N\subset M$:
\begin{itemize}\itemsep=0pt
	\item A reduction scheme is \emph{regular} when the constraint set $N\subset M$ is required to be a smoothly embedded submanifold.
	\item Otherwise, it is \emph{singular}.
\end{itemize}

For example, the Marsden--Weinstein reduction \cite{MarsdenWeinstein74} is a regular symmetry-based geometric reduction scheme, while the Arms--Cushman--Gotay \cite{ArmsCushmanGotay91} reduction is a singular geometric symmetry-based scheme. In Table \ref{tab:SymplecticReductionSchemes}, we indicate a few instances of symplectic reduction schemes, arranged according to the classification above, we provide a survey of these approaches in Appendix \ref{Sec:SymSingRed}.
\begin{table}[h!]
	\centering \renewcommand{\arraystretch}{1.2}
	\caption{Some well-known symplectic reduction schemes.}
	\label{tab:SymplecticReductionSchemes}\vspace{1mm}

	\begin{tabular}{l|c|c|}
		\cline{2-3}
		&  {constraints-based}
		&  {symmetry-based}
		\\ \hline
		\multicolumn{1}{|c|}{\begin{tabular}[c]{@{}l@{}}  geometric\\   reduction\end{tabular}}
		& \begin{tabular}[c]{@{}l@{}}
			Arms--Gotay--Jennings \cite{ArmsGotayJennings90}
		\end{tabular}
		& \begin{tabular}[c]{@{}l@{}}
			Marsden--Weinstein \cite{MarsdenWeinstein74}
			\\ Arms--Cushman--Gotay \cite{ArmsCushmanGotay91}
		\end{tabular}
		\\ \hline
		\multicolumn{1}{|c|}{\begin{tabular}[c]{@{}l@{}}  observable\\   reduction\end{tabular}}
		& 			Dirac \cite{Dirac64}
		& \begin{tabular}[c]{@{}l@{}}
			\'{S}niatycki--Weinstein \cite{SniatyckiWeinstein83}
			\\ $L_\infty$ [Section~\ref{Sec:CoreSection}]
		\end{tabular}
		\\ \hline
	\end{tabular}
\end{table}

\subsection*{Summary and results}

We define a symmetry-based observable reduction scheme in the multisymplectic setting.
The parameters of the reduction consist of a suitably compatible (pre-)multisymplectic manifold $(M,\omega)$, constraint set $N\subset M$, and Lie algebra action $\g\curvearrowright M$. The result is an $L_\infty$-algebra $\Ham_\infty(M,\omega)_N$.

Key features of this reduction include:
\begin{itemize}\itemsep=0pt
	\item The reduction is applicable for any subset $N\subset M$, without conditions on smoothness or type of singularity.
	\item The action $\g\curvearrowright M$ is required to preserve $N$ in a weak sense: It is the ideal of functions vanishing on $N$, rather than $N$ itself, that must be preserved.
	\item The action need not preserve $\omega$: It suffices for $\omega$ to satisfy a strictly weaker condition of \emph{reducibility}.
\end{itemize}

The reduction naturally applies to constraint sets $N\subset M$ induced by covariant moment maps, and hence also induces a symmetry-based observable reduction scheme that interacts with geometric multisymplectic reduction.
It may come as a surprise that our construction showcases a novel behavior even in the symplectic setting. Namely, it produces a reduced Poisson algebra that is not isomorphic to the observable algebras obtained by other singular reduction schemes, for a non regular constraint set $N$.
We emphasize that the proposed $L_\infty$ observable reduction scheme does not rely on the non-degeneracy of the closed form,
however, since many relevant cases of interest involve non-degenerate differential forms, we refer to the scheme as multisymplectic even though it is valid in larger generality.

The paper is organized as follows.
We begin in Section \ref{Sec:MultiReview} providing a brief review of multisymplectic geometry, including (pre-)multisymplectic manifolds $(M,\omega)$, and their observable $L_\infty$-algebras $\Ham_\infty(M,\omega)$. We conclude this section by reviewing covariant moment maps and geometric multisymplectic reduction.

In Section \ref{Sec:CoreSection}, we present our main construction. Here we introduce the notions of reducibility for differential forms, vector fields, and multisymplectic observables. We prove that the space of reducible observables ${\Ham_\infty(M,\omega)_{[N]}}$ is an $L_\infty$-algebra and show that the space of observables $I_{\Ham_\infty}(N)$ that should vanish after reduction is an $L_\infty$-ideal. Hence their quotient is a well-defined $L_\infty$-algebra:

\setcounter{section}{2}
\setcounter{Theorem}{20}
\begin{Definition}
	The \emph{reduction} of $\Ham_\infty(M,\omega)$ with respect to $\g\curvearrowright (N\subset M)$ is the $L_\infty$-algebra
	\[
	\Ham_\infty(M,\omega)_N = \frac{\Ham_\infty(M,\omega)_{[N]}}{I_{\Ham_\infty}(N)}.
	\]
\end{Definition}

We then turn to symmetry-based multisymplectic observable reduction. We specialize our reduction procedure to the case where $N\subset M$ is a level set of a covariant moment map, and we compare the reduced space of observables with the geometric reduction $(M_N,\omega_N)$. Whenever the latter exists, we have:

\setcounter{section}{2}
\setcounter{Theorem}{37}
\begin{Theorem}
	The geometric reduction map
	\begin{align*}
		r_N\colon \	\Ham_\infty(M,\omega)_{[N]}	&\to		\Ham_\infty(M_N,\omega_N),		\\
		(v,\alpha)			&\mapsto	(v_N,\alpha_N)	,		\\
		\alpha		&\mapsto	\alpha_N
	\end{align*}
	is a strict $L_\infty$-morphism with kernel $I_{\Ham_\infty}(N)$. In particular, there is a natural inclusion of $L_\infty$-algebras
	\[
	\Ham_\infty(M,\omega)_N = \frac{\Ham_\infty(M,\omega)_{[N]}}{I_{\Ham_\infty}(N)} \;\xhookrightarrow{\;\bar{r}_N\;}\; \Ham_\infty(M_N,\omega_N).
	\]
\end{Theorem}

In Section \ref{Sec:Examples}, we turn to interesting and illuminating classes of examples of reduction. We first look at the important special case of symplectic manifolds. First, we show that our reduction scheme applied to the Poisson algebra $C^\infty(M,\omega)$ is canonically Poisson.

\setcounter{section}{3}
\setcounter{Theorem}{0}
\begin{Theorem}
	If $\g\curvearrowright M$ is tangent to $N\subset M$, and if the symplectic structure $\omega\in\Omega^2(M)$ is reducible, then the reduction $L_\infty(M,\omega)_N$ inherits a natural Poisson algebra structure from $C^\infty(M,\omega)$.
\end{Theorem}

We then establish that, while the $L_\infty$-reduction procedure is distinct from the Dirac, \'Snia\-tycki--Weinstein, and Arms--Cushman--Gotay schemes, they all coincide whenever the Marsden--Wein\-stein quotient exists:

\setcounter{section}{3}
\setcounter{Theorem}{16}
\begin{Theorem}
	Let $G\action (M,\omega)$ be a symplectic Hamiltonian action and suppose that $0\in\g^*$ is a regular value of the moment map $\mu\colon M\to\g^*$.
	If $G\curvearrowright M$ is free and proper, then the {\rm[$L_\infty$]}, {\rm[\'SW]}, {\rm[D]}, and {\rm [ACG]} reductions are equal.
	In particular, each is isomorphic to the Poisson algebra $C^\infty(M_0,\omega_0)$ of smooth functions on the Marsden--Weinstein reduced space.
\end{Theorem}

In Section \ref{Ssec:coordinate_cross}, we compute the $L_\infty$-reduction for a simple but illustrative example in which the symplectic form is merely reducible but not invariant.

In Sections \ref{Ssec:multicotangent} and \ref{Ssec:Multiphase}, we consider the celebrated examples of multicotangent bundles and multiphase spaces, which underlie the covariant Hamiltonian approach to classical field theories. Even in this natural and motivating setting, a complete description of the reduced space appears to be highly nontrivial. In spite of this, we show that a particular natural class of observables associated to the underlying covariant configuration bundle is always reducible:

\setcounter{section}{3}
\setcounter{Theorem}{20}
\begin{Theorem}
	If $v\in\X(M)$ preserves $\X_\g(E)$, the $C^\infty(M)$-module generated by the fundamental distribution $\underline\g$, and if $v\in\X(E)$ is $(E\to\Sigma)$-projectable, then $(\tilde{v},\iota_{\tilde v}\theta)\in\Ham_\infty^0\bigl(\Lambda_1^nT^*E,\omega\bigr)$ is reducible.
\end{Theorem}

Our theoretical development is supplemented in Section~\ref{Ssec:2dim-ScalarField} with explicit computations of reducible observables for a particular scalar field theory.

We conclude in Section \ref{Sec:Outlook} with a discussion of three topics for future development. First, we consider the setting of spaces of connections under the action of a gauge group, following~\mbox{\cite[Section 10]{CalliesFregierRogersZambon16}}. Second, we consider the application of our work to the historical source for multisymplectic geometry, classical field theory. Finally, we consider the extension of our methods and results to multisymplectic quantization.

In Appendix \ref{Sec:SymSingRed}, we review classical symplectic reduction with respect to a moment map, and a survey of three approaches to the observable reduction of symplectic manifolds (cf.\ Table~\ref{tab:SymplecticReductionSchemes}). In addition to contextualizing our reduction scheme in terms of symplectic predecessors, this exposition will serve as a reference in Section~\ref{Ssec:LinftyRed-SympCase} for comparison with our new multisymplectic reduction scheme.

\setcounter{section}{1}
\setcounter{Theorem}{0}

\subsection*{Notation and conventions}

Our sign conventions are chosen to be broadly consistent with \cite{CalliesFregierRogersZambon16}.

All manifolds are assumed to be $C^\infty$ and paracompact.
We denote by $G\action M$ the smooth action of a Lie group on a manifold.
Except where explicitly stated otherwise, all actions are on the left.
Induced actions on spaces of forms $G\curvearrowright\Omega(M)$ are given in the usual manner by inverse pullback.
The induced infinitesimal action $\g\curvearrowright M$ is the assignment of fundamental vector fields
\[
	\morphism{\underline{\,\cdot\,}}
	{\mathfrak{g}}
	{\mathfrak{X}(M),}
	{\xi}
	{\underline\xi,}
\]
where
\[
	\underline\xi_x = \dfrac{\d}{\d t}\,\exp{(-t \xi)}\hspace{1pt}x \,\Big|_{t=0}
	\qquad \forall x \in M.
\]
Consequently, $\xi\mapsto\underline\xi$ is a Lie algebra homomorphism. We write $\L_\xi$ and $\iota_\xi$ for $\L_{\underline\xi}$ and $\iota_{\underline\xi}$ throughout.

On a pre-$n$-plectic manifold $(M,\omega)$, a Hamiltonian form $\alpha\in\Omega^{n-1}(M)$ and any associated Hamiltonian vector field $v_\alpha \in\X(M)$ are related by the identity $\d\alpha=-\iota_{v_\alpha}\omega$.
The relation between a $\omega$ and a premultisymplectic potential $\theta\in\Omega^n(M)$ is $\omega=\d\theta$.

The Leibniz bracket on $\Omega_\ham^{n-1}(M,\omega)$ is given by $\{\alpha,\beta\}= \L_{v_\alpha}\beta$.
In particular, in the symplectic setting the Poisson bracket is $\{f,h\}=\omega(v_f,v_h)$.
Consequently, the assignment of Hamiltonian vector fields $C^{\infty}(M,\omega)\to \mathfrak X(M)$ is a Lie algebra homomorphism.

We denote by $\alpha_\xi=\langle\alpha,\xi\rangle$ the contraction of a $\g^*$-valued form $\alpha\in\Omega(M,\g^*)$ with an element $\xi\in\g$. Consequently, $\ad_\xi^*\alpha_\zeta = -\alpha_{[\xi,\zeta]}$ where $\ad^*\colon \g\curvearrowright\g^*$ is the coadjoint action.

A comprehensive glossary of symbols is provided in Table \ref{tab:symbols}.

\begin{table}[t]
	\centering\renewcommand{\arraystretch}{1.22}
	\caption{Index of notation.}	\label{tab:symbols}

\vspace{1mm}

	\begin{tabular}{lll}
 \hline
		symbol			&meaning			&reference		\\ \hline
		%
		$\omega$		&(pre)multisymplectic form		&Definition~\ref{def:MultisymplecticManifold}		\\
		$\mu$			&covariant moment map			&Definitions~\ref{def:covariant_moment_map} and \ref{def:momentmap}		\\
		$\widetilde\mu$		&comoment map				&Definitions~\ref{def:multisymplectic_comoment_map} and \ref{def:comomentmap}		\\
		$\X(M,\omega)$		&multisymplectic vector fields		&Definition~\ref{def:multisymplectic_vector_field}			\\
		$\X_\ham(M,\omega)$	&Hamiltonian vector fields		&Definition~\ref{def:multisymplectic_Hamiltonian_vector_field}	\\
		$L_\infty(M,\omega)$	&associated $L_\infty$-algebra		&Definition~\ref{def:RogersAlgebra}			\\
		$\varsigma(k)$		&total Koszul sign			&Definition~\ref{def:RogersAlgebra}	\\
		$l_k$			&$k$-th multibracket of $L_\infty(M,\omega)$	&Definition~\ref{def:RogersAlgebra}	\\
		$\Ham_\infty(M,\omega)$	&$L_\infty$-algebra of observables
										&Definition~\ref{def:algebra_of_observables}	\\
		$\tilde{l}_k$		&$k$-th multibracket of 	$\Ham_\infty(M,\omega)$	&Definition~\ref{def:algebra_of_observables}	\\
		$\{\,,\}$		&Leibniz bracket on $\Omega_\ham^{n-1}(M,\omega)$	&Definition~\ref{def:LeibObs}	\\
		$\theta$		&multisymplectic potential		&Lemma~\ref{lem:potinducesmom}, Section~\ref{Ssec:multicotangent}	\\
		$(M_\phi,\omega_\phi)$	&reduction of $(M,\omega)$ at level $\phi$		&Theorem~\ref{thm:CaseyReduction}	\\
		$I_N$			&constraints ideal in $C^\infty(M)$ 	&Definition~\ref{def:constraintideal}	\\
		$\X_N(M)$		&vector fields tangent to $N$		&Definition~\ref{def:tangentvf}		\\
		$I_\X(N)$		&vanishing ideal in $\X(M)$		&Definition~\ref{def:tangentvf}		\\
		$\g\curvearrowright(N\subset M)$
					&action tangent to $N$			&Definition~\ref{def:Lie_algebra_action_tangent_to_N}	\\
		$I_\Omega(N)$		&vanishing de Rham ideal		&Definition~\ref{def:vanishing_de_Rham_ideal}		\\
		$\Omega(M)_{[N]}$	&reducible forms			&Definition~\ref{def:reducible_form}		\\
		$\X(M)_{[N]}$		&reducible vector fields		&Definition~\ref{def:reducible_vector_field}	\\
		$\fgmodule(M)$		&fundamental submodule			&Definition~\ref{def:reducible_vector_field}	\\
		$\Ham_\infty(M,\omega)_{[N]}$	
					&reducible observables			&Definition~\ref{def:reducibleObservables} 	\\
		$I_{\Ham_\infty}(N)$	&vanishing observable ideal		&Definition~\ref{def:VanishIdeal}		\\
		$\Ham_\infty(M,\omega)_N$
					&$L_\infty$-reduction along $N$		&Definition~\ref{def:LinftyReductions}		\\
		$\Ham_\infty(M,\omega)_\phi$
					&$L_\infty$-reduction at level $\phi$	&Definition~\ref{def:observable_reduction_wrt_moment_map}		\\
		$\underline\g$		&fundamental distribution		&Lemma~\ref{lem:RedIngredientsRegularCase}, Remark~\ref{rem:understanding_reducibility} 	\\
		$r_N$			&geometric reduction map		&Definition~\ref{def:geometric_reduction_map}\\
		$(M_N,\omega_N)$	&reduction of $(M,\omega)$ along $N$	&Theorem~\ref{def:geometric_reduction}	\\
		$[\Ham_\infty,R]$ 	&residue defect				&Definition~\ref{def:residue_defect}		\\
		$F_N$			&first-class function set		&Definition~\ref{def:IclassFunc}		\\
		$\mathcal{Q}$		&Casimir functions			&Definition~\ref{def:Casimir_function}	\\
		$I_\mu$			&momentum ideal				&Definition~\ref{def:momentumideal}		\\
		$\Lambda^k_1T^*E$	&multiphase space			&Section~\ref{Ssec:Multiphase}		\\
				$C^\infty(M,\omega)$	&Poisson algebra of observables	&Appendix~\ref{Sec:SymSingRed}		\\
		$\mathcal{N}(S)$	&normalizer of $S\subset C^\infty(M,\omega)$	&Theorem~\ref{thm:SW-red}	\\
		$\mathcal{O}(N)$	&Dirac observables			&Definition~\ref{def:Dirac_observables}			\\
		$S^G$			&fixed-point set of $G\curvearrowright S$	&Appendix~\ref{Ssec:ACG-red}		\\
	\end{tabular}
\end{table}

\section{Review of multisymplectic geometry}\label{Sec:MultiReview}

In this section, we provide relevant background on multisymplectic geometry. Our aim is to define the associated $L_\infty$-algebra $L_\infty(M,\omega)$ and the $L_\infty$-algebra of observables $\Ham_\infty(M,\omega)$, and to recall the geometric multisymplectic reduction procedure of Theorem \ref{thm:CaseyReduction}.

\subsection[Multisymplectic manifolds and the L\_infty-algebra of observables]{Multisymplectic manifolds and the $\boldsymbol{L_\infty}$-algebra of observables}

Fix a smooth manifold $M$.

\begin{Definition}[multisymplectic manifold \cite{CantrijnIbortdeLeon99}]\label{def:MultisymplecticManifold}
	A \emph{pre-$n$-plectic structure} on $M$ is a closed $(n+1)$-form $\omega\in\Omega^{n+1}(M)$. Without reference to the degree $n+1$, we say that $\omega$ is a \emph{premultisymplectic structure} on $M$. If additionally $\omega$ is nondegenerate in the sense that
	\begin{displaymath}
		\morphism{\omega^\flat}
		{TM}{\Lambda^nT^*M,}
		{v}{\iota_v\omega}
	\end{displaymath}
	is an inclusion of vector bundles, then we call $\omega$ an \emph{$n$-plectic} or a \emph{multisymplectic structure} on~$M$.
\end{Definition}

Symplectic manifolds, manifolds equipped with volume forms, multicotangent bundles, and multiphase spaces, which we will discuss in the sequel, are all examples of multisymplectic manifolds.

Multisymplectic manifolds form a category with morphisms given as follows.

\begin{Definition}
	A \emph{multisymplectic map} is a smooth map of multisymplectic manifolds
\[
\Psi\colon \ (M,\omega) \to (M',\omega')
\]
 such that $\Psi^* \omega' = \omega$. A \emph{multisymplectomorphism} is a multisymplectic diffeomorphism.
\end{Definition}

The multisymplectomorphisms thus comprise the global symmetry group of $(M,\omega)$. The infinitesimal symmetries are given as follows.

\begin{Definition}\label{def:multisymplectic_vector_field}
	We say that $v\in\X(M)$ is a \emph{multisymplectic vector field} when $\L_v\omega = 0$. We denote the space of multisymplectic vector fields by $\X(M,\omega)$.
\end{Definition}

Observe that $v$ is multisymplectic precisely when $\d\iota_v\omega=\L_v\omega=0$, that is, precisely when~$\iota_v\omega$ is closed. As in the symplectic setting, we distinguish those vector fields $v$ for which $\iota_v\omega$ is exact.

\begin{Definition}\label{def:multisymplectic_Hamiltonian_vector_field}
	Consider a pre-$n$-plectic manifold $(M,\omega)$. We say that $v\in\X(M)$ is a \emph{Hamiltonian vector field} when $\d\alpha=-\iota_v\omega$ for some $\alpha\in\Omega^{n-1}(M)$. In this case, we say that $\alpha$ is a~\emph{Hamiltonian form} associated to $v$, and we write $\X_\ham(M,\omega)$ and $\Omega_\ham^{n-1}(M,\omega)$ for the spaces of Hamiltonian vector fields and forms, respectively.
\end{Definition}

\begin{Remark}
	The identity $\d\alpha +\iota_v\omega=0$ is known as the \emph{Hamilton--De Donder--Weyl $($HDDW$)$ equation} (see, e.g., \cite{Helein04,RyvkinWurzbacher18} and references therein).
\end{Remark}

\begin{Remark}
	When $\omega$ is multisymplectic, the vector field $v$ associated to $\alpha$ is unique. In this case, we say that $v$ is \emph{the} Hamiltonian vector field associated to $\alpha$.
\end{Remark}

Observe that
\begin{align*}
	&\mathfrak{X}(M,\omega) =
	\bigl\lbrace v \in \mathfrak{X}(M) \mid  \iota_v\omega \in \Omega_{\mathrm{cl}}^n(M) \bigr\rbrace,
	\\
	&\mathfrak{X}_{\ham}(M,\omega) =
	\bigl\lbrace v \in \mathfrak{X}(M) \mid  \iota_v\omega \in \Omega_{\mathrm{ex}}^n(M) \bigr\rbrace,
\end{align*}
from which it follows that Hamiltonian vector fields are multisymplectic. Given $v\in\X_\ham(M,\omega)$, an associated Hamiltonian form $\alpha$ is determined up to closed forms.

The space $\Omega_\ham^{n-1}(M,\omega)$ is not generally a Poisson algebra. However, it can be extended to $\Omega_\ham^{n-1}(M,\omega)\oplus\Omega^{\leq n-2}(M)[n{-}1]$, which possesses a natural $L_\infty$-algebra structure. When $\omega$ is degenerate there are in fact two distinct constructions, $L_\infty(M,\omega)$ and $\Ham_\infty(M,\omega)$.

\begin{Definition}[{\cite[Theorem~5.2]{Rogers10}, see also \cite{BarinchFulpLadaStasheff98}}]\label{def:RogersAlgebra}
	Given a pre-$n$-plectic manifold $(M,\omega)$, the associated $L_\infty$-algebra $L_{\infty}(M,\omega)=\bigl(L,\{l_k\}_{k\geq 1}\bigr)$ comprises
	\begin{itemize}\itemsep=0pt
		\item
		the underlying graded vector space $L$, where
	\begin{equation*}
		L^i=\begin{cases}
			\Omega_{\ham}^{n-1}(M,\omega)
			& \text{if } i=0,
			\\
			\Omega^{n-1+i}(M)
		 	&  \text{if } 1-n \leq i\leq -1,
		 	\\			
			0 & \text{otherwise,}
		\end{cases}
	\end{equation*}

	\item
		$n+1$ nontrivial multibrackets $\bigl\lbrace l_k \colon L^{\wedge k} \to L\bigr\rbrace_{1\leq k\leq n+1}$, given by
	\begin{displaymath}
		l_1(\alpha) =
		\begin{cases}
			0 & \text{if~} |\alpha| = 0	,
			\\
			\d \alpha &  \text{if~} |\alpha| \leq -1,
		\end{cases}
	\end{displaymath}
	and, for $ 2 \leq k \leq n+1$, as
	\begin{displaymath}
		l_k(\alpha_1,\dots,\alpha_k) =
		\begin{cases}
			\varsigma(k) \iota( v_{\alpha_1}\wedge\dots\wedge v_{\alpha_k})~\omega
			&\text{if} \ |\alpha_i|=0 \text{ for } 1\leq i \leq k,
			\\
			0 & \text{otherwise.}		
		\end{cases}
	\end{displaymath}
	\end{itemize}
In the above equation, $ v_{\alpha_k}$ denotes any Hamiltonian vector field associated to $\alpha_k\in \Omega^{n-1}_{\ham}(M,\omega)$ and \smash{$\varsigma(k) = - (-1)^{\frac{k(k+1)}{2}}$} is the total Koszul sign. Note that we contract multivector fields according to the rule
	 \[
		 \iota( v_{\alpha_1}\wedge\dots\wedge v_{\alpha_k}) \omega=\iota_{ v_{\alpha_k}}\cdots\iota_{ v_{\alpha_1}}\omega=\omega( v_{\alpha_1},\dots, v_{\alpha_k}, \,\cdot\,, \dots,\,\cdot\,).
	 \]
\end{Definition}

\begin{Definition}[{$L_\infty$-algebra of observables \cite[Theorem 4.7]{CalliesFregierRogersZambon16}}]\label{def:algebra_of_observables}
	The \emph{$L_{\infty}$-algebra of observables} $\Ham_\infty(M,\omega)=\bigl(\Ham,\bigl\{\tilde{l}_k\bigr\}_{k\geq 1}\bigr)$ associated to the pre-$n$-plectic manifold $(M,\omega)$ consists of
	\begin{itemize}\itemsep=0pt
		\item the underlying graded vector space $\Ham$, where
		\[
			\Ham^i=\begin{cases}
				\{(v,\alpha)\mid \d\alpha=-\iota_v\omega\}\subset \X_\ham(M,\omega)\oplus\Omega_{\ham}^{n-1}(M,\omega)
				&  \text{if } i=0,
				\\
				\Omega^{n-1+i}(M)
				&   \text{if } 1-n \leq i\leq -1,
				\\			
				0 & \text{otherwise,}
			\end{cases}
		\]

	\item $n+1$ multibrackets $\bigl\lbrace \tilde{l}_k\colon \Ham^{\wedge k} \to \Ham\bigr\rbrace_{1\leq k\leq n+1}$, where all non trivial components are given, in the unary case, by
		\[
			\tilde{l}_1(\alpha) =
			\begin{cases}
				(0,\d\alpha)	& \text{if~} |\alpha| = -1,	\\
				\d\alpha	& \text{if~} |\alpha| <-1,
			\end{cases}
		\]
		in the binary case by
		\[
			\tilde{l}_2 ((v_1,\alpha_1),(v_2,\alpha_2) ) =
				 ([v_1,v_2], \varsigma(k) \iota(v_1\wedge v_2)~\omega )
				,
		\]
		and for each $k$-ary multibracket with $3\leq k\leq n+1$ are given by
		\[
			\tilde{l}_k ((v_1,\alpha_1),\dots,(v_k,\alpha_k) ) =
				\varsigma(k) \iota(v_1\wedge\dots\wedge v_k)~\omega
		\]
		with  $(v_1,\alpha_1),\dots,(v_k,\alpha_k)\in\Ham^0$.
	\end{itemize}
\end{Definition}

\begin{Remark}
We briefly outline some basic definitions in homotopy algebras, tailoring the exposition to our purposes and referring to \cite{LadaMarkl95} and \cite[Section~3]{Lazarev14} for further background material.

Recall that a $L_\infty$-algebra consists of a graded vector space $L$ together with homogeneous, graded skew-symmetric, multilinear maps $\lambda_k\colon L^{\otimes k} \to L$ (with $k\geq 1$) satisfying the so-called \emph{higher Jacobi relations}, see \cite[Definition~2.1]{LadaMarkl95}.
We denote collectively as $\lambda$ the collection $\{\lambda_k\}_{k\geq 1}$ of all multilinear maps, also called \emph{multibrackets}.

Given two $L_\infty$-algebras $(L,\lambda)$ and $(M,\mu)$, a \emph{strict $L_\infty$-morphism} $\psi\colon  (L,\lambda)\to(M,\mu)$ is a~graded map from $L$ to $M$ such that $\mu_k\circ \psi^{\otimes k} = \psi \circ \lambda_k$ for each $k\geq 1$. In particular, \emph{$L_\infty$-subalgebra} of $(L,\lambda)$ is a graded vector subspace $V$ such that the injection map $j\colon V\hookrightarrow L$ is a strict $L_\infty$-morphism $\bigl(V,\bigl\{\lambda_k\circ j^{\otimes k}\bigr\}_{k\geq 1}\bigr)\to (L,\lambda)$. Both the kernel and the image of a strict $L_\infty$-morphism are the kernel and the image of the underlying graded map endowed with the $L_\infty$-algebra structure inherited from its domain and codomain.

A \emph{$($strict$)$ $L_\infty$-ideal} of $(L,\lambda)$ is a graded vector subspace $I\subset L$ for which ${\lambda_k(x,y_1,\dots,y_{k{-}1})\! \in\! I}$ for every $x\in I$, $y_i \in L$, and $k\geq 1$. In particular, $I\subset L$ is an $L_\infty$-subalgebra. Moreover, there is a natural correspondence between strict $L_\infty$-ideals $I\subset L$ and kernels of strict $L_\infty$-morphisms $\psi_I\colon L\to L/I$. It is readily shown that the linear quotient $L/I$ inherits the multibrackets $\lambda_{L/I,k}([y_1],\dots,[y_k])=[\lambda_k(y_1,\dots, y_k)]$ for any $[y_i]\in L/I$.

Expanding on this, given a short exact sequence in the category of $L_\infty$-algebras with strict $L_\infty$-morphisms,
\[
	0\to N\to M \xrightarrow{\pi} L \to 0
	,
\]
we have that $N \cong \ker\pi$ is an ideal and that $\pi$ induces an isomorphism $\bar\pi\colon M/N \cong L$. That is, the first isomorphism theorem for Lie algebras generalizes in a straightforward manner to the setting of strict $L_\infty$-morphisms.

We emphasize here that this discussion applies only to the case of strict $L_\infty$-morphisms. We do not consider non-strict $L_\infty$-morphisms -- including, prominently, homotopy moment maps -- at any point in this paper.
\end{Remark}

\subsection{Covariant moment maps and multisymplectic reduction}\label{Ssec:MultiRed}

In addition to extending to an $L_\infty$-algebra, the space $\Omega^{n-1}_{\ham}(M,\omega)$ possesses a natural Leibniz algebra structure. Compatibility conditions in terms of this structure are invoked in the geometric multisymplectic reduction procedure.

\begin{Definition}[Leibniz algebra \cite{Loday93}]
	A \emph{$($left$)$ Leibniz algebra} comprises a vector space $V$ and a bilinear map $\{\,,\}\colon V\times V\to V$ satisfying the Leibniz equation
	\[
		\{\alpha,\{\beta,\gamma\}\} = \{\{\alpha,\beta\},\gamma\} + \{\beta,\{\alpha,\gamma\}\}
	\]
	for any $\alpha,\beta,\gamma \in V$.
\end{Definition}

\begin{Definition}[Leibniz algebra of observables]\label{def:LeibObs}
	The space $\Omega_\ham^{n-1}(M,\omega)$ of Hamiltonian forms on a multisymplectic manifold $(M,\omega)$ possesses a natural Leibniz bracket $\{\,,\}$, given by
	\[
		\{\alpha,\beta\} = \L_{v_\alpha}\beta
	\]
	for $\alpha,\beta\in\Omega_\ham^{n-1}(M,\omega)$.
\end{Definition}

When $(M,\omega)$ is premultisymplectic, there is no Leibniz algebra structure on $\Omega_\ham^{n-1}(M,\omega)$, however there still is one on $\Ham^0(M,\omega)$, given by
\[
	\{(u,\alpha),(v,\beta)\}=([u,v], \L_u\beta).
\]

\begin{Remark}[on the two notions of higher observables]
	It is shown in \cite[Proposition~5.2]{Dehling17} that a certain natural extension of the Leibniz algebra of Definition \ref{def:LeibObs}, incorporating all forms of degree strictly less than $n-1$, and the construction of Definition \ref{def:RogersAlgebra} are equivalent as weak $L_\infty$-algebras at least up to the $3$-plectic case. The $2$-plectic case was previously established in~\cite{BaezHoffnungRogers10}.
\end{Remark}

\begin{Definition}\label{def:covariant_moment_map}
	A \emph{covariant moment map} for $G\action(M,\omega)$ is a differential form $\mu\!\in\!\Omega^{n-1}(M,\g^*)$ such that
		\begin{enumerate}\itemsep=0pt
			\item[(i)] $\d\mu_\xi = -\iota_\xi\omega$ for all $\xi\in\g$,
			\item[(ii)] $\mu\colon \Lambda^{n-1}TM\to\g^*$ is $G$-equivariant,
		\end{enumerate}
	where $\mu_\xi=\langle\mu,\xi\rangle$ is the contraction on $\g^*\otimes\g$ and where $G\curvearrowright\g^*$ is the coadjoint action.
\end{Definition}

As in the symplectic case, there is an attendant notion of a comoment map.

\begin{Definition}\label{def:multisymplectic_comoment_map}
	A (covariant) \emph{comoment map} for $G\action(M,\omega)$ is a linear map
\[
\widetilde{\mu}\colon \ \g\to\Omega_\ham^{n-1}(M,\omega)
\] satisfying
	\begin{enumerate}\itemsep=0pt
		\item[(i)] $\d\,\widetilde{\mu}(\xi) = -\iota_{\underline{\xi}}\omega$,
		\item[(ii)] $\widetilde{\mu}([\xi,\zeta]) = \{\widetilde{\mu}(\xi),\widetilde{\mu}(\eta)\}$
	\end{enumerate}
	for all $\xi,\zeta\in\g$.
\end{Definition}

In parallel with the symplectic case, $\tilde\mu$ is a lift in the category of Leibniz algebras of the action $\g\curvearrowright M$ by the assignment of Hamiltonian vector fields:
\[
\begin{tikzpicture}
	\node (A) at (0,0) {$\g$};
	\node (B) at (2,0) {$\X(M)$};
	\node (C) at (2,2) {$\Omega_\ham^{n-1}(M)$};

	\node (D) at (0,-.6) {$\xi$};
	\node (E) at (2,-.6) {$\underline\xi$};

	\node (F) at (3.1,2) {$\alpha$};
	\node (G) at (3.1,0) {$v_\alpha$};

	\draw[->] (A) to (B);
	\draw[->,dashed] (A) to node[above left] {$\xi\mapsto\mu_\xi$} (C);
	\draw[->] (C) to (B);

	\draw[|->] (D) to (E);
	\draw[|->] (F) to (G);
\end{tikzpicture}
\]

To each moment map $\mu\in\Omega^{n-1}(M,\omega)$ there is an associated comoment map
\[
	\morphism{\widetilde\mu}
			{\g}
			{\Omega_\ham^{n-1}(M),}
			{\xi}
			{\mu_\xi,}
\]
and conversely when $G$ is connected.

We identify one class of covariant moment map, to which we return in Section~\ref{Ssec:multicotangent}.

\begin{Lemma}[{\cite[Section~8.1]{CalliesFregierRogersZambon16}}]\label{lem:potinducesmom}
	If $\theta\in\Omega^n(M)$ is a $G$-invariant potential for $\omega$, then $\xi\mapsto\mu_\xi=\iota_\xi\theta$ defines a covariant moment map for $G\curvearrowright(M,\omega)$.
\end{Lemma}

The main result we wish to recall in this subsection is as follows.

\begin{Theorem}[geometric multisymplectic reduction {\cite[Theorem~1]{Blacker21}}]\label{thm:CaseyReduction}
	Let $(M,\omega,G,\mu)$ be a~$n$-plectic Hamiltonian $G$-space with covariant moment map $\mu$, let $\phi\in\Omega^{n-1}(M,\g^*)$ be a~closed form, and let $M_\phi= \mu^{-1}(\phi)/G_\phi$. If $\mu^{-1}(\phi)\subset M$ is an embedded submanifold and $G$ acts freely on $\mu^{-1}(\phi)$, then there is a unique, closed $\omega_\phi\in\Omega^{n+1}(M_\phi)$ satisfying $j^*\omega = \pi^*\omega_\phi$, where $j\colon \mu^{-1}(\phi)\to M$ is the inclusion and $\pi\colon \mu^{-1}(\phi)\to M_\phi$ is the quotient.
\end{Theorem}

In the above theorem, we write $\mu^{-1}(\phi)$ for the equalizer $\{x\in M \mid  \mu_x=\phi_x\}$ of~$\mu$ and~$\phi$. Note that if $\phi$ is identified with its image in $\Lambda^{n-1}T^*M\otimes\g^*$, then $\mu^{-1}(\phi)$ is indeed the pre\-image.\looseness=1

The notion of covariant moment map turns out to be the right prerequisite for geometric reduction in the multisymplectic setting. Perhaps unexpectedly, this is also the case for $L_\infty$-reduction, as we will see in the following section.

\section[Reduction of L\_infty(M,omega) and Ham\_infty(M,omega)]{Reduction of $\boldsymbol{L_\infty(M,\omega)}$ and $\boldsymbol{\Ham_\infty(M,\omega)}$}\label{Sec:CoreSection}

\begin{figure}[t]	\centering
	\begin{tikzpicture}
		\node (A) at (0,10.2) {\fbox{\parbox{5.6cm}{\centering Ideal of functions vanishing on $N$ \\ $I_N$}}};
		\node (B) at (0,8.3) {\fbox{\parbox{4.4cm}{\centering Vector fields tangent to $N$ \\ $\X_N(M)$}}};
		\node (C) at (3,6.5) {\fbox{\parbox{4.2cm}{\centering Forms vanishing along $N$ \\ $I_\Omega(N)$}}};
		\node (D) at (3,4.2) {\fbox{\parbox{2.7cm}{\centering Reducible forms \\ $\Omega(M)_{[N]}$}}};
		\node (E) at (-3,6.3) {\fbox{\parbox{5.3cm}{\centering Fundamental $C^\infty$-module \emph{plus} \\ vector fields vanishing along $N$ \\ $\X_\g(M)+I_\X(N)$}}};
		\node (E') at (-3,4.2) {\fbox{\parbox{3.8cm}{\centering Reducible vector fields \\ $\X(M)_{[N]}$}}};
		\node (F) at (0,2.3) {\fbox{\parbox{3.7cm}{\centering Reducible observables \\ $\Ham_\infty(M,\omega)_{[N]}$}}};
		\node (G) at (0,0.5) {\fbox{\parbox{5.0cm}{\centering Vanishing ideal of observables \\ $I_{\Ham_\infty}(N)$}}};
		\node (H) at (0,-1.35) {\fbox{\parbox{4.3cm}{\centering Observable reduced space \\ $\Ham_\infty(M,\omega)_N$}}};

		\draw[->] (A)-- node[right=1pt]{{\scriptsize used to define}} (B);
		\draw[->] (B)-- (C);
		\draw[->] (C)-- node[right=.25cm]{\rotatebox{-90}{\small $\subset$}} (D);
		\draw[->] (B)-- (E);
		\draw[->] (D)-- (F);
		\draw[->] (E)-- node[left=.25cm]{\rotatebox{-90}{\small $\subset$}} (E');
		\draw[->] (E')-- (F);
		\draw[->] (F)-- node[right=.25cm]{\rotatebox{90}{\small $\subset$}} (G);
		\draw[->] (G)-- (H);

		\node[right=2.95cm] at (A) {\scriptsize\emph{associative ideal}};
		\node[right=2.30cm] at (B) {\scriptsize\emph{Lie algebra}};
		\node[right=2.25cm] at (C) {\scriptsize\emph{differential ideal}};
		\node[right=1.5cm] at (D) {\scriptsize\emph{differential graded subalgebra}};
		\node[left=2cm] at (E') {\scriptsize\emph{Lie subalgebra}};
		\node[left=2.8cm] at (E) {\scriptsize\emph{Lie ideal}};
		\node[right=2cm] at (F) {\scriptsize\emph{$L_\infty$-subalgebra}};
		\node[right=2.65cm] at (G) {\scriptsize\emph{$L_\infty$-ideal}};
		\node[right=2.3cm] at (H) {\scriptsize\emph{$L_\infty$-algebra}};
	\end{tikzpicture}
	\caption{Auxiliary spaces in the construction of $\Ham_\infty(M,\omega)_N$.} \label{fig:intermediate_spaces}
\end{figure}
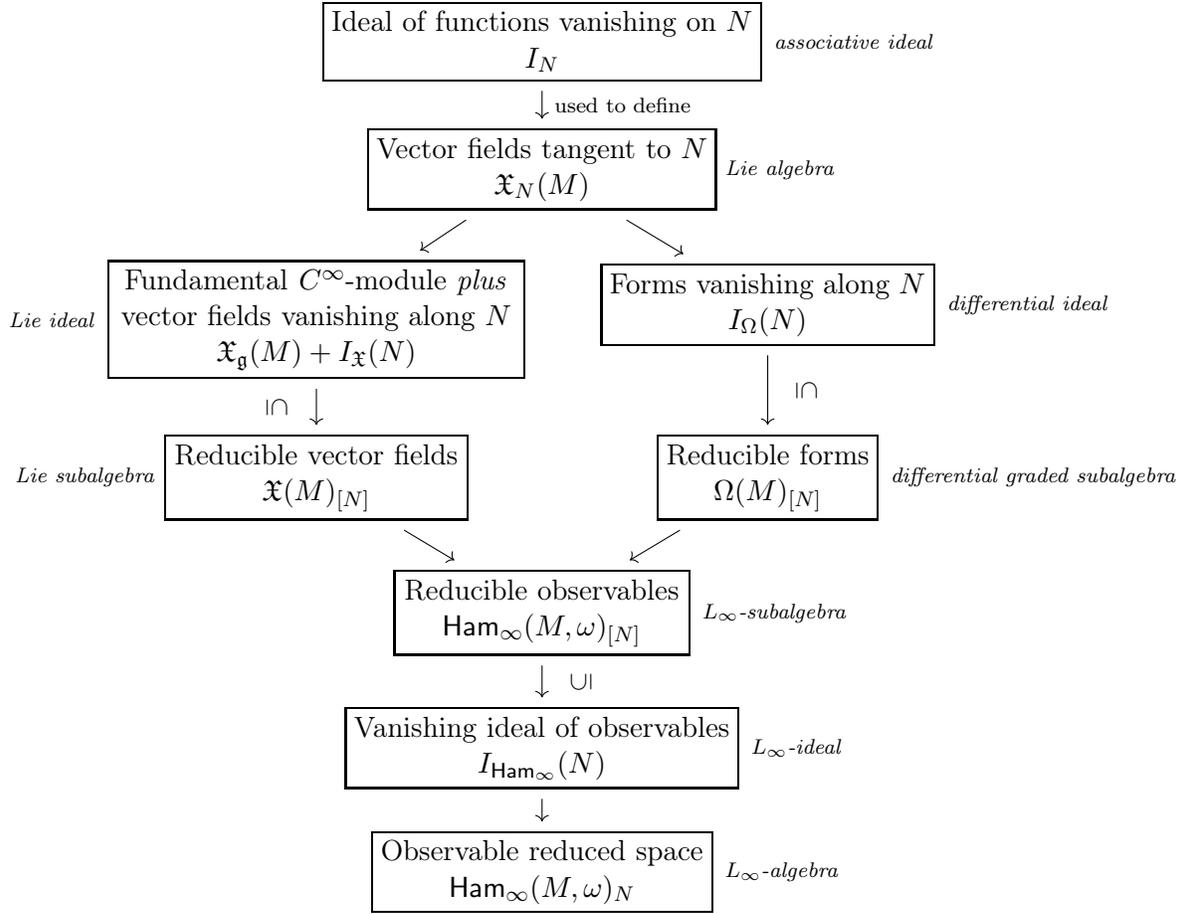

In this section, we present our main construction: The reduction of the $L_\infty$-algebra of observables $\Ham_\infty(M,\omega)$ along a subset $N\subset M$ by a compatible Lie algebra action $\g\curvearrowright M$, in the sense of Definition \ref{def:Lie_algebra_action_tangent_to_N}. We define this $L_\infty$-algebra in terms of a network of auxiliary spaces, as indicated in Figure~\ref{fig:intermediate_spaces}. As we aim for maximal generality, with the mildest possible geometric conditions, these intermediate spaces possess a rather algebraic character. Lemma \ref{lem:RedIngredientsRegularCase} recharacterizes these spaces in a sufficiently regular setting.

In Section~\ref{Ssec:LinftyRed-SympCase}, we show that in the $1$-plectic setting our reduced $L_\infty$-algebra $\Ham_\infty(M,\omega)_\phi$ coincides with the Poisson algebra of observables $C^\infty(M_\phi,\omega_\phi)$ on the symplectic reduced space when the conditions of the Marsden--Weinstein reduction theorem are met.

\subsection{Smooth structures on singular sets}

Let $M$ be a smooth manifold and $N\subset M$ an arbitrary subset.
To avoid the technicalities of making sense of smoothness in the singular setting, is convenient to adopt a more algebraic viewpoint.
The algebra of smooth functions that vanish on $N$ is instrumental in giving an algebraic characterization of the arbitrary subset $N$.

\begin{Definition}[constraints ideal]\label{def:constraintideal}
	We define a \emph{constraint} to be a smooth function vanishing on the constraint set $N$.
	The space of such functions defines an associative algebra ideal in~$C^\infty(M)$,
	\[
		I_N	= 	\{
					f \in C^\infty(M)
				\mid
					f|_N = 0
			 	\}.
	\]
\end{Definition}

An application of a smooth version of Urysohn's lemma provides that $I_N$ determines the closure of the subset $N\subset M$.
Specifically, $p \in \bar{N}$ if and only if $h(p)=0$ for every $h \in I_N$.
In particular, denoting by $\bar{N}$ the closure of $N$, we have that $I_N=I_{\bar N}$.

Finally, observe that when $M$ is endowed with a symplectic structure $\omega$, the associative ideal $I_N\subset C^\infty(M)$ is not generally a Lie subalgebra of $C^\infty(M,\omega)$.

\begin{Remark}\label{rem:generalizedSmootheology}
	There are several ways to specify a notion of smoothness for a possibly singular subset $N$ of the manifold $M$.
For instance, arbitrary constraint sets could be framed as smooth varieties \cite[Section~6]{ArmsCushmanGotay91}, as Sikorski's differential spaces \cite{Sikorski72,Sikorski72a} (see \cite[Section~3]{CushmanSniatycki01} for a quick review and \cite{Stacey08} for a more comprehensive account), or as stratified spaces \cite{SjamaarLerman91}.
	All of these frameworks pass by prescribing what is the algebra of smooth functions on the subset $N\subset M$, namely given by the algebra of Whitney \cite{Whitney34} smooth function
	\[
	C^\infty(N) = C^\infty(M)/I_{N} .
	\]
	Note that, when $N\subset M$ is closed and smoothly embedded, $C^\infty(N)$ can be interpreted as the algebra of smooth functions in the usual sense thus the above identification can be interpreted as the isomorphism induced by restriction map $C^\infty(M)\to C^\infty(N)$ (with kernel given by $I_N$).
	
	Similarly, given a smooth action $G\action M$, the algebra of smooth functions on the orbit spaces $M/G$ and $N/G$ can be given as
\[
	 C^\infty(M/G)= C^\infty(M)^G , \qquad C^\infty(N/G)= C^\infty(M)^G / I_N^G .
\]
\end{Remark}

We will extend this reasoning to differential forms and fields as a way to introduce the notion of tangency along the singular set $N$.

\begin{Definition}\label{def:tangentvf}
	We say that a vector field $u \in \mathfrak{X}(M)$ is \emph{tangent to $N$} when $u\, I_N\subset I_N$. The space of vector fields tangent to $N$ is denoted
	\begin{displaymath}
		\mathfrak{X}_N(M) =
		\lbrace
		u \in \mathfrak{X}(M)
		\mid
		u  I_N \subset I_N
		\rbrace
	\end{displaymath}	
	and the space of vector fields \emph{vanishing along $N$} is written
	\begin{displaymath}
		I_{\mathfrak{X}}(N) =
		\lbrace
		u \in \mathfrak{X}(M)
		\mid
		u \, C^\infty(M) \subset I_N
		\rbrace.
	\end{displaymath}
	We say that $u,v\in\X(M)$ are \emph{equal along $N$} whenever $u-v\in I_\X(N)$.
\end{Definition}

\begin{Definition}\label{def:Lie_algebra_action_tangent_to_N}
	We say that the Lie algebra action $\g\curvearrowright M$ is tangent to $N$ when $\underline\xi\in\X(M)$ is tangent to $N$ for all $\xi$. In this case, we write $\g\curvearrowright(N\subset M)$.
\end{Definition}

\begin{Remark}\label{rem:group_preservation_is_sufficient_for_tangency}
	If $\g\curvearrowright M$ is induced by an action $G\curvearrowright M$ preserving $N\subset M$, then $\g\curvearrowright M$ is tangent to $N$. Indeed, $G\curvearrowright(N\subset M)$ implies that every integral curve $\gamma$ of $\underline\xi\in\X(M)$ that meets $N$ at any point remains on $N$. In particular, each $f\in I_N$ vanishes along $\gamma$ so that $\underline\xi f=0$ along $N$.
\end{Remark}

\begin{Lemma}\label{lem:tangent_to_N_Lie_subalgebra}
	$\tgtvf\subset \X(M)$ is a Lie subalgebra and $I_{\mathfrak{X}}(N)\subset\X_N(M)$ is a Lie ideal.
\end{Lemma}

\begin{proof}
	If $u,v\in\X_N(M)$, then $v I_N\subset I_N$ yields $u(vI_N)\subset I_N$. Likewise $v(uI_N)\subset I_N$, and we conclude that $[u,v]\,I_N\subset I_N$.

	If additionally $v\in I_\X(N)$, then a similar argument yields $[u,v]\,C^\infty(M)\subset I_N$.
\end{proof}

Note that our characterization of tangency is distinctly algebraic. In general, any vector field whose flow preserves $N$ is tangent to $N$, however, there can be vector fields tangent to $N$ that do not preserve $N$. Let us illustrate this with a few examples.

\begin{Remark}
	Let $N\subset M$ be any subset. Since $I_N=I_{\bar N}$, any vector field tangent to the closure~$\bar N$ is tangent to $N$. In particular, when $N\subset M$ is dense (e.g., $N=\mathbb Q\subset  \mathbb R=M$), $I_N=\{0\}$ and all vector fields on $M$ are tangent to $N$.
\end{Remark}

\begin{Example}
	Even when $N$ is closed, there can be vector fields which are tangent to $N$ but whose flows do not preserve $N$. Let $N=(-\infty,0]\subset \mathbb R=M$. Then $I_N$ is the space of functions whose support is in $N$. Applying a vector field to a function cannot increase the support, so any vector field is tangent to $N$. In particular, the vector field $\partial_x$, the flow of which does not preserve $N$.
\end{Example}

\begin{Remark}
	We have chosen the algebraic notion of tangency since it is better adapted for the context of observable reduction. We could alternatively work with vector fields $v\in\X(M)$ that are \emph{geometrically tangent} to $N$ in the sense that the flow of $v$ locally preserves $N$. However, this space is less convenient from a technical standpoint. Indeed, merely to see that the commutator of two geometrically tangent vector fields is again geometrically tangent is far from obvious and follows from \cite{Sussmann73}. This said, we note that both notions of tangency coincide when $N\subset M$ is a~closed embedded submanifold.
\end{Remark}

\begin{Definition}\label{def:vanishing_de_Rham_ideal}
	Define the \emph{vanishing de Rham ideal of $N$} to be the differential ideal $I_\Omega(N)\subset\Omega(M)$ generated by all $\alpha\in\Omega^k(M)$, $k\geq 0$, for which
	\[
		\alpha(u_1,\dots,u_k) \in I_N
	\]
	whenever $u_1,\dots,u_k\in\X_N(M)$.
\end{Definition}

As with the case for vector fields, we say that $\alpha\in\Omega(M)$ \emph{vanishes along $N$} whenever ${\alpha\in I_\Omega(N)}$, and that $\alpha,\beta\in\Omega(M)$ are \emph{equal along $N$} whenever $\alpha-\beta\in I_\Omega(N)$.

It is easy to show that $I_\Omega(N)\subset\Omega(M)$ is an ideal with respect to the wedge product. We now show that it is closed under the exterior derivative.

\begin{Lemma}\label{lem:closure_of_I_Omega(N)}
	For all $v\in\X_N(M)$, the space $\tgtvanform$ is closed under $(i)$ $\iota_v$, $(ii)$ $\d$ and $(iii)$ $\L_v$.
\end{Lemma}

\begin{proof}
	Fix a homogeneous element $\alpha\in \tgtvanform$ of degree $k$, and suppose that $u_0,\dots,u_k\in\X(M)$ are tangent to $N$.

(i) The identity $\iota_{u_k}\cdots\iota_{u_2}\iota_v\,\alpha=0$ on $N$ is immediate.

(ii) We have
			\begin{align*}
				\d\alpha(u_0,\dots,u_k)
={}&	\sum_i (-1)^i u_i \,\alpha(u_0,\dots,\hat{u}_i,\dots,u_k)			\\
				&{}	+\sum_{i<j}(-1)^{i+j} \alpha([u_i,u_j],u_0,\dots,\hat{u}_i,\dots,\hat{u}_j,\dots,u_k),
			\end{align*}
			where each term in the first sum vanishes on $N$ as $\alpha(u_0,\dots,\hat{u}_i,\dots,u_k)\in I_N$
			and the same happens for each term in the second sum by Lemma \ref{lem:tangent_to_N_Lie_subalgebra}.

(iii) This follows from (i), (ii) and the identity $\L_v=\d\iota_v+\iota_v\d$.
\end{proof}

With the notation of Definition \ref{def:Lie_algebra_action_tangent_to_N}, consider a Lie algebra action $\g\curvearrowright(N\subset M)$.

\begin{Definition}\label{def:reducible_form}\samepage
	We call $\alpha\in\Omega(M)$ a \emph{reducible form} with respect to $\g\curvearrowright (N\subset M)$ when
	\begin{enumerate}\itemsep=0pt
		\item[(i)] $\L_\xi\alpha\in \tgtvanform$, and
		\item[(ii)] $\iota_\xi\alpha\in \tgtvanform$
	\end{enumerate}
	for every $\xi\in\g$. We write $\Omega(M)_{[N]}$ for the space of reducible forms.
\end{Definition}

\begin{Lemma}\label{lem:reducible_forms_closed_under_d}
	$\Omega(M)_{[N]}$ is closed under $\d$.
\end{Lemma}

\begin{proof}
	Fix a reducible form $\alpha\in\Omega^k(M)$. Since $\tgtvanform$ is closed under $\d$ by Lemma \ref{lem:closure_of_I_Omega(N)}, it follows that
	\[
		\L_\xi\d\alpha = \d\L_\xi\alpha
\qquad	\text{and} \qquad	
		\iota_\xi\d\alpha = \L_\xi\alpha - \d\iota_\xi\alpha
	\]
	both lie in $\tgtvanform$.
\end{proof}

\begin{Definition}\label{def:reducible_vector_field}
	We say that $v\in\X(M)$ is a \emph{reducible vector field} with respect to $\g\curvearrowright (N\subset M)$ when
	\begin{enumerate}[i.]\itemsep=0pt
		\item[(i)] $v$ is tangent to $N$, and
		\item[(ii)] $[v,\underline\xi]\in \fgmodule(M)+\vanvf$ for all ${\xi}\in\mathfrak g$,
	\end{enumerate}
	where $\fgmodule(M)$ denotes the $C^\infty(M)$-submodule of $\mathfrak{X}(M)$ generated by the fundamental vector fields $\underline\xi$ for  all $\xi\in\g$.
	We denote the space of reducible vector fields by $\X(M)_{[N]}$.
\end{Definition}

We note that condition (ii) of Definition \ref{def:reducible_vector_field} is equivalent to
\begin{itemize}\itemsep=0pt
	\item[(ii$^\prime$)] $[v,\underline\xi] \in \fgmodule(M)$ along $N$ for all $\xi\in\g$.
\end{itemize}

\begin{Remark}\label{rem:understanding_reducibility}
	More algebraically, $v$ is reducible when it preserves both $I_N$ and $\X_\g(M)+I_\X(N)$.
	The latter space may be understood as a generalization of the space of vector fields tangent to the fundamental distribution $\underline\g=\bigl\{\underline\xi_x\mid \xi\in\g,\,x\in M\bigr\}$ along $N$.
	Indeed, in the setting of free actions, $\Gamma(\underline\g)$ and $\X_\g(M)$ coincide.
\end{Remark}

\begin{Lemma}\label{lem:reducible_vector_fields_Lie_algebra}
	$\X(M)_{[N]}\subset\X(M)$ is a Lie subalgebra.
\end{Lemma}

\begin{proof}
	Let $u,v\in\X(M)_{[N]}$. Lemma \ref{lem:tangent_to_N_Lie_subalgebra} asserts that $[u,v]$ is tangent to $N$, and it remains only to establish condition (ii) of Definition \ref{lem:reducible_vector_fields_Lie_algebra}. Thus fix $\xi\in\g$ and consider the identity $[[u,v],\underline\xi] = [[u,\underline\xi],v] + [u,[v,\underline\xi]]$.
 Condition (ii) now follows if $u$ satisfies
	\[
		[u, \fgmodule(M)+\vanvf]\subset  \fgmodule(M)+\vanvf,
	\]
	and similarly for $v$. We proceed in two steps:
	\begin{itemize}\itemsep=0pt
		\item First we establish $[u, \fgmodule(M)]\subset  \fgmodule(M)+\vanvf$. Let $w=\sum_{i=1}^k f_i \underline \xi_i$ for some $\xi_i\in\g$. We have
		\[
			[u, w]	= \bigl[u, {\textstyle\sum_i} f_i \underline \xi_i\bigr]
				= \sum_{i=1}^k \bigl[u, f_i\underline\xi_i\bigr]
				= \sum_{i=1}^k f_i\bigl[u,{\underline\xi}_i\bigr]+(uf)\underline\xi_i.
		\]
		Terms of the form $f_i\bigl[u,{\underline\xi}_i\bigr]$ evidently lie in the $C^\infty(M)$-module $\fgmodule(M)+\vanvf$, while those of type $(uf)\underline\xi_i$ lie in $\fgmodule(M)$.

	\item Now we show $[u, \vanvf]\subset  \fgmodule(M)+\vanvf$. Here we apply the fact that $[\tgtvf,\vanvf]\subset \vanvf$. This follows as for any $f\in C^\infty(M)$, $u\in \tgtvf$, and $w\in \vanvf$, we have $[u,w]f=uwf-wuf\in u\,I_n-w\,C^\infty(M)\subset I_N$.
\hfill\qed
\end{itemize}
\renewcommand{\qed}{}
\end{proof}

\begin{Lemma}\label{lem:Xg+I_X_is_a_Lie_ideal}
	$\X_\g(M)+I_\X(N)\subset\X(M)_{[N]}$ is a Lie ideal.
\end{Lemma}

\begin{proof}
	Since Lemmas \ref{lem:tangent_to_N_Lie_subalgebra} and \ref{lem:reducible_vector_fields_Lie_algebra} imply that $I_\X(N)\subset\X(M)_{[N]}$ is a Lie ideal, it suffices to show that $\bigl[\X(M)_{[N]},\X_\g(M)\bigr]\subset \X_\g(M)+I_\X(N)$. Thus let $v\in\X(M)_{[N]}$, $f\in C^\infty(M)$, and $\xi\in\g$. We have
	\[
		[v,f\underline\xi] = f[v,\underline\xi] + (vf)\underline\xi \;\in \X_\g(M) + I_\X(N)
	\]
	since $[v,\underline\xi]\in\X_\g(M)+I_\X(N)$ by Definition \ref{lem:reducible_vector_fields_Lie_algebra}, since $\X_\g(M)+I_\X(N)$ is closed under multiplication by $C^\infty(M)$, and since $(vf)\underline\xi\in\X_\g(M)$. The result follows as this inclusion extends to arbitrary terms $\sum_if_i\underline\xi_i\in\X_\g(M)$ by linearity.
\end{proof}

\subsection{Reducible observables}

Fix an action $\g\curvearrowright(N\subset M)$ and suppose that the premultisymplectic form $\omega\in\Omega^{n+1}(M)$ is $\g$-invariant along $N$ in the sense that $\L_\xi\omega\in I_\Omega(N)$ for all $\xi\in\g$,

\begin{Definition}\label{def:reducibleObservables}
	We will say that $\alpha\in \Ham_\infty^{<0}(M,\omega)$ is a \emph{reducible observable} with respect to $\g\curvearrowright(N\subset M)$ whenever $\alpha$ is a reducible form, and that $(v,\alpha)\in\Ham_\infty^0(M,\omega)$ is a reducible observable whenever
	\begin{enumerate}[i.]\itemsep=0pt
		\item[(i)] $v$ is a reducible vector field, and
		\item[(ii)] $\alpha$ is a reducible form.
	\end{enumerate}
	We denote the space of reducible observables by $\Ham_\infty(M,\omega)_{[N]}$.
\end{Definition}

Even when the form is multisymplectic, conditions (i) and (ii) are independent.

\begin{Example}
	Let $M=\mathbb R^3$, $\omega=\d x\wedge \d y\wedge \d z$. Let $N=\{z=0\}$. Let $\mathfrak g$ be the $1$-dimensional abelian Lie algebra acting by $\xi=\partial_y$. The vector fields tangent to $N$ are $C^\infty(\mathbb R^3)$-generated by~$z\partial_z$, $\partial_x$, $\partial_y$.
	The form $\alpha=yz\,\d y$ is reducible, since $\L_{\partial_y}\alpha$ and $\iota_{\partial_y}\alpha$ are multiples of $z$. However, its Hamiltonian vector field $v_\alpha=y\partial_x$ is not reducible, since $[\partial_y,v_\alpha]=\partial_x$ which does not lie in~$\fgmodule(M)$.
\end{Example}

\begin{Lemma}
	$\Ham_\infty(M,\omega)_{[N]}$ is an $L_\infty$-subalgebra of $\Ham_\infty(M,\omega)$.
\end{Lemma}

\begin{proof}
	Lemma \ref{lem:reducible_forms_closed_under_d} implies that $\Ham_\infty(M,\omega)_{[N]}$ is closed under $\tilde{l}_1$.

Now suppose that $(v_1,\alpha_1),\dots,(v_k,\alpha_k)\in \Ham_\infty^0(M,\omega)_{[N]}$ for $k\geq 2$, that $u_{k+1},\dots,u_{n+1}\in\X_N(M)$, and that $\xi\in\g$. We have
\begin{gather*}
\iota_{u_{n+1}}\!\cdots\,\iota_{u_{k+1}}\L_\xi\iota_{v_k}\!\cdots\,\iota_{v_1}\omega
=	(\L_\xi\omega)(v_1,\dots,v_k,u_{k+1},\dots,u_{n+1})			             \\ \hphantom{\iota_{u_{n+1}}\!\cdots\,\iota_{u_{k+1}}\L_\xi\iota_{v_k}\!\cdots\,\iota_{v_1}\omega=}{}
+\omega\bigl([\underline\xi,v_1],v_2,\dots,v_k,u_{k+1},\dots,u_{n+1}\bigr)  +\cdots   \\  \hphantom{\iota_{u_{n+1}}\!\cdots\,\iota_{u_{k+1}}\L_\xi\iota_{v_k}\!\cdots\,\iota_{v_1}\omega=}{}
+\omega\bigl(v_1,\dots,[\underline\xi,v_k],u_{k+1},\dots,u_{n+1}\bigr)         \\ \hphantom{\iota_{u_{n+1}}\!\cdots\,\iota_{u_{k+1}}\L_\xi\iota_{v_k}\!\cdots\,\iota_{v_1}\omega}{}
=	(\L_\xi\omega)(v_1,\dots,v_k,u_{k+1},\dots,u_{n+1})		               	 \\ \hphantom{\iota_{u_{n+1}}\!\cdots\,\iota_{u_{k+1}}\L_\xi\iota_{v_k}\!\cdots\,\iota_{v_1}\omega=}{}
+(\d\alpha_2)\bigl([\underline\xi,v_1],v_3,\dots,v_k,u_{k+1},\dots,u_{n+1}\bigr)       \\ \hphantom{\iota_{u_{n+1}}\!\cdots\,\iota_{u_{k+1}}\L_\xi\iota_{v_k}\!\cdots\,\iota_{v_1}\omega=}{}
-(\d\alpha_1)\bigl([\underline\xi,v_2],v_3,\dots,v_k,u_{k+1},\dots,u_{n+1}\bigr)	     \\ \hphantom{\iota_{u_{n+1}}\!\cdots\,\iota_{u_{k+1}}\L_\xi\iota_{v_k}\!\cdots\,\iota_{v_1}\omega=}{}
-(\d\alpha_1)\bigl(v_2,[\underline\xi,v_3],\dots,v_k,u_{k+1},\dots,u_{n+1}\bigr)    -\cdots 	 \\ \hphantom{\iota_{u_{n+1}}\!\cdots\,\iota_{u_{k+1}}\L_\xi\iota_{v_k}\!\cdots\,\iota_{v_1}\omega=}{}
-(\d\alpha_1)\bigl(v_2,\dots,[\underline\xi,v_k],u_{k+1},\dots,u_{n+1}\bigr)
\in I_N,
	\end{gather*}
	since $\L_\xi\omega\in I_\Omega(N)$, since $[\underline\xi,v_i]\subset\fgmodule(M) + \vanvf$, and since $\d\alpha_i\in\Omega(M)_{[N]}$ by Lemma \ref{lem:reducible_forms_closed_under_d}. Similarly,
	\[
		\iota_{u_{n+1}}\!\cdots\,\iota_{u_{k+2}}\,\iota_\xi\,\iota_{v_k}\!\cdots\,\iota_{v_1}\omega
			= - \iota_{u_{n+1}}\!\cdots\,\iota_{u_{k+2}}\,\iota_\xi\,\iota_{v_k}\!\cdots\,\iota_{v_2}\d\alpha_1
			\in I_N
	\]
	since $\d\alpha_1\in\Omega(M)_{[N]}$. Thus, $\Ham_\infty(M,\omega)_{[N]}$ is closed under $\tilde{l}_k$ for $k\geq 3$.

	Finally, as $[v_1,v_2]$ is a Hamiltonian vector field for $\iota_{v_1\wedge v_2}\omega$, and as $[v_1,v_2]$ is reducible by Lemma \ref{lem:reducible_vector_fields_Lie_algebra}, we conclude that $\Ham_\infty(M,\omega)_{[N]}$ is closed under $\tilde{l}_2$.
\end{proof}

\subsection[L\_infty-reduction of the space of observables]{$\boldsymbol{L_\infty}$-reduction of the space of observables}

Again suppose that $\g\curvearrowright(N\subset M)$ and that the premultisymplectic form $\omega\in\Omega^{n+1}(M)$ is $\g$-invariant along $N$.

\begin{Definition}\label{def:VanishIdeal}
	The \emph{vanishing observable ideal} $I_{\Ham_\infty}(N)\subset \Ham_\infty(M,\omega)_{[N]}$ is linearly generated by the subspace
	\[
		I_\Omega(N)\cap\Omega^{<n-1}(M) \subset \Ham_\infty^{<0}(M,\omega)
	\]
	and all pairs $(v,\alpha)\in\Ham_\infty^0(M,\omega)_{[N]}$ for which
	\begin{enumerate}\itemsep=0pt
		\item[(i)] $v\in \fgmodule(M) + \vanvf$,
		\item[(ii)] $\alpha\in I_\Omega(N)$.
	\end{enumerate}
\end{Definition}

\begin{Lemma}\label{lem:Ham_infty_is_an_L-infty_ideal}
	$I_{\Ham_\infty}(N)$ is an $L_\infty$-ideal of $\Ham_\infty(M,\omega)_{[N]}$.
\end{Lemma}

\begin{proof}
	As Lemma \ref{lem:tangent_to_N_Lie_subalgebra} asserts that $\tgtvanform$ is closed under $\d$, it follows that $I_{\Ham_\infty}(N)$ is closed under $\tilde{l}_1$.

	Fix $k\geq 3$. If $\sigma_1,\dots,\sigma_k\in \Ham_\infty^0(M,\omega)_{[N]}$ with $\sigma_i=(v_i,\alpha_i)$, if $\sigma_1\in I_{\Ham_\infty}(N)$, and if $u_{k+1},\dots,u_{n+1}\in\X_N(M)$, then
	\[
		\iota_{u_{n+1}}\!\cdots\,\iota_{u_{k+1}}\tilde{l}_k(\sigma_1,\dots,\sigma_k)
			= {-}\varsigma(k)\,\iota_{u_{n+1}}\!\cdots\,\iota_{u_{k+1}}\iota_{v_k}\!\cdots\,\iota_{v_2}\,\d\alpha_1 \in I_N
	\]
	since $v_1,\dots,v_k\in \X_N(M)$ and since $\d\alpha_1\in \tgtvanform$ by Lemma \ref{lem:closure_of_I_Omega(N)}.

	Now take $k=2$. For any $\sigma_1,\sigma_2\in \Ham_\infty^0(M,\omega)$ with $\sigma_1\in I_{\Ham_\infty}(N)$, Lemma \ref{lem:Xg+I_X_is_a_Lie_ideal} provides $[v_1,v_2]\in \X_\g(M) + I_\X(N)$ from which $\tilde{l}_2(\sigma_1,\sigma_2)\in I_{\Ham_\infty}(N)$, as required.
\end{proof}

Thus we arrive at the following commutative diagram of graded vector spaces, where for a~graded vector space $V$ we write $V[k]$ for the $k$-th desuspension,
\[
	\begin{tikzcd}[column sep =1.0em]
		I_{\Ham_\infty}(N) \ar[r,hook,"L_\infty\text{-ideal}"] \ar[d,hook]&
		\Ham_\infty(M,\omega)_{[N]} \ar[r,hook,"L_\infty\text{-morph.}"] \ar[d,hook] &
		\Ham_\infty(M,\omega) \ar[d,hook]
		\\
		 (I_\Omega(N) ){[n{-}1]}\!\oplus\! (\fgmodule(M)\!+\!\vanvf ) \ar[r,hook]&
		 \bigl(\Omega(M)_{[N]} \bigr){[n{-}1]}\!\oplus\! \X(M)_{[N]} \ar[r,hook]&
		 (\Omega(M) ){[n{-}1]}\!\oplus\! \X(M).
	\end{tikzcd}
\]

We now present our main construction.

\begin{Definition}\label{def:LinftyReductions}
	The \emph{reduction} of $\Ham_\infty(M,\omega)$ with respect to $\g\curvearrowright (N\subset M)$ is the $L_\infty$-algebra
	\[
		\Ham_\infty(M,\omega)_N = \frac{\Ham_\infty(M,\omega)_{[N]}}{I_{\Ham_\infty}(N)}.
	\]
\end{Definition}

Consolidating the foregoing development, the graded vector space underlying the reduced $L_\infty$-algebra $\Ham_\infty(M,\omega)_N$ is explicitly given in degree $0$
by
\[
	\Ham_\infty^0(M,\omega)_N =
	\frac{
		\left\lbrace
			(\alpha,v) \in \Omega^{n{-}1}(M) \oplus\X(M)
		~\left\vert~
			\begin{array}{@{}l@{}}
				\iota_v \omega = -\d \alpha
				\\
				\iota_\xi \alpha \in I_{\Omega}(N)
				\\
				\L_\xi\alpha \in I_{\Omega}(N)
				\\
				\L_\xi v \in \fgmodule +\vanvf
				\quad\forall \xi \in \g
				\\
				v \in \X_N(M)
			\end{array}
		\right.
		\right\rbrace
	}{
		\left\lbrace
			(\alpha,v) \in \Omega^{n{-}1}(M)\oplus\X(M)
		~\left\vert~
			\begin{array}{@{}l@{}}
				\iota_v \omega = -\d \alpha
				\\
				\alpha \in I_{\Omega}(N)
				\\
				v \in \fgmodule +\vanvf
			\end{array}
		\right.
		\right\rbrace
	}.
\]

An analogous procedure applies in the setting of the associated $L_\infty$-algebra $L_\infty(M,\omega)$.

\begin{Definition}\label{def:associatedLinftyReductions}
	We define $L_\infty(M,\omega)_{[N]}\subset L_\infty(M,\omega)$ and $I_{L_\infty}(N)\subset L_\infty(M,\omega)$ to be the respective images of $\Ham_\infty(M,\omega)_{[N]}$ and $I_{\Ham_\infty}(N)$ under the natural projection
	\begin{align*}
		\Ham_\infty(M,\omega)	&\to		L_\infty(M,\omega),	\\
		(v,\alpha)		&\mapsto	\alpha,		\\
		\beta		&\mapsto	\beta
	\end{align*}
	for $|(v,\alpha)|=0$ and $|\beta|<0$. That is, $\alpha\in L_\infty^0(M,\omega)_{[N]}$ precisely when there exists a $v\in\X(M)$ for which $(v,\alpha)\in\Ham_\infty(M,\omega)_{[N]}$, and similarly for $I_{\Ham_\infty}(N)$. Straightforward adaptations of the results above show that $L_\infty(M,\omega)_{[N]}\subset L_\infty(M,\omega)$ is an $L_\infty$-subalgebra and that $I_{L_\infty}(N)\subset L_\infty(M,\omega)_{[N]}$ is an $L_\infty$-ideal. The \emph{reduction} of $L_\infty(M,\omega)$ with respect to $\g\curvearrowright (N\subset M)$ is the $L_\infty$-algebra
	\[
		L_\infty(M,\omega)_N = \frac{L_\infty(M,\omega)_{[N]}}{I_{L_\infty}(N)}.
	\]
\end{Definition}

\begin{Remark}\label{rem:reduced_spaces_of_vfs_and_forms}
	We note that Lemmas \ref{lem:closure_of_I_Omega(N)} and \ref{lem:reducible_forms_closed_under_d} yield a reduced complex of forms
	\[
		\Omega(M)_N = \frac{\Omega(M)_{[N]}}{I_\Omega(N)},
	\]
	while Lemmas \ref{lem:reducible_vector_fields_Lie_algebra} and \ref{lem:Xg+I_X_is_a_Lie_ideal} provide a reduced space of vector fields
	\[
		\X(M)_N = \frac{\X(M)_{[N]}}{\X_\g(M)+I_\X(N)}.
	\]
\end{Remark}

\begin{Remark}\label{rem:Ham_N_and_L_infty_N_are_distinct}
	When $\omega$ is multisymplectic, the strict $L_\infty$-morphism
	\begin{align*}
		L_\infty(M,\omega)	&\to		\Ham_\infty(M,\omega),	\\
		\alpha		&\mapsto	(v_\alpha,\alpha),	\\
		\beta		&\mapsto	\beta
	\end{align*}
	for $|\alpha|=0$ and $|\beta|<0$, is inverse to the natural projection of Definition \ref{def:associatedLinftyReductions} and exhibits the isomorphism $L_\infty(M,\omega)\cong\Ham_\infty(M,\omega)$, and in particular induces an identification $L_\infty(M,\omega)_{[N]}\cong \Ham_\infty(M,\omega)_{[N]}$.

	But even if $\omega$ is multisymplectic, it is not generally true that $L_\infty(M,\omega)_N \cong \Ham_\infty(M,\omega)_N$. This results follow from the fact that the differential form components of distinct elements $(v_\alpha,\alpha),(v_\beta,\beta)\in\Ham_{\infty}^0(M,\omega)_{[N]}$ may be identified in $\Ham_\infty(M,\omega)_N$, while the vector field components remain distinct. That is, $[\alpha]=[\beta]\in\Omega(M)_N$ while $[v_\alpha]\neq[v_\beta]\in\X(M)_N$, so that $[\alpha]=[\beta]\in L_\infty(M,\omega)_N$ while $[(v_\alpha,\alpha)]\neq [(v_\beta,\beta)]\in \Ham_\infty(M,\omega)_N$.
	
	We return to this topic in Remark \ref{rem:associated_geometric_reduction_map}.
\end{Remark}

\begin{Remark}
	When $\g=0$ and $N\subset M$ is any subset, our reduction scheme defines a \emph{restriction} of $\Ham_\infty(M,\omega)$ to an arbitrary subset $N\subset M$. Indeed, in this setting the conditions that $\g\curvearrowright(N\subset M)$ and that $\omega$ is reducible are trivially satisfied. We note that the details of this multisymplectic restriction are prefigured in \cite[Lemma 3.2]{ShabaziZambon16}.

	As another boundary case, when $N=M$ and $\g\curvearrowright M$ is any action with $\underline\g\subset\ker\omega$, where we recall that $\underline\g=\bigl\{\underline\xi_x\mid \xi\in\g,\,x\in M\bigr\}$, we obtain a \emph{quotient} of $\Ham_\infty(M,\omega)$. Note that this boundary case is only interesting in the strictly premultisymplectic setting, i.e., when $\ker\omega$ is non-trivial.
\end{Remark}

\subsection[L\_infty-reduction with respect to a covariant moment map]{$\boldsymbol{L_\infty}$-reduction with respect to a covariant moment map}\label{Ssec:LinftyReductionMomentMap}

In Section~\ref{Ssec:MultiRed}, we considered the geometric reduction of an $n$-plectic manifold $(M,\omega)$ with respect to a Hamiltonian action $G\curvearrowright (M,\omega)$ with moment map $\mu\in\Omega^{n-1}(M,\g^*)$ and a closed form $\phi\in\Omega^{n-1}(M,\g^*)$. Our aim in this subsection is to apply our observable reduction scheme to this setting.

We approach the problem in slightly greater generality. Fix a pre-$n$-plectic manifold $(M,\omega)$, a premultisymplectic action $\g\curvearrowright(M,\omega)$, a form $\mu\in\Omega^{n-1}(M,\g^*)$ satisfying
\begin{enumerate}[i.]\itemsep=0pt
	\item[(i)] $\d\mu_\xi = -\iota_\xi\omega$,
	\item[(ii)] $\L_\xi\mu_\zeta = \mu_{[\xi,\zeta]}$
\end{enumerate}
for all $\xi,\zeta\in\g$, and a closed form $\phi\in\Omega^{n-1}(M,\g^*)$. Consider the action $\g\curvearrowright\Omega^{n-1}(M,\g^*)$ given by
\[
	\xi\cdot \alpha = \L_\xi\alpha + \Ad_\xi^*\alpha
\]
and write
\[
	\g_\phi = \bigl\{\xi\in\g \mid  \L_\xi\phi_\zeta = \phi_{[\xi,\zeta]} \ \forall\zeta\in\g\bigr\}
\]
for the isotropy subalgebra of $\phi$.

\begin{Definition}\label{def:observable_reduction_wrt_moment_map}
	The \emph{reduction} of $\mathrm{Ham}_\infty(M,\omega)$ with respect to the Hamiltonian action $\g\curvearrowright(M,\omega)$, moment map $\mu\in\Omega(M,\g^*)$, and level $\phi\in\Omega_{\mathrm{cl}}(M,\g^*)$, is the reduction of $\Ham_\infty(M,\omega)$ with respect to $\g_\phi\action\bigl(\mu^{-1}(\phi)\subset M\bigr)$. We write $\Ham_\infty(M,\omega)_\phi = \Ham_\infty(M,\omega)_{\mu^{-1}(\phi)}$.
\end{Definition}

Our task is now to show that this construction is well defined.

\begin{Lemma}\label{lem:isotropy_action_is_tangent_to_phi_level_set}
	The action $\g_\phi\curvearrowright M$ is tangent to $\mu^{-1}(\phi)$.
\end{Lemma}

\begin{proof}
	When $\g_\phi\curvearrowright M$ is induced by the action of a connected Lie group $G_\phi\curvearrowright M$, the equivariance condition
	\[
		\L_\xi(\mu-\phi)_\zeta = (\mu-\phi)_{[\xi,\zeta]}
	\]
	for all $\xi,\zeta\in\g_\phi$ implies
	\[
		\bigl(g^{-1}\bigr)^*(\mu-\phi)_\zeta = (\mu-\phi)_{\Ad_g\zeta}
	\]
	for all $g\in G_\phi$. In particular, $G_\phi$ preserves $\mu^{-1}(\phi)=\{x\in M\mid (\mu-\phi)_x = 0\}$ and the result follows by Remark \ref{rem:group_preservation_is_sufficient_for_tangency}.

	When there is no group action, the above argument still holds since Lie's second fundamental theorem (see, e.g., \cite[Theorem~20.22]{Lee12}) provides that we may always integrate $\g_\alpha\curvearrowright M$ to a~local Lie group action, and this is sufficient for our purposes here.
\end{proof}

\begin{Lemma}\label{lem:omega_is_reducible_over_phi_level_set}
	The presymplectic form $\omega$ is reducible with respect to $\g_\phi\curvearrowright\bigl(\mu^{-1}(\phi)\subset M\bigr)$.
\end{Lemma}

\begin{proof}
	Fix $\xi\in\g_\alpha$. Since $\g\curvearrowright M$ is premultisymplectic, the equality $\L_\xi\omega=0$ trivially implies $\L_\xi\omega\in I_\Omega\bigl(\mu^{-1}(\phi)\bigr)$.

	It remains to show that $\iota_\xi\omega\in I_\Omega\bigl(\mu^{-1}(\phi)\bigr)$. By construction, $\mu-\phi=0$ at every point of~$\mu^{-1}(\phi)$ and thus, in particular, $(\mu-\phi)_\xi\in I_\Omega(N)$. As Lemma \ref{lem:closure_of_I_Omega(N)} provides that $I_\Omega\bigl(\mu^{-1}(\phi)\bigr)$ is closed under $\d$, it follows that
	\[
		\iota_\xi\omega = \d\mu_\xi = \d(\mu-\phi)_\xi \in I_\Omega\bigl(\mu^{-1}(\phi)\bigr).
\tag*{\qed}
\]
\renewcommand{\qed}{}
\end{proof}

In particular, it follows that $\L_\xi\omega\in I_\Omega(N)$ for all $\xi\in\g$, and we conclude from Lemmas \ref{lem:isotropy_action_is_tangent_to_phi_level_set} and \ref{lem:omega_is_reducible_over_phi_level_set} that Definition \ref{def:observable_reduction_wrt_moment_map} is indeed well defined.

\begin{Proposition}
	For every $\xi\in\g_\phi$, we have $\mu_\xi-\phi_\xi\in I_{\Omega}\bigl(\mu^{-1}(\phi)\bigr)$ and $\bigl(\underline\xi,\mu_\xi-\phi_\xi\bigr)\in \smash{I_{\Ham_\infty}\bigl(\mu^{-1}(\phi)\bigr)}$.
\end{Proposition}

\begin{proof}
	First note that $\underline\xi$ and $\mu_\xi-\phi_\xi$ comprise a Hamiltonian pair as $\d(\mu-\phi)_\xi = \d\mu_\xi = \iota_\xi\omega$.

	We showed that $(\mu-\phi)_\xi\in I_\Omega\bigl(\mu^{-1}(\phi)\bigr)$ in the proof of Lemma \ref{lem:omega_is_reducible_over_phi_level_set}, and that $\underline\xi$ is tangent to $\mu^{-1}(\phi)$ in Lemma \ref{lem:isotropy_action_is_tangent_to_phi_level_set}. The result follows as $\underline\xi$ is clearly an element of \smash{$\X_\g(M)+I_\X\bigl(\mu^{-1}(\phi)\bigr)$}.
\end{proof}

\subsection[Comparing geometric and L\_infty-reduction]{Comparing geometric and $\boldsymbol{L_\infty}$-reduction}

The aim of this section is to compare the $L_\infty$-reduction procedure to the multisymplectic reduction scheme of Theorem \ref{thm:CaseyReduction}. In particular, we exhibit a natural inclusion of the $L_\infty$-reduction $\Ham_\infty(M,\omega)_\phi$ into the $L_\infty$-algebra of observables $\Ham_\infty(M_\phi,\omega_\phi)$ associated to the reduced space $(M_\phi,\omega_\phi)$.

In fact, we will work in greater generality than the setting of Theorem \ref{thm:CaseyReduction}. For the extent of this subsection, fix a pre-$n$-plectic manifold $(M,\omega),$ a connected Lie group $G$, a smooth action $G\curvearrowright M$, and a $G$-invariant closed embedded submanifold $j\colon N\hookrightarrow M$ such that $\omega$ is reducible and $G\curvearrowright N$ is free and proper.

The inclusion $j\colon N\hookrightarrow M$ and the projection $\pi\colon N\twoheadrightarrow N/G$ each induce maps on the de Rham complex as depicted below:
\[
	\begin{tikzcd}
		N \ar[r,hook,"j"] \ar[d,two heads,"\pi"'] & M,  & \Omega(N) & \Omega(M). \ar[l,two heads,"j^\ast"']
		\\
		N/G &   & \Omega(N/G) \ar[u,hook,"\pi^\ast"]
	\end{tikzcd}
\]

We are now in a position to recast the auxiliary spaces of Figure \ref{fig:intermediate_spaces} in a more geometric fashion.

\begin{Lemma}\label{lem:RedIngredientsRegularCase}
	We have
	\begin{enumerate}\itemsep=0pt
		\item[$(i)$] $\tgtvf= \{v\in\X(M) \mid v|_N\in\mathfrak X(N)  \}$,
		\item[$(ii)$] $\vanvf= \{v\in\X(M) \mid v|_N=0  \}$,
		\item[$(iii)$] $\tgtvanform=\{\alpha\in\Omega(M)\mid j^*\alpha=0\}$,
		\item[$(iv)$] $\fgmodule(N)=\Gamma\bigl(\underline{\g}\bigr)$,
	\end{enumerate}
	where $v\vert_N = v \circ j \in \Gamma(j^\ast TM)$ is the restriction of $v\in\X(M)$ to $N$, and where
	\[
		\underline\g=\bigl\{\underline\xi_x\mid \xi\in\g,x\in N\bigr\}\subset TN
	\]
	is the fundamental distribution.
\end{Lemma}

\begin{proof}
	(i)--(iii) As $N\subset M$ is an embedded submanifold, it suffices to consider $\R^{\dim N}\subset\R^{\dim M}$ and the result now follows by a direct computation.

	(iv) Since the action is free, the fundamental vector fields associated to a basis of $\mathfrak g$ are linearly independent at every point. This implies that $\fgmodule$ is a regular foliation and hence generated by its associated distribution.
\end{proof}

Recall that $\alpha\in\Omega(N)$ is said to be \emph{$\g$-basic} when it is both $\g$-invariant and $\g$-horizontal. That~is,
\[
	\Omega_{\mathrm{bas}}(N) =
	\bigl\{
	\alpha \in \Omega(N)
	\mid
		\iota_\xi\alpha=0, \, 	\L_\xi\alpha=0 \, \forall \xi \in \g
	\bigr\}.
\]

We now show that, since $G$ is connected, the $\g$-basic forms on $N$ are precisely those induced by $\pi\colon N\to N/G$.

\begin{Lemma}\label{lem:BasicForms}
	For free and proper actions $G\curvearrowright N$, there is an isomorphism of de Rham complexes $\pi^*\colon \Omega(N/G)\xrightarrow{\sim}\Omega_{\mathrm{bas}}(N)$.
\end{Lemma}

\begin{proof}
	An application of \cite[Theorem~31.12]{Tu17} to the surjective submersion $\pi\colon N \to N/G$ provided by the quotient manifold theorem, yields $\operatorname{im}\pi^*=\Omega_{\mathrm{bas}}(N)$.
\end{proof}

In particular, given $\alpha\in\Omega(M)_{[N]}$, we have that $j^\ast \alpha$ is $\g$-basic, and hence there corresponds a unique $\alpha_N\in\Omega(N/G)$ for which $j^*\alpha=\pi^*\alpha_N$.
Likewise, the connectedness of $G$ provides that reducible vector fields are projectable along $\pi\colon N\to N/G$, and thus we may assign to each $v\in\X(M)_{[N]}$ a unique $v_N=\pi_*(v|_N)\in\X(N/G)$.

\begin{Definition}\label{def:geometric_reduction_map}
	Put $M_N=N/G$. The \emph{geometric reduction map} on forms is
	\begin{displaymath}
		\morphism{r_N}
		{\Omega(M)_{[N]}	}{\Omega(M_N),}
		{\alpha}{\alpha_N}
	\end{displaymath}
	and on vector fields is
	\begin{displaymath}
		\morphism{r_N}
		{\X(M)_{[N]}	}{\X(M_N),}
		{v}{v_N.}
	\end{displaymath}
\end{Definition}

\begin{Remark}
	As the equality $j^*\alpha=\pi^*\alpha_N$ implies $j^*\d\alpha = \pi^*\d(\alpha_N)$, it follows that $\d(\alpha_N)=(\d\alpha)_N$, and consequently that $r_N\colon \alpha\mapsto\alpha_N$ is a map of chain complexes. Similarly, as $v_N f = v|_N (\pi^*f)$ for all $f\in C^\infty(M_N)$, the map $r_N\colon v\mapsto \pi_*(v|_N)$ is easily seen to be a Lie algebra homomorphism.
\end{Remark}

As $\pi^*$ preserves $\d$, the form $\omega_N\in\Omega^{n+1}(M_N)$ is closed, and the following definition arises naturally. This will be our generalization of the geometric reduced spaces $(M_\phi,\omega_\phi)$ of Theorem~\ref{thm:CaseyReduction}.

\begin{Definition}\label{def:geometric_reduction}
	The \emph{geometric reduction} of the $n$-plectic manifold $(M,\omega)$ is the premultisymplectic manifold $(M_N,\omega_N)$.
\end{Definition}

Our aim now is to construct a strict $L_\infty$-morphism $r_N\colon \Ham_\infty(M,\omega)_{[N]}\to \Ham_\infty(M_N,\omega_N)$ on the level of observables in such a way that $\ker\,r_N=I_{\Ham_\infty}(N)$, so that $r_N$ descends to an inclusion
\[
	\bar{r}_N\colon \ \Ham_\infty(M,\omega)_N \hookrightarrow \Ham_\infty(M_N,\omega_N).
\]

\begin{Lemma}\label{lem:SesReducibleForms}
	There is a short exact sequence of chain complexes
	\[
		0\to I_\Omega(N)\hookrightarrow\Omega(M)_{[N]}\xrightarrow{r_N}  \Omega(M_N)\to 0.
	\]
\end{Lemma}

\begin{proof}
	If $\alpha\in\Omega(M)_{[N]}$, then $\alpha_N = 0$ precisely when $j^*\alpha=\pi^*\alpha_N=0$ by Lemma \ref{lem:BasicForms}, and $j^*\alpha=0$ precisely when $\alpha\in I_\Omega(N)$ by Lemma \ref{lem:RedIngredientsRegularCase}. Thus $\ker\,r_N= I_\Omega(N)$.

	The surjectivity of $r_N$ is follows from the isomorphism $\pi^*\colon \Omega(N/G)\xrightarrow{\sim}\Omega_{\mathrm{bas}}(M)$ and the surjectivity of $j^*\colon \Omega(M)\twoheadrightarrow\Omega(N)$.
\end{proof}

\begin{Lemma}\label{lem:SesReducibleFields}
	There is a short exact sequence of Lie algebras
	\[
		0\to\fgmodule(M)+\vanvf \hookrightarrow \X(M)_{[N]}\xrightarrow{r_N}\mathfrak X(M_N)\to 0.
	\]
\end{Lemma}

\begin{proof}
	Lemma \ref{lem:RedIngredientsRegularCase} yields
	\[
		\ker r_N = \{v\in\X(M) \mid v|_N \mathrm{~is~vertical~w.r.t.~}N\to N/G\} = \fgmodule(M)+\vanvf
	\]
	and surjectivity obtains as any vector field on $N/G$ may be lifted to $N$ and extended to $M$.
\end{proof}

We now show that the reduction maps $r_N$ of Lemmas \ref{lem:SesReducibleForms} and \ref{lem:SesReducibleFields} preserve Hamiltonian pairs.

\begin{Lemma}
	If $(\alpha,v)$ is a reducible Hamiltonian pair on $(M,\omega)$, then $(\alpha_N,v_N)$ is a Hamiltonian pair on $(M_N,\omega_N)$.
\end{Lemma}

\begin{proof}
	This follows from the identity
	\[
		\pi^*(\iota_{v_N}\omega_N) = \iota_{v|_N} j^* \omega = j^*(\iota_v\omega) = j^*(-\d\alpha) = \pi^*(-\d\alpha_N)
	\]
	and the injectivity of $\pi^*\colon \Omega(N/G)\xrightarrow{\sim}\Omega_{\mathrm{bas}}(N)$.
\end{proof}

We thus obtain our desired geometric reduction map on observables.

\begin{Theorem}\label{thm:MultiSymRed_ManifoldvsObservables}
	The geometric reduction map
	\begin{align*}
		r_N\colon \	\Ham_\infty(M,\omega)_{[N]}	&\to		\Ham_\infty(M_N,\omega_N),		\\
			(v,\alpha) 			&\mapsto	 (v_N,\alpha_N)		,	\\
			\alpha 			&\mapsto	 \alpha_N
	\end{align*}
	is a strict $L_\infty$-morphism with kernel $I_{\Ham_\infty}(N)$. In particular, there is a natural inclusion of $L_\infty$-algebras
	\[
		\Ham_\infty(M,\omega)_N = \frac{\Ham_\infty(M,\omega)_{[N]}}{I_{\Ham_\infty}(N)} \xhookrightarrow{\;\bar{r}_N\;} \Ham_\infty(M_N,\omega_N).
	\]
\end{Theorem}

\begin{proof}
	Lemmas \ref{lem:SesReducibleForms} and \ref{lem:SesReducibleFields} together yield $\ker r_N=I_{\Ham_\infty}(N)$, and it remains to show that $r_N$ is a strict $L_\infty$-morphism.

	For $k\geq 3$ and $\sigma^1,\dots,\sigma^k \in \Ham_\infty^0(M,\omega)_{[N]}$, with $\sigma^i=\bigl(\alpha^i,v^i\bigr)$, we have
	\begin{align*}
		\pi^* \tilde{l}_k\bigl(\sigma_N^1,\dots,\sigma_N^k\bigr)
			&=	\varsigma(k)  \iota_{v_N^k}\cdots\iota_{v_N^1} \omega_N		\\
		&=	\varsigma(k)  \iota_{v^k|_N}\cdots\iota_{v^1|_N} j^*\omega		\\
		&=	j^*(\varsigma(k)  \iota_{v^k}\cdots\iota_{v^1}\omega)		\\
		&=	\pi^* \tilde{l}_k\bigl(\sigma^1,\dots,\sigma^k\bigr)_N,
	\end{align*}
	from which we obtain
	\[
		\tilde{l}_k\bigl(\sigma_N^1,\dots,\sigma_N^k\bigr) = \tilde{l}_k\bigl(\sigma^1,\dots,\sigma^k\bigr)_N
	\]
	by the injectivity of $\pi^*\colon \Omega(N/G)\xrightarrow{\sim}\Omega_{\mathrm{bas}}(N)$. The case $k=2$ follows similarly, with the additional observation that $\bigl[v_N^1,v_N^2\bigr] = \bigl[v^1,v^2\bigr]_N$.

	For $k=1$ and $\alpha\in\Ham_\infty^{<0}(M,\omega)_{[N]}$, we have
	\[
		\pi^*\tilde{l}_1(\alpha_N) = \d j^*\alpha = j^*\d\alpha = \pi^*\tilde{l}_1(\alpha)_N
	\]
	and we conclude that $r_N$ is a strict $L_\infty$-morphism.
\end{proof}

\begin{Remark}
	In this regular setting, we may \emph{define} the reduced space $\Ham_\infty(M,\omega)_N$ to be the quotient $\Ham_\infty(M,\omega)_{[N]} / \ker r_N$ or, equivalently, to be the image of $r_N$ in $\Ham_\infty(M_N,\omega_N)$. The merit of Definition \ref{def:LinftyReductions} is that it applies even in the singular case in which the geometric reduction map $r_N$ of smooth manifolds fails to exist.
\end{Remark}

We will call $\bar{r}_N\colon \Ham_\infty(M,\omega)_N \hookrightarrow \Ham_\infty(M_N,\omega_N)$ the \emph{canonical embedding}. The following corollary is immediate.

\begin{Corollary}
	Let $(M,\omega,G,\mu)$ be a $n$-plectic Hamiltonian $G$-space with covariant moment map $\mu$, let $\phi\in\Omega^{n-1}(M,\g^*)$ be a closed form, and let $M_\phi= \mu^{-1}(\phi)/G_\phi$. If $\mu^{-1}(\phi)\subset M$ is an embedded submanifold, and if $G$ acts freely on $\mu^{-1}(\phi)$, then there is a natural inclusion
	\begin{align*}
		\bar{r}_\phi\colon \	\Ham_\infty(M,\omega)_\phi	&\xhookrightarrow{\hspace{10pt}}	\Ham_\infty\bigl(M_\phi,\omega_\phi\bigr),	\\
				[\sigma]		&\longmapsto				\sigma_\phi
	\end{align*}
	of $L_\infty$-algebras of observables.
\end{Corollary}

\begin{Remark}\label{rem:associated_geometric_reduction_map}
	A similar argument to that of Theorem \ref{thm:MultiSymRed_ManifoldvsObservables} yields a natural inclusion
	\[
		\bar{r}_N \colon \ L_\infty(M,\omega)_N \to L_\infty(M_N,\omega_N).
	\]
	Recalling the discussion of Remark \ref{rem:Ham_N_and_L_infty_N_are_distinct}, we see that when \emph{both} $(M,\omega)$ and $(M_N,\omega_N)$ are multisymplectic, the strict $L_\infty$-morphism
	\begin{align*}
		\bar{r}_N[L_\infty(M,\omega)_N]		&\xrightarrow{\;\sim\;}	\bar{r}_N[\Ham_\infty(M,\omega)_N],	\\
		\alpha_N					&\xmapsto{\;\phantom{\sim}\;}		\bigl(v_{\alpha_N},\alpha_N\bigr),	\\
		\beta_N					&\xmapsto{\;\phantom{\sim}\;}		\beta_N
	\end{align*}
	for $|\alpha_N|=0$ and $|\beta_N|<0$, provides the natural identification $L_\infty(M,\omega)_N\cong\Ham_\infty(M,\omega)_N$.
\end{Remark}

Lemmas \ref{lem:SesReducibleForms} and \ref{lem:SesReducibleFields} yield isomorphisms
\[
	\bar{r}_N\colon \ \Omega(M)_N \xrightarrow{\sim} \Omega(M_N)
\]
and
\[
	\bar{r}_N\colon \ \X(M)_N \xrightarrow{\sim} \X(M_N),
\]
 respectively, where $\Omega(M)_N$ and $\X(M)_N$ are the reduced spaces of Remark~\ref{rem:reduced_spaces_of_vfs_and_forms}.

As no such isomorphism is guaranteed on the level of observables, it is natural to inquire into the relation between $L_\infty$-reduction and the \emph{association of observables} functor
\[
	\Ham_\infty \colon \ (M,\omega) \mapsto \Ham_\infty(M,\omega).
\]
This question motivates the following definition, which heuristically measures the extent to which these procedures fail to commute.

\begin{Definition}\label{def:residue_defect}
	The \emph{residue defect} $[\Ham_\infty,R]$ of $G\curvearrowright(M,\omega)$ with respect to $N\subset M$ is the cokernel
	\[
		\operatorname{coker} r_N = \frac{\Ham_\infty(M_N,\omega_N)}{\operatorname{im} r_N}
	\]
	of the geometric reduction map $r_N\colon \Ham_\infty(M,\omega)_{[N]}\to\Ham_\infty(M_N,\omega_N)$.
\end{Definition}

As Lemmas \ref{lem:SesReducibleForms} and \ref{lem:SesReducibleFields} ensure that $r_N^{<0}$ surjects onto $\Ham_\infty^{<0}(M_N,\omega_N)$, we may identify the residue defect with the component $\operatorname{coker}r_N^0 = \Ham_\infty^0(M_N,\omega_N) / \operatorname{im}r_N^0$ on Hamiltonian pairs.

The residue defect of $G\curvearrowright(M,\omega)$ with respect to $N\subset M$ measures the extent to which the reduced space $(M_N,\omega_N)$ exhibits observables that are not induced from $(M,\omega)$.

We will establish in Theorem \ref{thm:AllReductionsComparison-sympleticcase} that in the suitably regular symplectic setting $\operatorname{coker} r_N=0$. In the $n$-plectic case with $n>1$, the determination of it appears to be highly nontrivial.

\section{Applications and examples}\label{Sec:Examples}

In this section, we investigate our reduction formalism in the setting of various natural examples.

\subsection[L\_infty-reduction in the 1-plectic case]{$\boldsymbol{L_\infty}$-reduction in the 1-plectic case}\label{Ssec:LinftyRed-SympCase}

In this subsection, we discuss how definition the $L_\infty$-reduction scheme manifests in the $1$-plectic, i.e., symplectic, setting.

Let $(M,\omega)$ be a symplectic manifold, be $N\subset M$ be a closed subset, and consider a Lie algebra action $\mathfrak{g}\action (N\subset M)$.
The auxiliary spaces of the $L_\infty$-reduction scheme are
\begin{equation}\label{diagram:symplectic-ingredients}
	\begin{tikzcd}[
		ampersand replacement = \&
		,row sep = 0ex
		,/tikz/column 2/.append style={anchor=base west}
		]
		L_\infty(M,\omega) \ar[r,phantom,"="] \&
		C^\infty(M),
		\\
		L_\infty(M,\omega)_{[N]} \ar[r,phantom,pos=0.42,"="] \ar[u,hook]\&
		\left\lbrace
		f\in C^\infty(M)
		\,\left\vert \,
		\begin{array}{@{}l r l@{}}
			\underline\xi\hspace{1pt} f \in I_N &&\forall \xi \in \mathfrak{g}\\
			v_f h\in I_N && \forall h \in I_N
			\\
			{[v_f,\underline{\xi}]} \in \fgmodule(M)  & \text{along N},
			& \forall \xi \in \g
		\end{array}
		\right\rbrace,\right.
		\\
		I_{L_\infty}(N) \ar[r,phantom,pos=0.55,"="] \ar[u,hook]\&
		\left\lbrace
		f\in C^\infty(M)
		\,\left\vert \,
		\begin{array}{@{}l l@{}}
			\underline\xi\hspace{1pt} f \in I_N & \forall \xi \in \mathfrak{g}\\
			v_f \,h \in I_N & \forall h \in I_N
			\\
			f \in I_N\\
			v_f \subset \fgmodule(M) & \text{along N}
		\end{array}
		\right\rbrace.\right.
	\end{tikzcd}
\end{equation}

A particular feature of the symplectic setting is that the space $L_\infty(M,\omega)=C^\infty(M)$ naturally possesses the structure of an associative algebra with respect to which the $l_2$ bracket is a Poisson structure. We now show that this structure descends to the reduced space $L_\infty(M,\omega)_N$.

\begin{Theorem}
	If $\g\curvearrowright M$ is tangent to $N\subset M$, and if the symplectic structure $\omega\in\Omega^2(M)$ is reducible, then the reduction $L_\infty(M,\omega)_N$ inherits a natural Poisson algebra structure from $C^\infty(M,\omega)$.
\end{Theorem}

\begin{proof}
	Since $I_N$ is an associative ideal in $C^\infty(M)$, it follows that $I_{L_\infty}(N) = L_\infty(M,\omega)_{[N]}\cap I_N$ is an associative ideal in $L_\infty(M,\omega)_{[N]}$. Moreover, we have previously established that $I_{L_\infty}(N)\subset L_\infty(M,\omega)_{[N]}$ is an $L_\infty$-ideal (see Definition \ref{def:associatedLinftyReductions} and Lemma \ref{lem:Ham_infty_is_an_L-infty_ideal}), and thus, in particular, an ideal with respect to the $l_2$ bracket. The result follows as these two facts together imply that $I_{L_\infty}(N)\subset L_\infty(M,\omega)_{[N]}$ is a Poisson ideal.
\end{proof}

Thus, in  the presence of the additional hypothesis that $\g\curvearrowright(N\subset M)$ and $\omega$ is a reducible $2$-form, diagram \eqref{diagram:symplectic-ingredients} is in the category of Poisson algebras.

\begin{Remark}\label{rem:comparingAlongN}
	Observe that if $v\in\X_\ham(M,\omega)$, then $v_f \in I_{\mathfrak{X}}(N)$ precisely when $\{f, C^\infty(M)\} = v_f\, C^\infty(M)\subset I_N$.
\end{Remark}

In order to better understand the algebra defined above, it is helpful to consider certain subclasses of well-behaved functions.
\begin{Definition}[first class function]\label{def:IclassFunc}
	The function $f\in C^\infty(M)$ is said to be \emph{first class} if its Poisson bracket with any constraint vanishes on $N$. That is, the set of first class functions is
	\[
	 F_N =
		\lbrace
 		f \in C^\infty(M)
 		\mid
		\lbrace f, I_N \rbrace \subset I_N
		\rbrace,
	\]
	where $I_N$ denotes the ideal of functions vanishing on $N$.
\end{Definition}

Equivalently, $F_N$ is the Lie algebra normalizer of the subspace $I_N\subset C^\infty(M)$.
By virtue of the Jacobi identity, it follows that $F_N$ is a Poisson subalgebra of $C^\infty(M,\omega)$.

\begin{Definition}\label{def:Casimir_function}
	We define the space of \emph{Casimir functions} along $N$ to be
	\begin{equation}\label{eq:calQ}
	 \mathcal{Q}=
		\lbrace
			f \in C^\infty(M)
		 \mid
			\{ f, C^\infty(M)\} \subset I_N
		\rbrace	
		.
	\end{equation}
\end{Definition}

That is, in light of Remark \ref{rem:comparingAlongN},
\[
	f \in \mathcal{Q} \iff v_f\in I_\X(N).
\]

\begin{Lemma}[properties of $\mathcal{Q}$]\label{lem:Qproperties}
	The subspace $\mathcal{Q}$ enjoys the following properties:
	\begin{enumerate}[i.]\itemsep=0pt
		\item[$(i)$] $\mathcal{Q}\subset C^\infty(M,\omega)$ is a Poisson subalgebra $($but not an associative or Lie ideal$)$,
		\item[$(ii)$] $\mathcal{Q}\subset F_N$ is a Lie ideal,
		\item[$(iii)$] $\mathcal{Q}\cap I_N$ is an associative ideal in $C^\infty(M)$ and a Poisson ideal in $F_N$.		
	\end{enumerate}
\end{Lemma}

\begin{proof}
	This follows from an easy verification that
	\[
		\{fh,C^\infty(M)\} = f\{h,C^\infty(M)\} + h\{f,C^\infty(M)\}
	\]
	and
	\[
		\{\{f,h\},C^\infty(M)\} = \{\{h,C^\infty(M)\},f\}- \{\{f,C^\infty(M)\},h\}
	\]
	lie in $I_N$ whenever $f\in\mathcal{Q}$ or $\mathcal{Q}\cap I_N$ and $h\in C^\infty(M)$, $F_N$ or $\mathcal{Q}$, depending on the statement to be proved.
\end{proof}

We now specialize to the case where the action $G\action (M,\omega)$ is Hamiltonian, with moment map $\mu\colon M\to \g^*$.  Consider as the constraints set the preimage of the ensuing moment map $N=\mu^{-1}(0)$.
Momenta yields a relevant subset of vanishing functions.

\begin{Definition}[momentum ideal]\label{def:momentumideal}
	The \emph{momentum ideal} is the associative ideal $I_\mu\subset C^\infty(M)$ generated by the momenta $\mu_{\xi}$ for any $\xi\in \mathfrak{g}$.
	Namely
	\[
		I_\mu
		=
		\langle \mu_\xi \rangle_{\xi \in \mathfrak{g}}^{\text{asso}}
		=
		\Biggl\{
				\sum_{i=1}^n f_i ~ \mu_{\xi_i}
			\,\Bigg|\,
				n\geq0,\, f_i \in C^\infty(M),\, \xi_i\in \mathfrak{g},\,  1 \leq i \leq n
		\Biggr\}
		.
	\]
\end{Definition}

\begin{Theorem}\label{thm:LinftyRedSymp}
	Fix a symplectic Hamiltonian action $G\action (M,\omega)$ with moment map $\mu\colon M\to\g^*$ and put $N=\mu^{-1}(0)$.
	Then
	\begin{gather*}
			L_\infty(M,\omega)_{[N]} = F_N \cap \mathcal{Q}', \\
			I_{L_\infty}{(N)} = I_\mu + (I_N \cap \mathcal{Q} ),
	\end{gather*}
	where $I_\mu$, $F_N$ are given by Definitions~{\rm \ref{def:momentumideal}} and~{\rm \ref{def:IclassFunc}}, respectively, $\mathcal{N}\bigl(I_\mu\bigr)$ is the Lie algebra normalizer of $I_\mu$ inside of $C^\infty(M)$ $($see equation~\eqref{eq:ImuNorm}$)$, $\mathcal{Q}$ is given in equation~\eqref{eq:calQ}, and
	\begin{equation}
	\label{eq:calQprime}
	 \mathcal{Q}'=
		\lbrace
			f \in C^\infty(M)
		\mid
			\{ f, \mu_\xi\}\in I_\mu + \mathcal{Q}\  \forall \xi \in \g
		\rbrace		.			
	\end{equation}
	Furthermore, the reduced Poisson algebra is
\begin{align}\label{eq:1plred}
		\frac{L_\infty(M,\omega)_{[N]}}{I_{L_\infty}{(N)}} = \frac{F_N\cap\mathcal{Q}'}{I_\mu+I_N\cap \mathcal{Q}}.
\end{align}
\end{Theorem}

\begin{proof}
	Consider $f\in L_\infty(M,\omega)_{[N]}$.
	By the very definition of Poisson bracket, the first condition (see diagram \eqref{diagram:symplectic-ingredients}) in the definition of $L_\infty(M,\omega)_{[N]}$ reads as $\L_\xi f = \{\mu_\xi,f\} \in I_N$ for any $\xi \in \mathfrak{g}$.
	Hence the Jacobi identity implies that $\lbrace I_\mu, f\rbrace \subset I_N$.
	Similarly the second condition implies that $\L_{v_f} h = \lbrace f,h \rbrace \in I_N$ for any $h\in I_N$, hence $\{f,I_N\}\subset I_N$, i.e., $f$ is a first class function. In particular, since $I_\mu\subset I_N$, the latter condition implies the former one.
	The third condition implies that for any $\xi \in \g$ there exist two families $g_i \in C^\infty(M)$ and $\xi_i \in \g$ such that
	\[
		[v_f,\underline{\xi}] = \sum_i g_i \underline{\xi}_i \qquad \text{along $N$.}
	\]
	The left-hand side can be recast as
	$[v_f,\underline{\xi}] = \bigl[v_f,v_{\mu_\xi}\bigr] = v_{\{f,\mu_\xi\}}$ since according to our conventions the comoment map, the assignment of Hamiltonian vector fields, and the infinitesimal action are each Lie algebra morphism (see Remark \ref{rem:comomentAsLift}).
	From the expression of the Hamiltonian vector of the product of two functions, we have
	\[
		\sum_i g_i\underline{\xi}_i =
		\sum_i g_i v_{\mu_{\xi_i}} =
		v_{\sum_i g_i \underline{\xi}_i} - \sum_i \mu_i v_{g_i}  ,
	\]
	where the last term vanishes along $N=\mu^{-1}(0)$.
	According to Remark \ref{rem:comparingAlongN}, equality along $N$ implies
	\[
		v_{\{f,\mu_\xi\}} - v_{\sum_i g_i \underline{\xi}_i} \in I_{\mathfrak{X}}(N),
	\]
	from which $\{f,\mu_\xi\} - \sum_i g_i \underline{\xi}_i \in \mathcal{Q}$ and so
	\[
		\{f,\mu_\xi\} \subset I_\mu + \mathcal{Q} \qquad \forall \xi \in \g.
	\]
	Introducing $\mathcal{Q}'$ as in equation~\eqref{eq:calQprime}, the latter implies that
	\[
		L_\infty(M,\omega)_{[N]} = F_N \cap \mathcal{Q'}.
	\]
	Consider now $f \in I_{L_\infty}(N)$. The fourth condition in diagram \eqref{diagram:symplectic-ingredients} implies the existence of two families $g_i \in C^\infty(M)$ and $\xi_i \in \g$ such that, along $N$,
	\[
		v_f = \sum_i g_i  v_{\mu_i} =
		v_{\textstyle\sum_i g_i\mu_i} - \sum_i \cancel{\mu_i}v_{g_i} ,
	\]
	where $\mu_i = \mu_{\xi_i}$.
	By a similar argument as that of Remark \ref{rem:comparingAlongN}, we obtain $f\in I_\mu + \mathcal{Q}$.
	In other terms,
	\[
		I_{L_\infty}(N)=L_\infty(M,\omega)_{[N]} \cap (I_\mu+ I_N \cap\mathcal{Q}).
	\]
	The last claim follows by noting that both $I_\mu$ and $I_N\cap\mathcal{Q}$ lie inside $F_N\cap I_N \cap \mathcal{Q}'$ (see Lem\-ma~\ref{lem:Qprimeproperties}).
\end{proof}

\begin{Lemma}[properties of $\mathcal{Q}'$]\label{lem:Qprimeproperties}
	The subspace $\mathcal{Q}'$ enjoys the following properties:
	\begin{enumerate}\itemsep=0pt
		\item[$(i)$] $\mathcal{Q}'$ is an associative subalgebra of $C^\infty(M)$,
		\item[$(ii)$] $F_N\cap\mathcal{Q}'=\{f\in F_N \mid  \{f,I_\mu\}\subset I_\mu + \mathcal{Q}\cap I_N \}$ and it is a Poisson subalgebra of $C^\infty(M,\omega)$,
		\item[$(iii)$] $I_\mu\hookrightarrow \mathcal{N}\bigl(I_\mu\bigr)\hookrightarrow \mathcal{Q}'$ as associative algebras,
		\item[$(iv)$] $I_N\cap Q \hookrightarrow F_N\cap I_N \cap \mathcal{Q}' \hookrightarrow \mathcal{Q}'$ as associative algebras,
		\item[$(v)$] $I_\mu+I_N\cap Q$ is a Poisson ideal inside of $\mathcal{Q}'\cap F_N$.
	\end{enumerate}
\end{Lemma}
\begin{proof}
		(i)
		Consider $f$, $g$ in $\mathcal{Q}'$. By the distributive property of the Poisson bracket, we have for any $\xi \in g$ that
		\[
			\{f g, \mu_\xi\} \in \mathcal{Q}'\cdot\bigl\{\mathcal{Q}',\mu_\xi\bigr\}
			\subset
			\mathcal{Q}'\cdot I_\mu + (I_N\cap \mathcal{Q})\cdot \mathcal{Q}'
			\subset
			I_\mu+ I_N \cap \mathcal{Q},
		\]
		where in the last inclusion we employed Lemma \ref{lem:Qproperties}.

		(ii)
		The equality expresses the fact that the defining condition for $\mathcal{Q}'$ may be given in terms of $I_\mu$.
		This follows from the fact that $f \in F_N$ implies $\{f,\mu_\xi\}\subset I_N$, and from the fact that $Q\cap I_N\subset C^\infty(M)$ is an associative ideal by statement (iii) of Lemma \ref{lem:Qproperties}.
		As the intersection of the associative algebra $\mathcal{Q}'$ and the Poisson algebra $F_N$ it is, in particular, a Poisson algebra.

		(iii)
			The inclusion is immediate as $\mathcal{N}\bigl(I_\mu\bigr)$ is a Lie algebra normalizer.
		
        (iv)
			Clearly $I_N\cap Q$ includes in $I_N\cap F_N$.
			To prove the inclusion in $\mathcal{Q}'$ notice first that $\{(I_N\cap \mathcal{Q}),I_\mu\}\subset I_N$.
			The Jacobi identity implies that $\{(I_N\cap \mathcal{Q}),I_\mu\}$ lies in $\mathcal{Q}$ since
			\[
				 \{\{I_N\cap\mathcal{Q},I_\mu\},C^\infty(M) \}
				\subset
				 \{\{I_N\cap\mathcal{Q},C^\infty(M)\},I_\mu \}
				+
				 \{\{I_\mu,C^\infty(M)\},I_N\cap\mathcal{Q}  \}
				\subset
				I_N
			\]
	by the definition of $\mathcal{Q}$ and the inclusion $I_\mu\subset F_N$ (see \cite[Proposition~5.1]{ArmsGotayJennings90}).
\end{proof}

\begin{Remark}
	Observe that $F_N$ comprises precisely those functions with Hamiltonian vector field tangent to $N$.
	Thus we see from statement (ii) of Lemma \ref{lem:Qprimeproperties} that the reducible observables in the symplectic case consists of those functions in $F_N$ that satisfy a slightly stronger condition of being preserved along $N$. Specifically, for any $f \in\mathcal{Q}'\cap F_N$ and $\xi \in \g$, the Lie derivative $\L_\xi f \in I_N$ is a linear combination of an element $I_\mu$ and a constraint with Hamiltonian vector field vanishing along $N$.
\end{Remark}

\begin{Remark} The above formula \eqref{eq:1plred} for $L_\infty$-reduction can be applied to any triple $(P,I,J)$, where $P$ is a Poisson algebra, $I\subset J$ are associative ideals in $P$ and $I$ is stable under the Poisson bracket (i.e., also a Poisson subalgebra).  In our setting $P=C^{\infty}(M)$, $I=I_\mu$, $J=I_N$, however the proof of Lemma \ref{lem:Qprimeproperties} assuring that the quotient \eqref{eq:1plred} is a Poisson algebra carries over directly to general triples $(P,I,J)$. In particular, we get a reduction procedure for moment maps on Poisson manifolds and actions of connected Poisson--Lie groups on Poisson manifolds (cf., e.g.,~\cite{Esposito2014}).
\end{Remark}

We proceed now to compare the $L_\infty$-reduction procedure with the symplectic reduction schemes surveyed in Appendix \ref{Sec:SymSingRed}.
We denote by MW, SW, D, ACG, and AGJ the reduction procedures introduced by Marsden--Weinstein \cite{MarsdenWeinstein74},  \'Sniatycki--Weinstein \cite{SniatyckiWeinstein83}, Dirac \cite{Dirac64}, Arms--Cushman--Gotay \cite{ArmsCushmanGotay91}, and Arms--Gotay--Jennings \cite{ArmsGotayJennings90} respectively.
The hypotheses required by these constructions are summarized in Table \ref{tab:reduction_hypotheses}.

\begin{table}[h!]
	\centering\renewcommand{\arraystretch}{1.12}
    \caption{Reduction scheme hypotheses.}
	\label{tab:reduction_hypotheses}\vspace{1mm}

	\begin{tabular}{|r|ccc|}
		\hline
		reduction		&
		\multicolumn{3}{c|}{hypothesis}
		\\
		scheme & action $G\action M$ &constraint set $N$ &$\omega\in \Omega^2(M)$
		\\[2pt] \hline
		AGJ
		& none
		&
		``well behaved''\footnotemark
		& symplectic
		\\ \hline
		D
		& none
		& first class
		& symplectic
		\\ \hline
		$L_\infty$
		& any
		& $\g$-stable
		& \begin{tabular}[c]{@{}l@{}}
			presymplectic \\
			reducible
		 \end{tabular}
		\\ \hline
		\'SW
		& Hamiltonian
		& $N=\mu^{-1}(0)$
		& symplectic
		\\ \hline
		ACG
		& Hamiltonian
		& $N=\mu^{-1}(0)$
		& symplectic
		\\ \hline
		MW
		& \begin{tabular}[c]{@{}l@{}}
			Hamiltonian\\
			free and proper on $\mu^{-1}(0)$
		  \end{tabular}
		& \begin{tabular}[c]{@{}l@{}}
			$N=\mu^{-1}(0)$\\
			$0$ regular value for $\mu$
		  \end{tabular}
		& symplectic
		\\ \hline
	\end{tabular}
\end{table}
\footnotetext{For example, strongly coisotropic and locally conical, see \cite[Proposition~3.3]{ArmsGotayJennings90}.}

Let us focus on the \emph{symmetry-based} observable reductions, that is, we assume  $N=\mu^{-1}(0)$ for a moment map $\mu$. In this case, both the [\'SW], the [ACG], and the $L_\infty$-reduction schemes apply. In particular, we have the following morphisms of Poisson algebra.

\begin{Proposition}[{[$L_\infty$] vs.\ [\'SW]}]\label{prop:LinftyVsSW}
	Fix a symplectic Hamiltonian action $G\action (M,\omega)$ with moment map $\mu\colon M\to\g^*$, and put $N=\mu^{-1}(0)$.
	When $\mathcal{N}\bigl(I_\mu\bigr)\hookrightarrow F_N$ we have a canonical Poisson morphism from the {\rm [\'SW]} reduction to the $L_\infty$-reduction.
	When $I_N\cap\mathcal{Q}\hookrightarrow I_\mu$ the $L_\infty$-reduction embeds into the {\rm [\'SW]} reduction.
	When $I_\mu=I_N$, the {\rm [\'SW]} and {\rm [$L_\infty$]} reduced spaces are isomorphic.
\end{Proposition}
\begin{proof}
	According to Section \ref{Ssec:LinftyReductionMomentMap}, $N=\mu^{-1}(0)$ implies that $\omega$ is reducible.
	First observe that under the above hypothesis the [\'SW] and [$L_\infty$] reductions are well-defined.
	We exhibit all pertinent relations between the considered spaces in the following diagram commutative diagram in the category of Poisson algebra:
	\[
		\begin{tikzcd}[column sep = large]
			I_{\mu} \ar[r,hook,"\text{\tiny (Lie.) ideal}"'] \ar[d,equal] &[1em]
			\mathcal{N}\bigl(I_\mu\bigr)\ar[r,two heads] &[2em]
			\dfrac{\mathcal{N}\bigl(I_\mu\bigr)}{I_\mu}
			\\
			I_{\mu} \ar[d,hook]\ar[r,hook,"\text{\tiny (Lie.) ideal}"'] &
			\mathcal{N}\bigl(I_\mu\bigr)\cap F_N \ar[u,hook]\ar[d,equal]\ar[r,two heads]&
			\dfrac{\mathcal{N}\bigl(I_\mu\bigr)\cap F_N}{I_{\mu}}		\ar[u,hook]\ar[d,two heads]
			\\
			I_{\mu} + I_N\cap \mathcal{Q}\cap \mathcal{N}\bigl(I_\mu\bigr) \ar[d,hook]\ar[r,hook,"\text{\tiny (Lie.) ideal}"'] &
			\mathcal{N}\bigl(I_\mu\bigr)\cap F_N \ar[d,hook]\ar[r,two heads]&
			\dfrac{\mathcal{N}\bigl(I_\mu\bigr)\cap F_N}{I_{\mu} + I_N\cap \mathcal{Q}\cap \mathcal{N}\bigl(I_\mu\bigr)}		\ar[d,hook]
			\\
			I_\mu+ I_N\cap\mathcal{Q} \ar[r,hook,"\text{\tiny (Lie.) ideal}"'] &
			\mathcal{Q}'\cap F_N \ar[r,two heads]&
			\dfrac{\mathcal{Q}'\cap F_N}{I_\mu+ I_N\cap\mathcal{Q}},
		\end{tikzcd}
	\]
	where the rightmost bottom arrow is given by the second isomorphism theorem of Lie algebras.
	This induces the following diagram at the level of quotient spaces, where dashed arrow denotes morphisms of associative algebras~--~not Poisson~-- (cf.\ Lemma~\ref{lem:ImuVSIN}),
	\[
		\begin{tikzcd}[column sep = large,
			cells={font=\everymath\expandafter{\the\everymath\displaystyle}},]
			\frac{\mathcal{N}\bigl(I_\mu\bigr)\cap F_N}{I_{\mu}}\ar[r,hook,sloped]
			\ar[dr,two heads,sloped]
			&
			\frac{\mathcal{N}\bigl(I_\mu\bigr)}{I_{\mu}} \ar[r,hook,dashed,sloped]
			&
			\frac{\mathcal{Q}'}{I_{\mu}}
			\ar[dr,two heads, dashed,sloped]
			\\
			&
			\frac{\mathcal{N}\bigl(I_\mu\bigr)\cap F_N}{I_{\mu} + I_N\cap \mathcal{Q}\cap \mathcal{N}\bigl(I_\mu\bigr)}	
			\ar[r,hook,,sloped]
			&
			\frac{\mathcal{Q}'\cap F_N}{I_{\mu}+I_N\cap\mathcal{Q}}
			\ar[r,hook,dashed,sloped]
			&
			\frac{\mathcal{Q}'}{I_{\mu}+I_N\cap\mathcal{Q}},
		\end{tikzcd}	
	\]
	where we identify the [\'SW] reduction in the middle of the first row (see Theorem \ref{thm:SW-red}) and the~[$L_\infty$] reduction in the middle of the second row (see Theorem \ref{thm:LinftyRedSymp}).
	The latter diagram means that the two considered reduction schemes yield different Poisson algebras in general and they posses a nontrivial intersection inside of the associative algebra $\mathcal{Q}'/(I_{\mu}+I_N\cap\mathcal{Q})$.

	When $F_N \hookrightarrow \mathcal{N}\bigl(I_\mu\bigr)$, the top left horizontal arrow in the previous diagram is an isomorphism and we obtain the sought map.
	Similarly, when $I_N\cap\mathcal{Q}\hookrightarrow I_\mu$ we have that $F_N\cap\mathcal{Q}'\cong \mathcal{N}\bigl(I_\mu\bigr)\cap F_N$, hence the two bottom Poisson map are indeed identifications.
	Finally, the condition $I_\mu=I_N$ implies the previous two; hence, the two reduced Poisson algebras coincide.
\end{proof}

\begin{Remark} The BFV-BRST approach to symplectic reduction leads to a differential graded algebra, whose zero-th homology coincides with the [\'SW]-reduction, cf., e.g., \cite{Stasheff1992}. Hence, in the setting of Proposition \ref{prop:LinftyVsSW}, the [$L_\infty$]-reduction also coincides with this homology group.
\end{Remark}

\begin{Proposition}[{[$L_\infty$] vs.\ [ACG]}]\label{prop:LinftyVsACG}
	Consider a symplectic Hamiltonian action $G\action (M,\omega)$ with moment map $\mu\colon M\to\g^*$.
	Assume that $N=\mu^{-1}(0)$ where $\mu$ is a moment map associated to the action.
	If $I_\mu= I_N$, then the {\rm [ACG]} reduction naturally embeds in the {\rm [$L_\infty$]} reduced algebra.
	If moreover $F_N\cong C^\infty(M)^G +I_N$, then the {\rm [ACG]} and {\rm [$L_\infty$]} reductions are isomorphic.
\end{Proposition}
\begin{proof}
	We indicate the relevant relations between all spaces in the following commutative diagram in the category of Poisson algebras:
	\[
		\begin{tikzcd}[
			column sep = large,
		cells={font=\everymath\expandafter{\the\everymath\displaystyle}},]
			I_{N}^G
			\ar[r,hook,"\text{\tiny (Lie.) ideal}"']
			&[1em] C^\infty(M)^G  \ar[r,two heads]
			&[3em] \frac{C^\infty(M)^G}{I_N^G}
			\\
			I_N^G\cap F_N
			\ar[r,hook,"\text{\tiny (Lie.) ideal}"']
			\ar[ur,phantom,very near start, "\urcorner"]
			\ar[u,hook]
			& C^\infty(M)^G\cap F_N  \ar[r,two heads]\ar[u,hook]
			&[3em] \frac{C^\infty(M)^G\cap F_N}{I_{N}^G\cap F_N}\ar[u,hook]
			\\
			I_\mu^G + I_N^G\cap \mathcal{Q}
			\ar[u,hook]
			\ar[r,hook,"\text{\tiny (Lie.) ideal}"']
			\ar[dr,phantom,very near start, "\lrcorner"] \ar[d,hook]
			& C^\infty(M)^G\cap F_N  \ar[r,two heads]\ar[u,equal]\ar[d,hook]
			&[3em] \frac{C^\infty(M)^G\cap F_N}{I_\mu^G + I_N^G\cap \mathcal{Q}}\ar[d,hook]\ar[u,two heads]
			\\
			I_\mu + I_N\cap \mathcal{Q}
			\ar[r,hook,"\text{\tiny (Lie.) ideal}"']
			& \mathcal{Q}'\cap F_N  \ar[r,two heads]
			&[3em] \frac{\mathcal{Q}'\cap F_N}{I_\mu + I_N \cap \mathcal{Q}},
		\end{tikzcd}
	\]
	where the top and bottom squares encode the second Lie algebra isomorphism theorem.
	The inclusion of $C^\infty(M)^G$ in $\mathcal{Q}'$ follows from $C^\infty(M)^G\subset \mathcal{N}\bigl(I_\mu\bigr)$.
	The inclusion of ${I_{\mu}}^G + {I_{N}}^G\cap \mathcal{Q} \subset {I_{N}}^G\cap F_N$ follows from the inclusion of $I_\mu$ and $\mathcal{Q}$ in $F_N$ (\cite[Proposition~5.1]{ArmsGotayJennings90} and Lemma \ref{lem:Qproperties}).
	Note that the bottom right object coincides with the [$L_\infty$] reduction (see Theorem \ref{thm:LinftyRedSymp}), and the top right object corresponds to the [ACG] reduction as recalled in Definition \ref{def:ACGreduction}.
	The first claim follows by noticing that the rightmost two top vertical arrows are indeed identifications whenever $I_\mu= I_N$.
	If, moreover,  $F_N\cong C^\infty(M)^G +I_N$, the bottom right vertical arrow yields an isomorphism.
\end{proof}

\begin{Remark}
	Note that the inclusion from Proposition \ref{prop:LinftyVsACG} follows also when $C^\infty(M)^G\hookrightarrow \mathcal{Q}$, i.e., if the Hamiltonian vector field associated to any $G$-invariant smooth function vanishes along~$N$.
\end{Remark}

When the constraint set $N=\mu^{-1}(0)$ is also first class, it is possible to compare the~[$L_\infty$] and~[D] reduction schemes.
\begin{Proposition}[{[$L_\infty$] vs.\ [D]}]\label{prop:LinftyVsD}
	Let $G\curvearrowright(M,\omega)$ be a symplectic Hamiltonian action with moment map $\mu\colon M\to\g^*$.
	Further assume that $N=\mu^{-1}(0)$ is first class.
	There is a natural mapping from the $L_\infty$-reduced Poisson algebra to the {\rm [D]} reduction.
	If $I_N=I_\mu+I_N\cap\mathcal{Q}$, then the {\rm [D]} and {\rm [$L_\infty$]} reduced spaces coincide.
\end{Proposition}

\begin{proof}
	The relevant mappings are indicated in the following commutative diagram in the category of Poisson algebras:
	\[
		\begin{tikzcd}[
			column sep = large,
		cells={font=\everymath\expandafter{\the\everymath\displaystyle}},]
			I_{\mu} + I_N\cap \mathcal{Q} \ar[r,hook,"\text{\tiny (Lie.) ideal}"']
			& \mathcal{Q}'\cap F_N \ar[r,two heads,] \ar[d,equal]
			&[2em] \frac{\mathcal{Q}'\cap F_N}{I_{\mu} + I_N\cap \mathcal{Q}}
			\ar[d,two heads]
			\\
			I_N\cap \mathcal{Q} \ar[r,hook,"\text{\tiny (Lie.) ideal}"']\ar[d,hook]
			\ar[dr,phantom,very near start, "\lrcorner"]
			\ar[u,hook]
			& \mathcal{Q}'\cap F_N \ar[r,two heads,] \ar[d,hook]
			&[2em] \frac{\mathcal{Q}'\cap F_N}{I_N\cap \mathcal{Q}} \ar[d,hook]
			\\
			I_N\ar[r,hook,"\text{\tiny (Lie.) ideal}"']
			& F_N \ar[r,two heads,]
			&	
			\frac{F_N}{I_{N}}.
		\end{tikzcd}
	\]
	The top and the bottom embody the [$L_\infty$] and [D] reductions.
	The two bottom squares encode the second isomorphism theorem of Lie algebras.
	In the case that $I_N=I_\mu+I_N\cap\mathcal{Q}$, we have in particular that $\mathcal{Q}'\cap F_N = F_N$, hence the three Poisson algebras in the rightmost column are identified.
\end{proof}

\begin{Remark}
Albeit the $L_\infty$-reduction procedure draws its original inspiration from the~[SW] scheme, from the proof of Proposition \ref{prop:LinftyVsD} and the simple observation that $\mathcal{O}(N)\supset F_N\cap\mathcal{Q}'$ transpires how the $L_\infty$-reduction is rather a refinement of the [D] reduction.
More specifically, such a refinement is obtained by considering the subset of Dirac's observables consisting of suitably locally preserved first-class functions and dividing out by a subset of functions vanishing on~$N$. These choices, in particular, imply that this construction yields a Poisson algebra regardless of the condition of being~$N$ a first-class constraint.

The laxer relationship with the [SW] reduction is unsurprising as that scheme relies fundamentally on the associative structure of $C^\infty(M)$ which, in principle, does not carry to $\Ham_\infty(M,\omega)$.
\end{Remark}

An observable reduction scheme cannot be considered entirely satisfactory if it were not to agree with the [MW] reduction in the presence of a regular constraint set.
Crucially, all the reduction schemes mentioned in Appendix \ref{Sec:SymSingRed} satisfy this compatibility property (see \cite[Theorem~1]{SniatyckiWeinstein83} and \cite[Proposition~3.6]{ArmsGotayJennings90} for further details).
More formally, under the hypotheses of the Marsden--Weinstein reduction theorem, the reduced Poisson algebra $C^\infty(M,\omega)_0$ is isomorphic to the Poisson algebra $C^\infty(M_0,\omega_0)$ on the reduced symplectic manifold.

\begin{Theorem}\label{thm:AllReductionsComparison-sympleticcase}
	Let $G\action (M,\omega)$ be a symplectic Hamiltonian action and suppose that $0\in\g^*$ is a regular value of the moment map $\mu\colon M\to\g^*$.
	If $G\curvearrowright M$ is free and proper, then the {\rm [$L_\infty$]}, {\rm [\'SW]}, {\rm [D]}, and {\rm [ACG]} reductions are equal.
	In particular, each is isomorphic to the Poisson algebra $C^\infty(M_0,\omega_0)$ of smooth functions on the Marsden--Weinstein reduced space.
\end{Theorem}

\begin{proof}
	First observe that the above hypothesis coincide with those of Theorem~\ref{thm:MW-red}, hence the [MW] geometric reduction is well defined.
	Let us denote by $(M_0,\omega_0)$ the symplectic manifold obtained by the Marsden--Weinstein reduction procedure.
	The regular value condition guarantees that~$\mu^{-1}(0)$ is an embedded submanifold of $M$. In particular, we have $I_\mu= I_N$ (see \cite[Theorem~1]{SniatyckiWeinstein83} and \cite[Proposition~5.12]{ArmsGotayJennings90} for details).
	According to Propositions \ref{prop:LinftyVsSW} and \ref{prop:LinftyVsD} the latter condition implies that the [$L_\infty$], [\'SW] and [D] reduction schemes yield the same Poisson algebra.
	In particular, all of them coincide with $C^\infty(M_0,\omega_0)$ in virtue of \cite[Theorem~1]{SniatyckiWeinstein83}.
	The isomorphism between the latter and the [ACG] reduction is given in \cite[Theorem~1]{Sniatycki05}.
\end{proof}

The following diagram subsumes the relationship between the considered reduction schemes:
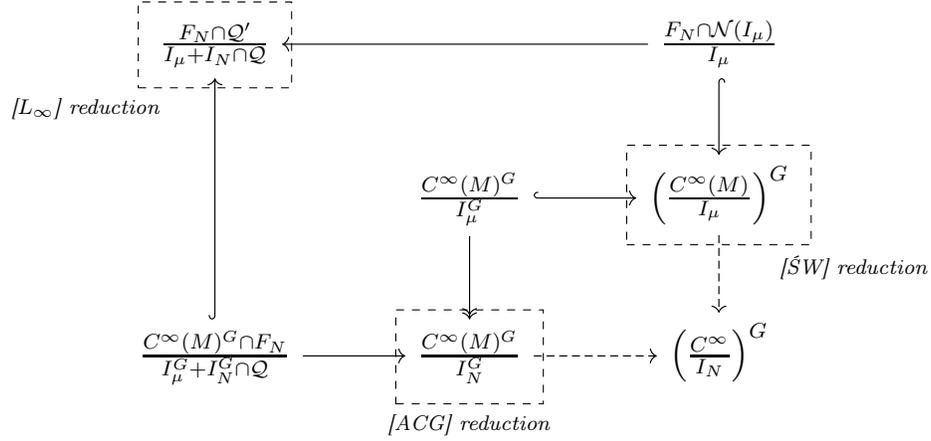
\begin{displaymath}
	\begin{tikzcd}[
			row sep = large,
			/tikz/execute at end picture={
	    			\node (large) [rectangle,dashed, draw, fit=(SW),label=below right:{\scriptsize{\emph{{\rm [\'SW]} reduction}}}] {};
	    			\node (large) [rectangle,dashed, draw, fit=(Linf),label=below left:{\scriptsize{\emph{$[L_\infty]$ reduction}}}] {};
	    			\node (large) [rectangle,dashed, draw, fit=(ACG),label=below:{\scriptsize{\emph{{\rm [ACG]} reduction}}}] {};
	    		}
	    	]	
		|[alias=Linf]|\dfrac{F_N\cap \mathcal{Q}'}{I_\mu + I_N\cap\mathcal{Q}} & &
		\dfrac{F_N\cap\mathcal{N}\bigl(I_\mu\bigr)}{I_\mu}\ar[ll]\ar[d,hook]
		\\
		& \dfrac{C^\infty(M)^G}{I_\mu^G} \ar[r,hook]\ar[d,two heads]&
		|[alias=SW]|\left(\dfrac{C^\infty(M)}{I_\mu}\right)^G \ar[d,dashed]
		\\
		\dfrac{C^\infty(M)^G\cap F_N}{I_\mu^G + I_N^G\cap\mathcal{Q}} \ar[uu,hook]\ar[r]&
		|[alias=ACG]|\dfrac{C^\infty(M)^G}{I_N^G} \ar[r,dashed] &
		\left(\dfrac{C^\infty}{I_N}\right)^G.
	\end{tikzcd}
\end{displaymath}
As before, the dashed (resp.\ solid) arrows indicate morphisms of associative (resp.\ Poisson) algebras.

\subsection[The coordinate cross in R\^{}2]{The coordinate cross in $\boldsymbol{\R^2}$}\label{Ssec:coordinate_cross}

Let us consider an elementary example in which the $L_\infty$-reduction of a symplectic manifold is possible even though the action fails to preserve the symplectic form.

Consider the plane $M=\mathbb R^2$ and the union of the coordinate axes $N=\{ (x,y) \mid  xy=0\}$. The vanishing ideal is $I_N=xy\cdot C^\infty(M)$. We consider the vector field $\underline\xi =x\partial_x+y\partial_y$ as the action of a one-dimensional abelian Lie algebra on $M$. The vector fields tangent to $N$ are generated by $x\partial_x$ and $y\partial_y$. In particular, $\omega=\d x\wedge \d y$ is reducible, as it is closed and as $\iota_\xi\omega= x \d y-y\d x$ vanishes on $N$ when contracted with $x\partial_x$ or $y\partial_y$.

As described in the preceding section, the reduction is the quotient
\[
	\frac{ \bigl\{ f\mid  v_f\in \tgtvf,\, \L_\xi f\in I_N,\, [\underline{\xi},v_f]\in \fgmodule(M) + \vanvf \bigr\}
	}{
	\{
	f\mid  f\in I_N,\, v_f\in \fgmodule(M) + \vanvf
	\}
	}.
\]

We first consider the denominator. Since $f\in I_N$, we have $f=xyh$ for some smooth function~$h$. Thus $\d f=y\bigl(h+x\frac{\partial h}{\partial x}\bigr)\d x +x\bigl(h+y\frac{\partial h}{\partial y}\bigr)\d y$, from which $v_f=y\bigl(h+x\frac{\partial h}{\partial x}\bigr)\partial_y - x{\bigl(h+y\frac{\partial h}{\partial y}\bigr)}\partial_x$. Such expressions are required to lie in $\fgmodule(M)+\vanvf$. Now, the equality of ideals $\vanvf=I_N \mathfrak X(M)$ in this example yields $v_f=h(y\partial_y-x\partial_x)$ up to an element of $\vanvf$. As multiples of $xy$ in $h$ are redundant, we arrive at a decomposition $h(x,y)=h^0 + xh^x(x)+yh^y(y)+xyh^{\rm res}$. Then $v_f= \bigl(h^0+xh^x(x)+yh^y(y)\bigr)(y\partial_y-x\partial_x)$. Multiplying and removing all multiples of $xy$, we obtain
\[
v_f=h^0(y\partial_y-x\partial_x) + y^2h^y(y)\partial_y - x^2h^x(x)\partial_x .
\]
Since $x$ and $y$ are fully decoupled, there are no terms in $\vanvf$. It remains to check when this lies in $\fgmodule(M)$. This requires that $h^0$, $h^x$ and $h^y$ be identically zero, so that $h(x,y)=xyh^{\rm res}(x,y)$ and thus $f=(xy)^2h^{\rm res}(x,y)$ for some smooth function $h^{\rm res}(x,y)$.

Now for the numerator. Since multiples of $(xy)^2$ vanish in the quotient, we introduce the ansatz
\[
	f(x,y)= f^0 + xf^x(x) +yf^y(y) + xyf^{xy|x}(x) + xyf^{xy|y}(y).
\]
Taking exterior derivatives gives
\begin{gather*}
	\d f= \bigl(f^x(x)+x(f^x)'(x)+ yf^{xy|x}(x) +xy\bigl(f^{xy|x}\bigr)'(x)   +yf^{xy|y}(y)\bigr)\d x		\\ \hphantom{\d f=}{}
		+ \bigl(f^y(y)+y(f^y)'(y)+ xf^{xy|y}(y) +xy\bigl(f^{xy|y}\bigr)'(y)   +xf^{xy|x}(x)\bigr)\d y.
\end{gather*}
Before calculating $v_f$, let us look at $\L_\xi f=\d f(\underline{\xi})$,
\begin{gather*}
	\d f(\underline{\xi})
        =	xf^x(x)+x^2(f^x)'(x)+ xyf^{xy|x}(x) +x^2y\bigl(f^{xy|x}\bigr)'(x)   +xyf^{xy|y}(y)	\\ \hphantom{\d f(\underline{\xi})=}{}
		+	yf^y(y)+y^2(f^y)'(y)+ xyf^{xy|y}(y) +xy^2\bigl(f^{xy|y}\bigr)'(y)   +xyf^{xy|x}(x).
\end{gather*}
As this function lies in $I_N$, it follows that $f^x$ and $f^y$ must vanish. In terms of our expression for~$f$,
\[
	v_f=    \bigl(yf^{xy|x}(x) +xy\bigl(f^{xy|x}\bigr)'(x)   +yf^{xy|y}(y)\bigr)\partial_y - \bigl(xf^{xy|y}(y) +xy\bigl(f^{xy|y}\bigr)'(y)   +xf^{xy|x}(x)\bigr)\partial_x.
\]
Since $\tgtvf$ is generated by $x\partial_x$ and $y\partial_y$, it follows that $v_f\in\tgtvf$. We now turn to the commutator $[\underline\xi,v_f]$. Since $\underline\xi$ preserves $I_N$, it also preserves $\vanvf=I_N\mathfrak X(M)$. We may thus disregard all terms which are multiples of $xy$ and write
\begin{align*}
	[\underline\xi,v_f]\equiv y^2\bigl(f^{xy|y}\bigr)'(y)\partial_y - x^2\bigl(f^{xy|x}(x)\bigr)'\partial_x.
\end{align*}
This is a multiple of $\underline\xi$ precisely when \smash{$\bigl(f^{xy|y}\bigr)'=\bigl(f^{xy|x}\bigr)'=0$}, that is, precisely when~$f^{xy|y}$ and~$f^{xy|x}$ are constant. Thus, up to elements in the denominator, $f$ takes the form
\[
f=f^0 + xyf^{xy}
\]
for arbitrary constants $f^0, f^{xy}\in\R$. In particular, the reduced $L_\infty$-algebra is isomorphic to the abelian Lie algebra $\mathbb R^2$.

\subsection{Multicotangent bundles}\label{Ssec:multicotangent}

As a precursor to the primary example of multiphase spaces, let us examine the intermediate construction of multicotangent bundles.

Fix a manifold $E$. The $n$-th \emph{multicotangent bundle} of $E$ is the manifold $\Lambda^nT^*E$. As $T^*E$ carries a canonical $1$-form, so $\Lambda^nT^*E$ carries a \emph{canonical $n$-form} $\theta\in \Omega^n(\Lambda^nT^*E)$. Explicitly,
\[
	\theta_\eta(v_1,\dots,v_n)=\eta(\pi_*v_1,\dots,\pi_*v_n),
\]
where $\eta\in\Lambda^nT^*E$, $v_1,\dots,v_n\in T_\eta\Lambda^nT^*E$, and $\pi\colon \Lambda^nT^*E\to E$ is the projection. Local coordinates $(e_i)_{i\leq n}$ on $U\subset E$, induce local coordinates $\bigl(e_i, p^I\bigr)$ on $\pi^*U\subset \Lambda^nT^*E$, where $I=1\leq i_1<\dots<i_n\leq n$ is a multiindex of length $n$, and where $p^I$ represents the coefficient of $\d e_I=\d e_{i_1}\wedge \cdots\wedge \d e_{i_n}$. In these terms,
\[
	\theta= {\textstyle\sum} p^I\,\d e_I = {\textstyle\sum} p^I\,\d e_{i_1}\wedge \cdots\wedge \d e_{i_n}.
\]
The \emph{canonical multisymplectic structure} on $\Lambda^nT^*E$ is $\omega=\d\theta$. Note that any diffeomorphism $\psi$ of $E$ naturally extends to a diffeomorphism $\Psi=\bigl(\psi^{-1}\bigr)^*$ of $\Lambda^nT^*E$. On the infinitesimal level, this induces an inclusion of Lie algebras
\begin{align*}
	\mathfrak X(E)	&\hookrightarrow	\mathfrak X(\Lambda^nT^*E),	\\
	v		&\mapsto		\tilde v
\end{align*}
called the \emph{prolongation map} along $\pi\colon \Lambda^nT^*E\to E$. From $\Psi^*\theta=\theta$ for all diffeomorphisms $\psi\colon E\to E$, we easily derive $\L_{\tilde v}\theta=0$ for all vector fields $v\in\X(E)$, whence $\iota_{\tilde v}\omega + \d\iota_{\tilde v}\theta = 0$ so that $\tilde{v}$ is a Hamiltonian vector field for $\iota_{\tilde v}\theta\in\Omega^{n-1}(\Lambda^nT^*E)$. In particular, any smooth action $\g\curvearrowright E$ prolongs to a $\theta$-preserving action $\g\curvearrowright \Lambda^nT^*E$ and thus Lemma \ref{lem:potinducesmom} provides a \emph{canonical moment map} $\mu\in\Omega^{n-1}(\Lambda^nT^*E,\g^*)$ according to $\mu_\xi = \iota_\xi\theta$.

From
\[
	\alpha\in\mu^{-1}(0) \iff \iota_\g\theta_\alpha = 0 \iff \iota_{\underline\g_E}\alpha = 0,
\]
we obtain the $0$-level set $\mu^{-1}(0)\subset\Lambda^nT^*E$ as the annihilator
\[
	N=\mu^{-1}(0)
	= \bigl\{\alpha\in\Lambda^nT^*E \mid  \iota_{\underline\g_E}\alpha = 0 \bigr\}.
\]

Notwithstanding the apparent difficulty of computing the full reduced space of observables, certain elements are always present.

\begin{Proposition}\label{prop:reducible_observables_on_multicotangent_bundles}
	If $v\in\X(E)$ preserves $\X_\g(E)$, then $(\tilde v,\iota_{\tilde v}\theta)\in\Ham_\infty^0(\Lambda^nT^*E,\omega)$ is reducible.
\end{Proposition}

\begin{proof}
	More formally, our condition on $v$ is that $[v, \X_\g(E)]\subset \X_\g(E)$. For the reducibility of $\tilde v\in \X(\Lambda^nT^*E)$, we must check two conditions.
	
	First, we show that $[\tilde v, \underline \xi]\in \X_\g(\Lambda^nT^*E)+ I_\X(N)$ for all $\xi \in \g$. We will use the fact that $[\tilde v, \underline \xi]$ is the Hamiltonian vector field of $\iota_{[\tilde v, \underline \xi]}\theta$, and that $[\tilde v, \underline\xi]$ is the prolongation of the vector field  $[v, \underline\xi_E]\in \X(E)$. By hypothesis, $[v,\underline\xi_E] = \sum_i f_i\underline\zeta_{Ei}$ for some $f_i\in C^\infty(E)$ and $\zeta_i\in\g$, from which $\iota_{[\tilde v, \underline \xi]}\theta=\sum f_i\iota_{\underline \zeta_i}\theta$. We have
	\[
		\d\Bigl(\sum f_i\iota_{\underline\zeta_i}\theta\Bigr)=\sum f_i \d\iota_{\zeta_i}\theta_i +\d f_i\wedge \iota_{\zeta_i}\theta.
	\]
	At every point $x\in E$, each term $f_i\d\iota_{\zeta_i}\theta_i$ is related by the injective map $v\mapsto\iota_v\omega$ to $f_i \zeta_i$; while, at every point $x\in N=\{x\mid \iota_{\underline{\zeta}}\theta=0~\forall\zeta\}$, each term $\d f_i\wedge\iota_{\zeta_i}\theta$ vanishes. Therefore, the injectivity of $v\mapsto\iota_v\omega$ yields $[\tilde v, \underline \xi]\in \X_\g(\Lambda^nT^*E)+I_\X(N)$.
	
	Second, we establish that $\tilde v$ is tangent to $N$. In fact, this is true even in the stronger sense that the flow $\Psi_t$ of $\tilde v$ preserves $N$, that is, that $\eta\in N$ implies $\Psi_t(\eta)\in N$. We will show that $(\iota_\zeta\theta)_{\Psi_t\eta}=0$ if and only if $\Psi_t^* (\iota_\zeta\theta)_\eta=0$. Since $\Psi_t$ preserves $\theta$ this is equivalent to $ \bigl(\iota_{(\Psi_t)_*\underline{\zeta}}\theta\bigr)_{\eta}=0$. Denote by $\psi$ the flow of $v$ on $E$ and observe that $\theta$ depends only on $(\psi_{t})_*\underline{\zeta}=\pi_*(\Psi_{t})_*\underline{\zeta}$. Invoking our hypothesis that $[v, \X_\g(E)]\subset \X_\g(E)$, we note that $\psi_t$ preserves $\X_\g(E)$, and hence that ${(\psi_{t})_*\underline{\zeta}}_{\pi(\eta)}$ is a linear combination of fundamental vectors of $\g$. This provides $\bigl(\iota_{(\Psi_t)_*\underline{\zeta}}\theta\bigr)_{\eta}$, from which we obtain $(\iota_\zeta\theta)_{\Psi_t\eta}$.
	
	For the reducibility of $\iota_{\tilde{v}}\theta$, we require that $\iota_\xi\iota_{\tilde{v}}\theta \in I_\Omega(N)$ and that $\L_\xi\iota_{\tilde{v}}\theta = [\L_\xi,\iota_{\tilde{v}}]\,\theta = \iota_{[{\underline{\xi}},\tilde{v}]}\theta\in I_\Omega(N)$ for all $\xi\in\g$.  The first containment is an immediate consequence of the equality $N= \bigl\{\alpha\in\Lambda^nT^*E \mid  \iota_{\underline\g_E}\alpha = 0 \bigr\}$. The second follows from this characterization of $N$ and the from the property that $[{\underline{\xi}},\tilde{v}]\in\X_\g(\Lambda^nT^*E)+I_\X(N)$.
\end{proof}

\subsection{Multiphase spaces}\label{Ssec:Multiphase}

We now consider multiphase spaces, arguably the foremost class of examples of multisymplectic manifolds, playing an essential role in the covariant Hamiltonian description of first-order classical field theories. See \cite{Helein12,Kijowski73} for a shorter exposition and \cite{GotayIsenbergMarsden04, GotayIsenbergMarsdenMontgomery98} for a comprehensive treatment.

The \emph{multiphase space} associated to a surjective submersion $\pi\colon E\to\Sigma$ is
\[
	\Lambda_1^nT^*E = \{\eta \in \Lambda^nT^*E\mid \iota_u\iota_v\eta=0\ \forall u,v\in T_{\rm vert}E\}.
\]
Traditionally, we consider $\Sigma$ to embody the independent variables of a classical field theory, and the fiber directions of $E$ the dependent variables.

When the dimension of $\Sigma$ is sufficiently great, $\Lambda^n_1T^*E$ is a multisymplectic submanifold of $\Lambda^nT^*E$ and inherits the constructions of Section~\ref{Ssec:multicotangent}.

We begin with the extension of vector fields. As it is not generally true that the prolongation of $v\in\X(E)$ to $\Lambda^nT^*E$ is tangent to $\Lambda_1^nT^*E$, our first task is to identify those $v\in\X(E)$ that do prolong to $\Lambda_1^nT^*E$.

\begin{Lemma}\label{lem:projectable_vfs_prolong_to_multiphase_space}
	If $v\in\X(E)$ is $(E\to \Sigma)$-projectable, then the prolongation $\tilde{v}\in\X(\Lambda^n T^*E)$ is tangent to $\Lambda_1^n T^*E$.
\end{Lemma}

\begin{proof}
	Suppose $v\in\X(E)$ projects along $\pi\colon E\to\Sigma$ to some $w\in\X(\Sigma)$. Then there are smooth families of diffeomorphisms $(\Psi_t)_t$ and $(\psi_t)_t$ such that $v=\frac{\d}{\d t}\Psi_t|_{t=0}$, $w=\frac{\d}{\d t}\psi_t|_{t=0}$, and
	\[
		\pi\circ\Psi_t = \psi_t\circ\pi.
	\]
	Since $\tilde{v} = \frac{\d}{\d t}\Psi_{-t}^*|_{t=0}$ it suffices to show that each $\Psi_t^*$ preserves $\Lambda_1^nT^*E$. To see this is so, fix $t\in\R$, put $\Psi=\Psi_t$, let $\alpha\in\Lambda_1^nT^*E$ and $u,u'\in T_{\rm vert}E$, and observe that
	\[
		\iota_u\iota_{u'} (\Psi^*\alpha) = \Psi^*(\iota_{\Psi_*u}\iota_{\Psi_*u'}\alpha) = 0.
	\]
	Note that the last equality follows as $\alpha\in\Lambda_1^nT^*E$, and as $\pi_*\Psi_*u=\psi_*(\pi_*u)=0$ provides $\Psi_*u,\Psi_*u'\in T_{\rm vert}E$.
\end{proof}

\begin{Lemma}\label{lem:projectable_vfs_are_Hamiltonian_in_multiphase_space}
	If $v\in\X(E)$ is $(E\to\Sigma)$-projectable, then $\tilde{v}\in\X\bigl(\Lambda_1^nT^*E\bigr)$ is a Hamiltonian vector field for $\iota_{\tilde{v}}\theta\in\Omega^{n-1}\bigl(\Lambda_1^nT^*E\bigr)$.
\end{Lemma}

\begin{proof}
	This is a straightforward consequence of Lemma \ref{lem:projectable_vfs_prolong_to_multiphase_space} and the expository review of Section~\ref{Ssec:multicotangent}.
\end{proof}

We now again turn our attention to the canonical moment map induced by $\theta=\theta|_{\Lambda_1^nT^*E}\in \Omega^n(\Lambda_1^nT^*E)$ for the action $G\curvearrowright\Lambda_1^nT^*E$ lifted from $G\curvearrowright (E\to \Sigma)$. If $G\curvearrowright\Lambda_1^nT^*E$ lifts the action $G\curvearrowright (E\to \Sigma)$, then it is easy to show that $G$ preserves $\theta$.

For the zero level set of the moment map, the same argument as with the multicotangent bundles applies and we obtain
\[
	N=\mu^{-1}(0) = \bigl\{\alpha\in\Lambda_1^nT^*E \mid  \iota_{\underline\g_E}\alpha = 0\bigr\}.
\]

We thus arrive at a natural class of reducible observables on $\Lambda_1^nT^*E$.

\begin{Theorem}
	If $v\in\X(E)$ preserves $\X_\g(E)$, the $C^\infty(E)$-module generated by the fundamental distribution $\underline\g$, and if $v\in\X(E)$ is $(E\to\Sigma)$-projectable, then $(\tilde{v},\iota_{\tilde v}\theta)\in\Ham_\infty^0\bigl(\Lambda_1^nT^*E,\omega\bigr)$ is reducible.
\end{Theorem}

\begin{proof}
	In light of Lemma \ref{lem:projectable_vfs_are_Hamiltonian_in_multiphase_space}, this follows by similar argument as Proposition \ref{prop:reducible_observables_on_multicotangent_bundles}.
\end{proof}

\subsection{2-dimensional scalar fields}\label{Ssec:2dim-ScalarField}

Let us consider a concrete example of a multiphase space for a field theory with configuration bundle $E=\mathbb R^3\to \Sigma=\mathbb R^2$ possessing two independent variables $(\sigma_1,\sigma_2)$ and one dependent variable $q$. In this setting, a general element of the multicotangent bundle $\Lambda^2T^*E$ has the form
\[
	p\,\d\sigma_1\wedge \d\sigma_2+p^1\,\d\sigma_1\wedge\d q+p^2\,\d\sigma_2\wedge\d q
\]
and we take $\bigl(\sigma_1,\sigma_2,q, p, p^1,p^2\bigr)$ as our preferred coordinates on $\Lambda^2T^*E$. The multiphase space $\Lambda_1^2T^*E$ in this setting is equal to the multicotangent bundle.

The canonical $2$-form on the multiphase space is given by
\[
	\theta=p\,\d\sigma_1\wedge \d\sigma_2+p^1\d\sigma_1\d q+p^2\d\sigma_2\d q
\]
and the multisymplectic form by $\omega=\d\theta$.

Let us now fix an infinitesimal action by $\mathfrak g=\mathbb R$ on $E$, generated by the vector field $v=(q)^2\partial_{q}\in\mathfrak X(E)$. Via natural prolongation, this action induces a multisymplectic action on $M=\Lambda_1^2T^*E$, given by  $V=(q)^2\partial_{q}-2q\bigl(p^1\partial_{p^1}+p^2\partial_{p^2}\bigr)$. By construction, $\theta$ is an invariant potential of $\omega$, hence it induces a covariant moment map.

Let us now apply the reduction procedure. Since $\g$ is one-dimensional and the action on $M$ is generated by $V$, we can set $N=\{x\mid (\iota_V\theta)_x=0\}$. Since $\iota_V\theta=q^2\bigl(p^1\d\sigma_1+ p^2\d\sigma_2\bigr)$, this means that $N=\{q=0\}\cup \bigl\{p^1=p^2=0\bigr\}$ is the intersection of two vector spaces of unequal dimensions. Let us describe the spaces relevant for reduction, where we use $\langle\cdots\rangle$ to denote the $C^\infty(M)$-span of a collection of elements:
\begin{itemize}\itemsep=0pt
	\item $I_N=\bigl\langle qp^1 qp^2\bigr\rangle$,
	\item $\tgtvf=\bigl\langle q\partial_{q}, p^i\partial_{p^j}, \partial_{\sigma_1}, \partial_{\sigma_2},  \partial_p  \bigr\rangle$,
	\item $\vanvf=I_N\cdot\mathfrak X(M)$,
	\item $\tgtvanform^1=I_N\cdot \langle \d{\sigma_1}, \d{\sigma_2}, \d p \rangle + \bigl\langle p^1,p^2\bigr\rangle \d q+ q\cdot \bigl\langle \d p^1,\d p^2\bigr\rangle$.
\end{itemize}

Using this, we can try to determine the reduction.
In degree $0$, the reduction is given by $\frac{\{f|\L_V(f)\in I_N\}}{I_N}$. The condition in the numerator means that
\[
q^2\frac{\partial f}{\partial q}-qp^1\frac{\partial f}{\partial p^1}-qp^2\frac{\partial f}{\partial q} \in I_N.
\] \big(To~avoid possible confusion, note that here \smash{$q^2$} denotes the square of $q$, rather than an index.\big) The second and third term are always in $I_N$, so we can concentrate on the first one. There the condition implies that $f$ has to take the form
\[
	f=f^0\bigl(\sigma_i,p^i,p\bigr)+ qp^1f^1+ qp^2f^2,
\]
where $f^1$, $f^2$ are arbitrary functions and $f^0$ is a function of $\sigma_i$, $p^i$, $p$. Since the second and third term lie in the denominator, the reduction in degree $0$ will be isomorphic to the space of possible functions $f^0$.

Even in this relatively elementary case, calculating the full reduction is very difficult. There is, however, a class of observables for which the calculations are tractable: namely, $1$-forms of the type $\theta(v)$ for $v\in \mathfrak X(E)$. All such forms are observables by the construction of the multiphase space, and their Hamiltonian vector fields are prolongations $\tilde w\in \mathfrak X(M)$ of $w$. A general vector field $w$ on $E$ has the form $w=w^1\partial_{\sigma_1}+w^2\partial_{\sigma_2}+w^q\partial_{q} $
for arbitrary functions $w^1$, $w^2$, $w^q$ of $\sigma_1$, $\sigma_2$, $q$.
Since the multiphase space here is equal to the multicotangent bundle, we need only to check the condition $[v,w ]\in \X_\g(E)$. We obtain
\[
	[v,w]=\bigl[q^2\partial_q,v\bigr]=q^2\frac{\partial w^1}{\partial q}\partial_{\sigma_1}+\frac{\partial w^2}{\partial q}\partial_{\sigma_2}+\biggl(\frac{\partial w^1}{\partial q}-2qw^q\biggr)\partial_{q} .
\]
For this to lie in $\X_\g(E)$, $w^1,w^2$ must be independent of $q$ and $w^q$ must be a multiple of~$q$. A~lengthy calculation shows that an observable of the form $(\iota_{\tilde w}\theta, \tilde w)$ reduces to zero, if and only if $w^1=w^2=0$ and $w^q$ is a multiple of $q^2$. So we have a subspace in the reduction that is isomorphic to
\[
\bigl\{(\iota_{\tilde w}\theta, \tilde w)\mid w=w^1(\sigma_1,\sigma_2)\partial_{\sigma_1}+w^2\partial_{\sigma_2}(\sigma_1,\sigma_2)+q\hat w^q(\sigma_1,\sigma_2)\partial_{q}
\bigr\}.
\]

\section{Outlook}\label{Sec:Outlook}

We briefly consider several natural directions in which to extend our research.

\subsection{Interpretation in terms of the constraint algebraic formalism}

One way to interpret the $L_\infty$ reduction scheme presented in this article, is by the constraint triple formalism introduced
in \cite{DippelEspositoWaldmann2019}. Recasting the above result in the setting of constraint triples is part of the ongoing work \cite{ongoing}, but we will outline the idea here: Heuristically, a constraint triple is a triple of objects intended to describe subquotients, i.e., the quotient of a subspace of some initial space. This is exactly what happens in Marsden--Weinstein--Meyer reduction, in the setting of which a subset $\mu^{-1}(0)\subset M$ is quotiented by the action of a Lie group. A~constraint algebra is a triple $A = (A_t , A_n , A_0)$, where $A_t$ is
an algebra, $A_n$ a subalgebra and $A_0$ an ideal in~$A_n$.  The idea of the constraint triple
approach is working with $A$ instead of $\frac{A_n}{A_0}$ in order to avoid the appearance of singularities.
The constraint formalism allows us to talk about constraint modules, constraint derivations,
constraint tensor products, constraint duals and so~on. With respect to this formalism, the
above reduction procedure of $L_\infty$-algebras naturally appears as the $L_\infty$-algebra associated
to a constraint Lie--Rinehart algebra acting via Cartan calculus on a~constraint differential
graded algebra of differential forms.

Going further, it would be interesting to apply a BFV-BRST-type formalism to resolve singular properties of the $L_\infty$-reduced observables. However, we expect such a procedure to be combinatorially more intricate, since in our case already the unreduced space of observables is graded and homotopical in nature.

\subsection{Moduli spaces of flat connections}

Let $P\to M$ be $G$-principal bundle over an $(n+1)$-dimensional compact manifold $M$ admitting a flat connection, let $\G=\Ad P$ be the gauge group, let $\A$ be the space of connections on $P$, let
\[
	N = \{ A\in\A \mid  F_A = 0\}
\]
be the subspace of flat connections, and let $\M=N/\G$ be the moduli space of flat connections. Given any invariant polynomial $q\in S^{n+1}(\g^*)^G$, it is shown in \cite[Section 10]{CalliesFregierRogersZambon16} that the assignment
\[
	\omega_q(v_1,\dots,v_{n+1}) \mapsto \int_M q(v_1,\dots, v_{n+1}),		\qquad v_1,\dots,v_{n+1}\in\Omega^1(M,\ad P)\cong T_A\A
\]
defines a presymplectic form on $\A$, where we identify $\X(\A)$ with a space of suitably smooth maps $v\colon \A\to\Omega^1(M,\ad P)$ as $\A$ is an $\Omega^1(M,\ad P)$-affine space, and moreover that $\M$ arises as the geometric multisymplectic reduction of $\G\curvearrowright(N\subset\A)$. This extends a classic observation of Atiyah--Bott \cite[Section 9]{AtiyahBott83} in the case in which $M$ is a surface. In each case, the there is an associated moment map derived from the curvature $F\colon \A\to\Omega^2(M,\ad P)$.

It would be interesting to rigorously adapt our framework to this infinite-dimensional setting. In particular, it would be interesting to identify suitable analogues of the auxiliary spaces of Section \ref{Sec:CoreSection} and to examine the reduced $L_\infty$-algebra in particular examples.

\subsection{Classical field theories}

One of the original motivations behind the development of multisymplectic geometry was the search for a rigorous mathematical framework for prequantum field theories.
The driving principle was to look for a suitable extension of the geometric mechanics' framework, based on symplectic geometry (see \cite{AbrahamMarsden78}), from point-like particles with a finite set of configuration coordinates to systems with infinite denumerable degrees of freedom.

Cornerstones of the former are the philosophy of encoding phase spaces as cotangent bundles of configuration manifolds and the fact that cotangent bundles are canonically symplectic.
Extending such ideas to classical field theories led to the introduction of multiphase spaces, as touched on in Section \ref{Ssec:Multiphase}.

At present, the application of recent developments in the theory of multisymplectic observables, moment maps, and reduction to its motivating subject of classical mechanics remains broadly underdeveloped.

A natural next step that could follow from the present work would be to provide a more explicit characterization of the reducible observables in the specific case of multiphase spaces.
In other terms, the latter would imply the displaying a multisymplectic analogue of the so-called \emph{cotangent bundle reduction}, see \cite{MarsdenMisiolekOrtegaPerlmutterRatiu07,OrtegaRatiu06}.
Although it is essentially a matter of restricting general constructions to a peculiar class, essentially comprising fibered spaces with an exact multisymplectic form as sketched in Sections~\ref{Ssec:multicotangent} and~\ref{Ssec:Multiphase}, such results could foster new applications of multisymplectic methods to the realm of the mathematical physics of continuous and field-like systems.

However, it should be clear from the example given in Section \ref{Ssec:2dim-ScalarField} that carrying out all computations required to display all reducible observables explicitly is a daunting task.
Accomplishing this merely for the prototypical model of the free scalar field on the Minkowski space-time would be a nontrivial result, particularly if it were to enable to retrieve geometrically the mechanical momenta usually introduced through variational methods in the physics literature.

\subsection{Quantization}

In the symplectic setting, a \emph{quantization} procedure is an assignment to certain symplectic manifolds $(M,\omega)$, equipped with additional structure, of a Hilbert space $\mathcal{H}(M,\omega)$. When~$(M,\omega)$ comes equipped with a compatible $G$-action, the $\mathcal{H}(M,\omega)$ inherits the structure of a~$G$-repre\-sentation.

Unsurprisingly, multisymplectic setting is more exotic. We refer to \cite{Schreiber21} for a general exposition and present below the prequantization construction of \cite{FiorenzaRogersSchreiber14,FiorenzaRogersSchreiber16}.

\begin{Definition}[\cite{FiorenzaRogersSchreiber14}]
	The cochain complex of sheaves
	\[
		C^\infty(-;U(1)) \xrightarrow{\d\hspace{.5pt}\mathrm{log}} \Omega^1(-) \xrightarrow{\d} \Omega^2(-) \xrightarrow{\d} \cdots \to \Omega^n(-) \to \Omega^{n+1}(-) \to \cdots,
	\]
	with $C^\infty(-;U(1))$ in degree $0$, will be called the \emph{Deligne complex} and will be denoted by the symbol $\underline{U}(1)_{\mathrm{Del}}$.
\end{Definition}

\begin{Definition}[\cite{FiorenzaRogersSchreiber14}]
	The \emph{$n$-stack of principal $U(1)$-$n$-bundles $($or $(n-1)$-bundle gerbes$)$ with connection} $\mathbf{B}^nU(1)_{\mathrm{conn}}$ is the $n$-stack presented via the Dold--Kan construction to the presheaf $\underline{U}(1)_{\mathrm{Del}}^{\leq n}[n]$ regarded as a presheaf of chain complexes concentrated in nonnegative degree.
\end{Definition}

\begin{Definition}[\cite{FiorenzaRogersSchreiber14}]
	Let $(M,\omega)$ be a pre-$n$-plectic manifold. A \emph{prequantization} of $(M,\omega)$ is a~lift
\[
	\begin{tikzpicture}
		\node (A) at (0,0) {$M$};
		\node (B) at (2,0) {$\Omega^{n+1}(-)_{\mathrm{cl}}$.};
		\node (C) at (2,2) {$\mathbf{B}^n U(1)_{\mathrm{conn}}$};

		\draw[->] (A)--node[below]{$\omega$}(B);
		\draw[->] (C)--node[right]{$F$}(B);

		\draw[->] (A) to (B);
		\draw[->,dashed] (A) to node[above left] {$\nabla$} (C);
		\draw[->] (C) to (B);
	\end{tikzpicture}
\]
\end{Definition}

The key observation is that this prequantization construction is given purely in terms of differential forms on $M$, a class of spaces for which we have just proposed a very general reduction scheme (see Remark \ref{rem:reduced_spaces_of_vfs_and_forms}). Consequently, given a subset $N\subset M$ and compatible action $\g\curvearrowright M$, a candidate construction for a \emph{reduced prequantization} that immediately suggests itself is the termwise reduction $\Omega\mapsto\Omega_N$ of the Deligne complex with respect to $\g\curvearrowright M$ along N.

Specializing to the $2$-plectic setting, in \cite{SevestreWurzbacher21} the prequantization of $(M,\omega)$ is realized in the framework of bundle gerbes, and in \cite{Rogers11} in terms of Courant algebroids. In each case, it would be interesting to adapt our reduction apparatus to reduce these prequantizations along a subset $N\subset M$ by an action $\g\curvearrowright(N\subset M)$. A rich field, further quantization schemes are proposed by Barron and Seralejahi \cite{BarronSerajelahi17}, Barron and Shafiee \cite{BarronShafiee19}, Krepski and Vaughan \cite{KrepskiVaughan21}, and by de Bellis, Samann, and Szabo \cite{DeBellisSamannSzabo10,DeBellisSamannSzabo11}.

Several natural questions arise. For example, it would be interesting to identify and investigate a $[Q,R]=0$ ``quantization commutes with reduction'' conjecture. Additionally, it is shown in \cite{FiorenzaRogersSchreiber14} that $\Ham_\infty(M,\omega)$ is suitably equivalent to the infinitesimal symmetries of a prequantization of $(M,\omega)$. It would be interesting to compare the reduction of a prequantization of $(M,\omega)$ with the $L_\infty$-reduction of its space of infinitesimal symmetries $\Ham_\infty(M,\omega)$.

\appendix

\section{Brief survey of symplectic observable reduction}\label{Sec:SymSingRed}

Fix a smooth manifold $M$ and an arbitrary closed subset $N\subset M$.
We often refer to $N$ as a~\emph{constraint set} since, in the context of time-evolving mechanical systems, it can be interpreted as the subspace of the phase space consisting of Cauchy data admissible by the physical constraints.

A (constraints-based) reduction scheme of  $(M,\omega)$  with respect to $N$ is a procedure yielding a~symplectic structure on a certain smooth quotient of $N$ (geometric reduction) or a certain Poisson algebra (observable reduction).
The latter could be eventually interpreted as representing the algebra of ``smooth functions'' over a certain topological space; however, it can happen that the reduced Poisson algebra does not arise from a symplectic structure on the reduced space.
Hence,  while a geometric reduction always implies an observable reduction, the converse is not true in general.

Our main interest will be the reduction with respect to symmetries (i.e., symmetry-based). In this case, we will be concerned with constraint sets induced by Hamiltonian group actions, namely $N=\mu^{-1}(0)$ with $\mu\colon M\to \mathfrak{g}^\ast $ a moment map associated to the action.

When $\mu^{-1}(0)\subset M$ is not a submanifold, or when the action $G\curvearrowright\mu^{-1}(0)$ fails to be free and proper, the reduced space $(M_0,\omega_0)$ prescribed by the Marsden--Weinstein reduction theorem may not exist.
One approach to addressing this pathology is to turn our attention from the underlying symplectic space $(M,\omega)$ to the Poisson algebra of observables $C^\infty(M,\omega)$. As we shall see, there are various natural constructions of a reduced space of observables even in the absence of a reduced space of points. We will designate such constructions \emph{observable reduction} schemes.

This section reviews some well-known approaches to reduction in the symplectic ($1$-plectic) setting, namely the Marsden--Weinstein \cite{MarsdenWeinstein74},  \'Sniatycki--Weinstein \cite{SniatyckiWeinstein83}, Dirac \cite{Dirac64}, Arms--Cush\-man--Gotay \cite{ArmsCushmanGotay91}, and Arms--Gotay--Jennings \cite{ArmsGotayJennings90} reduction schemes.

Throughout this section, $(M,\omega)$ denotes a symplectic manifold and $C^\infty(M,\omega)$ the associated Poisson algebra, that is the (unital, associative, commutative) algebra of smooth func\-tions~$C^\infty(M)$ endowed with the Lie bracket $\{\cdot,\cdot\}$ given by
\[
	\{f_1,f_2\} = \omega( v_{f_1}, v_{f_2})
	= \L_{ v_{f_1}} (f_2)
\]
for any $f_i\in C^\infty(M)$ with associated Hamiltonian vector field $v_{f_i}$.
We note that this bracket is indeed Poisson: in addition to the Jacobi identity it satisfies the Leibniz rule,
\[
	\{h, f_1\,f_2 \} = \{h,f_1\}f_2 + f_1\{h,f_2\}
\]
for all $h,f_1,f_2 \in C^\infty(M)$.

\subsection{Reduction of symplectic manifolds}

Before turning to observable reduction schemes, we first recall the original Marsden--Weinstein reduction theorem for symplectic Hamiltonian actions \cite{MarsdenWeinstein74}.

Let $G\action M$ be the action of a Lie group on a symplectic manifold $(M,\omega)$. We call this a \emph{symplectic action} when $\omega$ is preserved by $G$,
and it is in this setting that the action may additionally admit a moment map.

\begin{Definition}\label{def:momentmap}
	A \emph{moment map} for $G\action(M,\omega)$ is a smooth map $\mu\colon M\to \mathfrak{g}^\ast$ such that
		\begin{enumerate}\itemsep=0pt
			\item[(i)] $\d\mu_\xi= -\iota_\xi\omega$ for all $\xi\in\g$,
			\item[(ii)] $\mu\colon M\to\g^*$ is $G$-equivariant,
		\end{enumerate}
		where $\mu_\xi=\langle\mu,\xi\rangle$ is the contraction on $\g^*\otimes\g$ and where $G\curvearrowright\g^*$ is the coadjoint action.
\end{Definition}

By inverting the order of arguments, the moment map induces a Lie algebra homomorphism
\begin{displaymath}
	\morphism{\tilde\mu}
		{\g}
		{C^\infty(M,\omega),}
		{\xi}{\mu_\xi.}
\end{displaymath}

Thus, every moment map encodes a comoment map, defined as follows.

\begin{Definition}\label{def:comomentmap}
	A \emph{comoment map} for $G\action(M,\omega)$ is a linear map $\widetilde{\mu}\colon  \mathfrak{g}\to C^\infty(M,\omega)$ satisfying
	\begin{enumerate}\itemsep=0pt
		\item[(i)] $\d\,\widetilde{\mu}(\xi) = -\iota_\xi\omega$,
		\item[(ii)] $\widetilde{\mu}([\xi,\zeta]) = \lbrace\widetilde{\mu}(\xi),\widetilde{\mu}(\zeta)\rbrace$
	\end{enumerate}
	for all $\xi,\zeta\in\g$.
\end{Definition}

\begin{Remark}
	Condition (ii) of Definition \ref{def:momentmap} expresses the equivariance of $\mu$ with respect to $G\action M$ and the coadjoint action $G\action \mathfrak g^*$. This property implies that $\widetilde\mu$ is a Lie algebra homomorphism (condition (ii) of Definition \ref{def:comomentmap}).
	The converse is true when $G$ is a connected group.\looseness=-1
\end{Remark}

\begin{Remark}[comoment maps as lifts]\label{rem:comomentAsLift}
	An action admitting a moment map acts infinitesimally by Hamiltonian vector fields. The comoment map assigns to each $\xi\in\g$ a Hamiltonian function $\mu_\xi\in C^\infty(M)$ associated to the fundamental field $\underline\xi$. In this setting, $G\curvearrowright M$ is said to be a~\emph{Hamiltonian action}.

	More algebraically, $\widetilde\mu$ is a lift in the category of Lie algebras of the fundamental action $\xi\mapsto\underline\xi$ by the assignment of Hamiltonian vector fields $f\mapsto v_f$,
\[
	\begin{tikzpicture}
		\node (A) at (0,0) {$\g$};
		\node (B) at (2,0) {$\X(M)$};
		\node (C) at (2,2) {$C^\infty(M,\omega)$};

		\node (D) at (0,-.6) {$\xi$};
		\node (E) at (2,-.6) {$\underline\xi$};

		\node (F) at (3.15,2) {$\alpha$};
		\node (G) at (3.15,0) {$v_\alpha$};

		\draw[->] (A) to (B);
		\draw[->,dashed] (A) to node[above left] {$\widetilde\mu$} (C);
		\draw[->] (C) to (B);

		\draw[|->] (D) to (E);
		\draw[|->] (F) to (G);
	\end{tikzpicture}
\]
\end{Remark}

The Marsden--Weinstein symplectic reduction scheme is a rule that associates to each suitably compatible{\samepage
\begin{enumerate}\itemsep=0pt
	\item[(i)] symplectic Hamiltonian action $G\curvearrowright(M,\omega)$,
	\item[(ii)] moment map $\mu\colon M\to\g^*$, and
	\item[(iii)] element $\lambda\in\g^*$,
\end{enumerate}
a \emph{reduced} symplectic manifold $(M_\lambda,\omega_\lambda)$.}

\begin{Theorem}[Marsden--Weinstein symplectic reduction \cite{MarsdenWeinstein74}, see also \cite{Meyer73}]\label{thm:MW-red}
 Consider a symplectic action on the symplectic manifold $(M,\omega)$, with moment map  $\mu\colon M\to\g^\ast$.
 Let $\lambda\in\g^*$ be a regular value of $\mu$, denote by $N=\mu^{-1}(\lambda)$ the corresponding smoothly embedded level set.
 Assume also the action of the isotropy subgroup $G_\lambda$ restricted to $N$ to be free and proper, denote as $M_\lambda= \mu^{-1}(\lambda)/G_\lambda$ the corresponding quotient manifold.
 	Then there is a unique symplectic form $\omega_\lambda\in\Omega^{2}(M_\lambda)$ satisfying $j^*\omega = \pi^*\omega_\lambda$, where $j\colon \mu^{-1}(\lambda)\to M$ is the embedding of $N$ in $M$ and $\pi\colon \mu^{-1}(\lambda)\to M_\lambda$ is the canonical quotient.
\end{Theorem}

\begin{Definition}
	The symplectic manifold $(M_\lambda,\omega_\lambda)$ provided by Theorem \ref{thm:MW-red} is called the \emph{Marsden--Weinstein symplectic reduction} of $(M,\omega)$ with respect to the Hamiltonian action $G\action (M,\omega)$.
\end{Definition}

\begin{Remark}
In rough strokes, the Marsden--Weinstein reduction is a two-step process:
\begin{enumerate}\itemsep=0pt
	\item[(i)] first, we \emph{restrict} to a constraint set $N\subset M$,
	\item[(ii)] then, we \emph{descend} to the quotient $N/G$.
\end{enumerate}
At each stage, we rely on the assumption that $N\subset M$ is a smoothly embedded submanifold.

Several alternative reduction procedures have been introduced to account for the less well-behaved situation where $N$ is not smooth but is, for example, the preimage of a singular value of a moment map (hence the name \emph{singular reductions}).
\end{Remark}

\subsection{{[\'SW]} reduction}

The \'Sniatycki--Weinstein reduction is a symmetry-based observable reduction scheme defined for any singular constraint set \cite{SniatyckiWeinstein83}.

Such a procedure is guaranteed to produce a reduced Poisson algebra out of the regular Poisson structure naturally associated with the symplectic manifold on which the group acts.
Furthermore, in the case in which a geometric reduction does exist, such reduced Poisson algebra coincides with the canonical Poisson structure associated with the geometrically reduced symplectic manifold.

Let $(M,\omega)$ be a symplectic manifold and $G$ a connected Lie group with Lie algebra $\mathfrak{g}$ and dual $\mathfrak{g}^*$.
Consider a Hamiltonian action $G\action (M,\omega)$ with moment map $\mu\colon M\to \mathfrak{g}^*$.

The [\'SW] reduction procedure is based on the ideal generated by the momenta $\mu_\xi$ introduced in definition \ref{def:momentumideal}.
	We call $N= \mu^{-1}(0)$ the \emph{constraints locus} of the moment map $\mu$.
	This is a~subset of $M$, in general not a smooth submanifold, defined as the \emph{zero-level set} of all possible momenta (regarded as constraints).
	In general, $I_\mu$ is included in the ideal $I_N\subset C^\infty(M)$ of smooth functions vanishing along $N$, see Lemma \ref{lem:ImuVSIN} below.

\begin{Lemma}\label{lem:ImuPoisson}
	$I_\mu\subset C^\infty(M,\omega)$ is a $G$-stable Poisson subalgebra.
\end{Lemma}

\begin{proof}
	Every element of $I_{\mu}$ is a linear combination of products $f \mu_\xi$ for $f \in C^\infty(M)$ and $\xi \in \mathfrak{g}$. Two applications of the Leibniz rule yield
	\[
		\{f \mu_\xi, h \mu_\zeta\} = fh\, \{ \mu_\xi,\mu_\zeta \} + q
	\]
	for some $q\in I_{\mu}$.
	Since the equivariance of $\mu$ is equivalent to the condition $\{\mu_\xi,\mu_\zeta\}=\mu_{[\xi,\zeta]}$, see, e.g., \cite[Section~22.1]{CannasdaSilva01}, we conclude that the Poisson bracket lies in $I_{\mu}$.

	For $x\in M$, we have
	\[
		( g \cdot \mu_\xi )(x) =
		\mu_\xi\bigl(g^{-1}\cdot x\bigr)
		= \mu_{g^{-1}\cdot\xi}(x),
	\]
	where the second equality follows from the equivariance of $\mu$.
	Therefore, $G$ preserves $I_{\mu}$.
\end{proof}

\begin{Remark}
	Observe that the $G$-invariant Poisson subalgebra $I_\mu$ need not be a Lie algebra ideal with respect to the Poisson bracket.
	On the other hand, if we restrict to the Poisson subalgebra of $G$-invariant functions $C^\infty(M)^G$, it follows that $I_{\mu}^G = I_\mu \cap C^\infty(M)^G$ is a Poisson ideal since
	\[
		\bigl\lbrace	\mu_{\xi_i} , h \bigr\rbrace = \L_{\xi_i} h = 0
		\qquad \forall h \in C^\infty(M)^G.
	\]
\end{Remark}

From Lemma \ref{lem:ImuPoisson}, it follows that the action of $G$ on $C^\infty(M)$ induces an action on the quotient algebra $C^\infty(M)/I_\mu$ such that the projection algebra morphism $\rho\colon  C^\infty(M) \twoheadrightarrow C^\infty(M)/I_{\mu}$ is $G$-equivariant.

\begin{Theorem}[\'Sniatycki--Weinstein reduction \cite{SniatyckiWeinstein83}]\label{thm:SW-red}
 Consider a Hamiltonian action on the symplectic manifold $(M,\omega)$ and let $\mu\colon M\to\g^*$ be the corresponding moment map.
 Denote by $N=\mu^{-1}(0)$ the constraint set given as the zero locus of the moment map.
 The space of $G$-invariant elements of the quotient algebra $C^\infty(M)/I_\mu$, denoted as \smash{$\bigl(C^\infty(M)/I_\mu\bigr)^G$},
 forms a Poisson algebra together with the binary operator obtained by pushing forward the Poisson structure of $C^\infty(M,\omega)$ by the canonical projection $\rho\colon  C^\infty(M,\omega) \to C^\infty(M)/I_\mu$.
	Furthermore, we have an isomorphism of Poisson algebras\footnote{Observe that the quotient on the left-hand side is meant in the sense of associative algebras while the one on the right-hand side is a quotient in the category of Lie algebras.}
	\begin{displaymath}
		\biggl(\frac{C^\infty(M)}{I_{\mu}}\biggr)^G \cong
		\frac{\mathcal{N}\bigl(I_\mu\bigr)}{I_\mu},
	\end{displaymath}
	where $\mathcal{N}\bigl(I_\mu\bigr)$ is the Lie algebra normalizer of $I_\mu$ in $C^\infty(M,\omega)$, that is,
	\begin{equation}\label{eq:ImuNorm}
		\mathcal{N}\bigl(I_\mu\bigr) =
		\bigl\lbrace
 		f \in C^\infty(M)
 		\mid
		\bigl\lbrace f, I_\mu \bigr\rbrace \subset I_\mu
 		\bigr\rbrace 		.
	\end{equation}
\end{Theorem}

\begin{Definition}
	The Poisson algebra \smash{$\bigl(C^\infty(M)/I_\mu\bigr)^G$} of Theorem \ref{thm:SW-red} is called the \emph{\'Sniatycki--Weinstein reduction} of $C^\infty(M,\omega)$ with respect to the Hamiltonian action $G\action (M,\omega)$.
\end{Definition}

\begin{Remark}[observable reduction at nonzero orbits]
In the case of the [MW] reduction, the \emph{shifting trick} of Guillemin and Sternberg \cite{GuilleminSternberg90} establishes the equivalence of the geometric reduction at any $\lambda\in \mathfrak{g}^*$ with the reduction at $0$ for a suitably modified symplectic Hamiltonian action constructed out of the coadjoint orbit of $\lambda$.
A similar equivalence has been extended in~\cite{Arms96} to the case of not necessarily free and proper actions.
	We thus restrict our attention to the observable reduction at $0\in \mathfrak{g}$ without loss of meaningful generality.
\end{Remark}

\subsection{{[D]} reduction}

The Dirac reduction is a constraints-based observable reduction scheme defined on singular constraint sets satisfying a so-called \emph{first class} condition.

Recall that we denote by $F_N$ the set of \emph{first class function} (see Definition \ref{def:IclassFunc})
Elements of $F_N \cap I_N$ are called \emph{first class constraints}; those constraints that not are first class are said to be \emph{second class}.

\begin{Definition}[first class constraint set]
	A closed subset $N\subset M$ is said to be a \emph{first class set} if every associated constraint is first class, i.e., if $I_N\subset F_N$.
\end{Definition}

According to the Dirac's theory of constraints \cite{Dirac64} (see  also \cite{Sniatycki83}), it is useful to consider a~certain subclass of well behaved functions.
\begin{Definition}[Dirac observable]\label{def:Dirac_observables}
	We call a function $f\in C^\infty(M)$ a \emph{Dirac observable} if its Poisson bracket with any first class constraint vanishes on $N$. Formally, the set of Dirac observables is
	\begin{displaymath}
		\mathcal{O}(N) =
		\lbrace
 		f \in C^\infty(M)
 		\mid
		\lbrace f, F_N\cap I_N \rbrace \subset I_N
		\rbrace.
	\end{displaymath}
\end{Definition}

\begin{Remark}
	In general, $F_N \subset \mathcal{O}(N)$ and $\mathcal{O}(N)$ is not a Lie subalgebra of $C^\infty(M)$.
	When $N$ is first class, $\mathcal{O}(N)= F_N$, in particular Dirac observables form a Poisson subalgebra of $C^\infty(M)$.
\end{Remark}

\begin{Theorem}[Dirac reduction, {\cite{Dirac64} and \cite[Proposition~3.1]{ArmsGotayJennings90}}]\label{thm:Dirac}
	Consider a first class constraint set $N$.
	The quotient associative algebra $\mathcal{O}(N)/I_N$ forms a Poisson algebra together with the binary operator
	\begin{displaymath}
		\lbrace [h], [k] \rbrace = [\lbrace h, k \rbrace] \qquad \forall h,k \in \mathcal{O}(N),
	\end{displaymath}
	obtained by pushing-forward the Poisson structure of $\mathcal{O}(N)$ along the canonical projection
\[
[\,\cdot\,]\colon \ \mathcal{O}(N) \to \mathcal{O}(N)/I_N.
\]
\end{Theorem}

\begin{Definition}
	The Poisson algebra $\mathcal{O}(N)/I_N$ obtained by Theorem \ref{thm:Dirac} is called the \emph{Dirac reduction} of $C^\infty(M,\omega)$ with respect to $N$.
\end{Definition}	

Such a reduction can be interpreted geometrically regarding the above Poisson algebra as the set of ``smooth'' functions on the reduced topological space $N/{\sim}$, where $q\sim p$ if and only if $h(q)=h(p)$ for all Dirac observables $h$.
That means that states $q,p\in N$ are identified whenever they cannot be distinguished by means of the measurable quantities of $\mathcal{O}(N)$.

Assume that $N$ is the zero locus of a moment map $\mu$ associated to a certain group action $G\action M$.
We can consider at the same time the associative ideal generated by momenta $I_{\mu}=\langle \mu_\xi \rangle_{\xi \in \mathfrak{g}}^{\text{asso}}$ and the ideal of vanishing functions on $N$.
In general $I_{\mu}\neq I_{N}$, more precisely we have the following chain of inclusions.

\begin{Lemma}\label{lem:ImuVSIN}
	Let $\mu\colon M\to\g^*$ be a moment map for the symplectic action $G\action (M,\omega)$. Denote by $N=\mu^{-1}(0)$ the zero locus of the moment map and by $I_\mu$ and $I_N$ the ideals of Definitions~$\ref{def:momentumideal}$ and~$\ref{def:constraintideal}$ respectively.
	We have the following diagram in the category of vector spaces
	\[
	\begin{tikzcd}[column sep = 2.8em]
		& & &
		I_N \ar[dr,hook,dashed,"(asso.)","ideal"',sloped]
		&
		\\
		C^\infty(M)^G\cap I_\mu \ar[r,hook] \ar[d,hook] &
		I_{\mu} \ar[r,hook] \ar[d,hook,"(Pois.)"',"ideal"]
		& F_N\cap I_N \ar[ur,sloped, hook, dashed] \ar[dr,sloped, hook,"(Pois.)","ideal"',sloped]
		& &
		\mathcal{O}(N) \ar[dd,hook,dashed]
		\\[.5em]
		C^\infty(M)^G \ar[r,hook] &
		\mathcal{N}\bigl(I_\mu\bigr) \ar[dr,hook,sloped]& &
		\mathcal{N}(I_N)=F_N \ar[ur,sloped,hook,dashed]\ar[dl,hook,sloped] &
		\\
		& &
		C^\infty(M,\omega) \ar[rr,equal]& &
		C^\infty(M),
	\end{tikzcd}	
	\]
	where solid $($resp.\ dashed$)$ arrows denote Poisson $($resp.\ associative$)$ algebra morphisms.
\end{Lemma}

\begin{proof}
	According to Lemma \ref{lem:ImuPoisson}, $I_\mu\subset C^\infty(M,\omega)$ is a Poisson subalgebra.
	By the definition of the normalizer, $\mathcal{N}\bigl(I_\mu\bigr)$ is the largest Lie subalgebra of $C^\infty(M,\omega)$ such that $I_\mu\subset\mathcal{N}\bigl(I_\mu\bigr)$ is a~Lie algebra ideal.
	The Jacobi identity yields $\{F_N,F_N\}\subset F_N$, from which it follows that $F_N$ is a Poisson subalgebra.
	Similarly, $F_N\cap I_N$ is a Lie ideal in $F_N= \mathcal{N}(I_N)$.
	The assumption $N=\mu^{-1}(0)$ provides that $I_\mu$ is a Poisson subalgebra of first order constraints.
	By construction, $I_{N}$ is an associative ideal and both $F_N$ and $I_N$ lie in $\mathcal{O}(N)$.
	The proof that $I_{\mu}\subset F_N$ is an inclusion of Poisson algebras can be found in \cite[Proposition~5.1]{ArmsGotayJennings90}.
	The property of $C^\infty(M)^G$ to be a Lie subalgebra follows from the Jacobi identity for the Poisson bracket of $C^\infty(M)$.
	The commutation of the uppermost square is the trivial pullback given by the intersection operation in the category of vector spaces.
\end{proof}

We remark that some of the above inclusions are generally strict, for example, it is possible that $I_{\mu}\neq F_N\cap I_N \neq I_{N}$ (see \cite[Section~7]{ArmsGotayJennings90}).

\begin{Proposition}[{[\'SW] vs.\ [D]}]
	If the hypotheses of Theorems {\rm \ref{thm:SW-red}} and {\rm \ref{thm:Dirac}} hold, and if $I_\mu= I_N$, then the \'Sniatycki--Weinstein and Dirac reduced Poisson algebras coincide.
\end{Proposition}

\begin{proof}
	Consider a Hamiltonian group action $G\action (M,\omega)$ with moment map $\mu$.
	When $N=\mu^{-1}(0)$ is first class, the diagram of Lemma \ref{lem:ImuVSIN} condenses to the following open square in the category of Poisson algebras:
	\begin{displaymath}
		\begin{tikzcd}
			I_\mu \ar[r,hook] \ar[d,hook]&
			\mathcal{N}\bigl(I_\mu\bigr) \ar[r] &
			\dfrac{\mathcal{N}\bigl(I_\mu\bigr)}{I_\mu}
			\\
			I_N\ar[r,hook] &
			F_N = \mathcal{O}(N) \ar[r] &
			\dfrac{\mathcal{O}(N)}{I_N}.
		\end{tikzcd}
	\end{displaymath}
	In light of the equality $F_N = \mathcal{N}(I_N)$, the rows above are identical when $I_\mu=I_N$.
\end{proof}

\begin{Remark}[on the technical condition $I_\mu = I_N $]
	As discussed in \cite[Theorem~1]{SniatyckiWeinstein83} and \cite[Corollary~6.2]{ArmsGotayJennings90}, when $0$ is a weakly regular point, the momentum ideal and the constraint ideal with respect to $N=\mu^{-1}(0)$ coincide.
	In particular, when $N\subset M$ is a smoothly embedded submanifold we have that $I_\mu=I_N$.
	
	Understanding the relationship between $I_N$ and $I_\mu$ in general is a nontrivial problem in $C^\infty$ algebraic geometry. The special case of compact group actions is discussed in \cite[Sections~5 and~6]{ArmsGotayJennings90}. The case of free actions on paracompact manifolds is discussed in \cite[Lemma~2]{Sniatycki05}.
\end{Remark}

\subsection{{[ACG]} reduction}\label{Ssec:ACG-red}

The Arms--Cushman--Gotay reduction is another symmetry-based observable reduction scheme defined on singular constraint sets \cite{ArmsCushmanGotay91}.
For any subset $S\subset C^\infty(M)$, we denote by $S^G$ the subspace of $G$-fixed elements.

\begin{Definition}\label{def:ACGreduction}
	The Poisson algebra $C^\infty(M)^G / I_N^G$ is called the \emph{Arms--Cushman--Gotay reduction} of $C^\infty(M,\omega)$ with respect to the Hamiltonian action $G\action (M,\omega)$.
\end{Definition}

\begin{Remark}
	Despite introducing the ACG reduction exclusively in terms of observables, this scheme also admits a suitable interpretation as a geometric reduction of the symplectic space.
Namely, the ACG reduction is isomorphic to the unique Poisson structure  induced on the space of smooth functions on the variety $N/G$, interpreting $C^\infty(N/G)$ as in Remark \ref{rem:generalizedSmootheology}, from $C^\infty(M)$ by the following commutative diagram of suitably smooth mappings:
	\begin{displaymath}
	\begin{tikzcd}
		N \ar[r]\ar[d] & M \ar[d]\\
		N/G \ar[r]& M/G,
	\end{tikzcd}
	\end{displaymath}
see \cite[Theorem~1]{ArmsCushmanGotay91} for further details.
In the case of proper action, the latter Poisson structure is proved to be non degenerate \cite[Theorem~2]{ArmsCushmanGotay91}.
Hence, according to this interpretation, the ACG reduction scheme has been introduced as the \emph{universal symplectic reduction} of the symplectic manifold $(M,\omega)$.
\end{Remark}

\begin{Proposition}[{[ACG] vs.\ [\'SW], \cite[Theorem~6.1]{ArmsGotayJennings90}}]\label{prop:ACGvsSW}
	Consider an Hamiltonian group action $G\action (M,\omega)$.
	When $I_\mu^G \cong I_N^G$, there is a natural embedding of the ACG reduced Poisson algebra into the SW reduction.
	When $\mathcal{N}\bigl(I_\mu\bigr) \cong C^\infty(M)^G + I_\mu$, there is a natural surjection from the~\'SW reduced Poisson algebra to the ACG reduction.
	If both conditions applies, the two reduced Poisson algebras are isomorphic.
\end{Proposition}

\begin{proof}
When we consider a Hamiltonian group action $G\action (M,\omega)$ and a constraint set $N=\mu^{-1}(0)$, i.e., when Theorem~\ref{thm:SW-red} applies, we have the following commutative diagram in the category of Poisson algebras, where the dashed arrows indicate morphisms of associative algebras only,
\begin{displaymath}
	\begin{tikzcd}[column sep = large, 	
					cells={font=\everymath\expandafter{\the\everymath\displaystyle}},]
		I_{\mu} \ar[r,hook,"\text{\tiny (asso.) ideal}"']
		&[1em] C^{\infty}(M) \ar[r,two heads,dashed,"\rho"]
		& \frac{C^\infty(M)}{I_{\mu}}
		\\
		I_{\mu} \ar[u,equal] \ar[r,hook,"\text{\tiny (Lie) ideal}"']
		& \mathcal{N}\bigl(I_\mu\bigr) \ar[r,two heads,"\rho"] \ar[u,hook]
		& {\left(\frac{C^\infty(M)}{I_{\mu}} \right)^G} \ar[u,hook,dashed] \ar[r,equal,"\sim"]
		&[-2em] \frac{\mathcal{N}(I_{\mu})}{I_{\mu}}
		\\
		{I_{\mu}^G} \ar[u,equal]\ar[d,hook] \ar[r,hook,"\text{\tiny (Lie) ideal}"']
		\ar[ur,phantom,very near start, "\urcorner"]
		& {C^{\infty}(M)^G} \ar[r,two heads,"\rho"] \ar[u,hook]
		& \frac{C^\infty(M)^G}{I_{\mu}^G} \ar[r,equal,"\sim"] \ar[d,two heads]
		& \frac{C^\infty(M)^G + I_\mu}{I_{\mu}} \ar[u,hook]
		\\
		{I_N^G}  \ar[r,hook,"\text{\tiny (Lie) ideal}"']
		& {C^{\infty}(M)^G} \ar[r,two heads,"\rho"] \ar[u,equal]
		& \frac{C^\infty(M)^G}{I_{N}^G}.
		&
	\end{tikzcd}
\end{displaymath}
The top two squares encode the definition of the [\'SW] reduced Poisson algebra and the observation that $I_{\mu}$ is a normal Lie subalgebra of $\rho^{-1}\bigl[(\sfrac{C^\infty(M)}{I_{\mu}})^G\bigr]$, where $\rho$ is the canonical projection on the quotient (see \cite[Lemma~2]{SniatyckiWeinstein83} or Theorem \ref{thm:SW-red} above). The preimage
\begin{displaymath}
	\mathcal{N}\bigl(I_{\mu}\bigr)
	=
	\bigl\lbrace
		f \in C^\infty(M)	
	\mid
		\bigl\lbrace f, I_{\mu} \bigr\rbrace \subset I_{\mu}
	\bigr\rbrace
\end{displaymath}
is the Lie algebra normalizer of $I_\mu$ in $C^\infty(M,\omega)$.
We denoted by $I_\mu^G$ the intersection of the momentum ideal with the vector space of $G$-invariant smooth functions.
Observe that $C^\infty(M)^G \subset \mathcal{N}\bigl(I_\mu\bigr)$.
The middle two squares encode the second isomorphism theorem for Lie algebras.
Finally, the bottom two squares express the relation between the two quotients computed with respect to an ideal $I$ and a certain smaller ideal $I'\subset I$.

Overall, we end up with the following morphisms of Poisson algebras:
\[
	\text{[\'SW]}  \hookleftarrow \frac{C^\infty(M)^G}{I_\mu^G} \twoheadrightarrow \text{[ACG]},
\]
where the left (resp.\ right) mapping is an isomorphism when $\mathcal{N}\bigl(I_\mu\bigr)\cong C^\infty(M)^G + I_\mu$ (resp.\ $I_\mu^G=I_N^G$).
\end{proof}

\begin{Remark}
 When the acting group $G$ is compact, we can make use of the averaging trick to conclude that \smash{$\bigl(C^\infty(M)/ I_\mu\bigr)^G \cong C^\infty(M)^G / I_\mu^G$} (see \cite[Proposition~5.12]{ArmsGotayJennings90} or \cite[Proposition~5]{Sniatycki05} for details) hence $\mathcal{N}\bigl(I_\mu\bigr)= C^\infty(M)^G + I_\mu$ and Proposition \ref{prop:ACGvsSW} applies.
 In the case of compact group action the ACG reduction coincide with the Arms--Gotay--Jennings reduction \cite[Theorem~6.1]{ArmsGotayJennings90}
 with respect to the zero-level set of a moment map.
 The latter, in its more general incarnation, is a purely constraints-based reduction scheme defined on singular sets that satisfy mild technical conditions. For example, the constraint set $N$ must be strongly isotropic and locally conical (see \cite[Section~3]{ArmsGotayJennings90} for further details).
\end{Remark}

\begin{Remark}
In \cite{Sniatycki05} can be found a more general account of symplectic observable reduction schemes by framing them as two-steps procedures akin to the Marsden--Weinstein theorem.
In particular, taking $N=\mu^{-1}(\lambda)$ for a possibly nonzero $\lambda \in \g^*$, and denoting by $I_\mu$ the associative ideal generated by $\{\mu_\xi - \lambda\}_{\xi \in \g}$, we obtain the commutative square
\begin{displaymath}
	\begin{tikzcd}
		\dfrac{C^\infty(M)^G}{I_\mu^G} \ar[r]\ar[d] &
		\left(\dfrac{C^\infty(M)}{I_\mu}\right)^{G_\lambda} \ar[d,dashed]
		\\
		\dfrac{C^\infty(M)^G}{I_N^G} \ar[dashed,r] &
		\left(\dfrac{C^\infty(M)}{I_N}\right)^{G_\lambda},
	\end{tikzcd}
\end{displaymath}
where the solid (resp.\ dashed) arrows denote Poisson (resp.\ associative) morphisms, the top-right node encodes the SW reduction, and the bottom-left node gives the ACG reduction (see \cite[Remarks~1 to~3 and equation~(11)]{Sniatycki05} for complete details).
\end{Remark}

\subsection*{Acknowledgements}

C.B.\ would like to acknowledge the support of the Leonhard Euler International Mathematical Institute in Saint Petersburg, the Saint Petersburg State University, and the Ministry of Science and Higher Education of the Russian Federation agreement no.~075-15-2022-287.
A.M.M.\ thanks the Max Planck Institute for Mathematics in Bonn for its hospitality and financial support.
This work has received funding from the European Union's \emph{Horizon 2020 research and innovation programme} under the grant agreement no.~101034324 and has been partially supported by the \emph{Italian Group for Algebraic and Geometric Structures and their Application} (GNSAGA–INdAM).
L.R.~is supported by the CNRS project GraNum and by the RTG2491.
The authors thank Janina Bernardy, Christian Blohmann, Luca Vitagliano, and Marco Zambon for helpful and motivating conversations, and Tatyana Barron for indicating that the premultisymplectic form $\omega$ need only be $\g$-invariant, and not reducible, for the reduced space of observables to be defined. The authors would also like to thank the anonymous referees for their helpful and insightful comments.


\pdfbookmark[1]{References}{ref}
\LastPageEnding

\end{document}